\documentclass[journal]{IEEEtran}
\usepackage{amsmath,amssymb,amsfonts}
\usepackage{cases}
\usepackage{color,xcolor}
\usepackage{graphicx}
\usepackage{cite}
\usepackage{subfigure}
\usepackage{algorithm}
\usepackage{algorithmic}
\usepackage{epstopdf}
\usepackage{multirow}
\usepackage{rotating}
\usepackage{booktabs}
\usepackage{tikz}
\usepackage{acronym}
\usepackage{epstopdf}
\usepackage{pifont} 

\usepackage[demo,
            export]{adjustbox}  
\usepackage{tabularray}
\usepackage{colortbl}
\usepackage[colorlinks,
linkcolor=blue,
anchorcolor=blue,
citecolor=blue
]{hyperref}
\usepackage{amssymb}
\usepackage{array} 
\usepackage{catchfilebetweentags}
\definecolor{tableTitleColor}{rgb}{0.8,0.8,0.8}
\newtheorem{Theorem}{Theorem}[section]
\newtheorem{Lemma}[Theorem]{Lemma}

\newtheorem{Remark}[Theorem]{Remark}

\begin{document}

\title{Bi-Sparse Unsupervised Feature Selection}

\author{Xianchao Xiu, \IEEEmembership{Member,~IEEE}, Chenyi Huang, Pan Shang, and Wanquan Liu, \IEEEmembership{Senior Member,~IEEE}

\thanks{This work was supported in part by the National Natural Science Foundation of China under Grants 12371306, 12401430, 12271309, and the Postdoctoral Fellowship Program of CPSF under Grant GZB20240801. (\textit{Corresponding author: Pan Shang}.)}
\thanks{Xianchao Xiu and Chenyi Huang are with the School of Mechatronic Engineering and Automation, Shanghai University, Shanghai 200444, China (e-mail: xcxiu@shu.edu.cn; huangchenyi@shu.edu.cn).}
\thanks{Pan Shang is with the School of Mathematics and Statistics, Beijing Jiaotong University, Beijing 100044, China, and the Academy of Mathematics and Systems Science, Chinese Academy of Sciences, Beijing 100190, China (e-mail: pshang@amss.ac.cn).}
\thanks{Wanquan Liu is with the School of Intelligent Systems Engineering, Sun Yat-sen University, Guangzhou 510275, China (e-mail: liuwq63@mail.sysu.edu.cn).}
}

\maketitle

\begin{abstract}
To deal with high-dimensional unlabeled datasets in many areas, principal component analysis (PCA) has become a rising technique for unsupervised feature selection (UFS). However, most existing PCA-based methods only consider the structure of datasets by embedding a single sparse regularization or constraint on the transformation matrix. In this paper, we introduce a novel bi-sparse method called BSUFS to improve the performance of UFS. The core idea of BSUFS is to incorporate $\ell_{2,p}$-norm and $\ell_q$-norm into the classical PCA, which enables our method to select relevant features and filter out irrelevant noises, thereby obtaining discriminative features. Here, the parameters $p$ and $q$ are within the range of $[0, 1)$. Therefore, BSUFS not only constructs a unified framework for bi-sparse optimization, but also includes some existing works as special cases. To solve the resulting non-convex model, we propose an efficient proximal alternating minimization (PAM) algorithm using Stiefel manifold optimization and sparse optimization techniques. In addition, the computational complexity analysis is presented. Extensive numerical experiments on synthetic and real-world datasets demonstrate the effectiveness of our proposed BSUFS. The results reveal the advantages of bi-sparse optimization in feature selection and show its potential for other fields in image processing. Our code is available at \href{https://github.com/xianchaoxiu/BSUFS}{https://github.com/xianchaoxiu/BSUFS}.
\end{abstract}


\section{Introduction}

\IEEEPARstart{H}{igh-dimensional} features make data processing challenging  due to the presence of redundant information and additional noises \cite{jordan2015machine}. To address this challenge, feature selection techniques are proposed and widely studied. One can read the review paper \cite{zebari2020comprehensive} for more comprehensive backgrounds and advances. Among feature selection techniques, there is a special type called unsupervised feature selection (UFS), which is designed to select features for unlabeled datasets, i.e., label information is missing. Note that unlabeled datasets are cheap and easy to obtain in practice. Therefore, a large number of researchers have devoted themselves to developing UFS methods, which have shown outstanding performances in image processing \cite{bolon2020feature,li2015unsupervised,hu2024bi}, gene analysis \cite{ang2015supervised,garcia2020unsupervised,rahman2020review},  machine learning \cite{cai2018feature,dy2004feature,guo2024unsupervised}, and deep learning \cite{solorio2020review,10410215,10014652}.

According to \cite{zebari2020comprehensive,bolon2020feature}, existing UFS methods can be generally divided into three categories based on different search strategies, including filtering methods, wrapper methods, and embedded methods. It is worth pointing out that the embedded methods combine the advantages of filtering methods and wrapper methods by directly integrating feature selection into the learning procedure \cite{li2017feature}. Under some graph techniques, there are many representative UFS methods that are  proposed  to  optimize the similarity matrix,  such as Laplacian score (LapScore) \cite{he2005laplacian}, multi-cluster feature selection (MCFS) \cite{cai2010unsupervised}, unsupervised discriminative feature selection (UDFS) \cite{yang2011}, structured optimal graph feature selection (SOGFS)\cite{nie2016unsupervised}, robust neighborhood embedding (RNE)\cite{liu2020robust}, and non-convex regularized graph embedding and self-representation (NLGMS) \cite{bai2024precise}, to name a few.

As one of the most straightforward and accessible embedded methods, principal component analysis (PCA) stands out from UFS methods \cite{greenacre2022principal}. By introducing a linear transformation matrix (also projection matrix), PCA tends to characterize this matrix and is different from graph-based embedded methods. Despite its excellent performance, one of the main drawbacks is the lack of feature interpretability \cite{zou2018selective}. To overcome this limitation, Li et al. \cite{li2023sparse} combined PCA with the $\ell_{2,p}$-norm  regularization term on the transformation matrix. This novel PCA variant, referred to as SPCAFS, aims to achieve feature sparsity and still retain the principal information simultaneously. Here, $p\in(0,1)$ and thus SPCAFS is a non-convex version of sparse PCA with $\ell_{2,1}$-norm \cite{tian2022comprehensive}. In addition, Nie et al. \cite{nie2022learning} considered feature-sparsity constrained PCA (FSPCA), which is accomplished by enforcing the $\ell_{2,0}$-norm constraint. Note that $\ell_{2,0}$-norm is a direct choice to characterize the feature sparsity \cite{Zhang2022structured}.  FSPCA enables selecting a subset of features that are most representative and intrinsic. Recently, Zheng et al. \cite{zheng2023fast} proposed a sparse PCA method based on positive semidefinite projection (SPCA-PSD), which achieves excellent efficiency in clustering. Gao et al. \cite{gao2024principal} extended the $\ell_{2,p}$-norm to fuzzy elastic net for better characterizing the structure of the transformation matrix, which is called FEN-PCAFS. 
The above $\ell_{2,1}$-norm, $\ell_{2,p}$-norm, and $\ell_{2,0}$-norm can be used to capture row sparsity, which is closely related to the selected features. These sparse regularization terms attribute interpretability and robustness to PCA models, which is  beneficial for understanding the data structure \cite{hu2017group}.

\begin{figure*}[htbp]	
\hspace{-0.5cm}
	\subfigure[SPCAFS] 
	{
		\begin{minipage}{7cm}
			\centering          
			\includegraphics[scale=0.43]{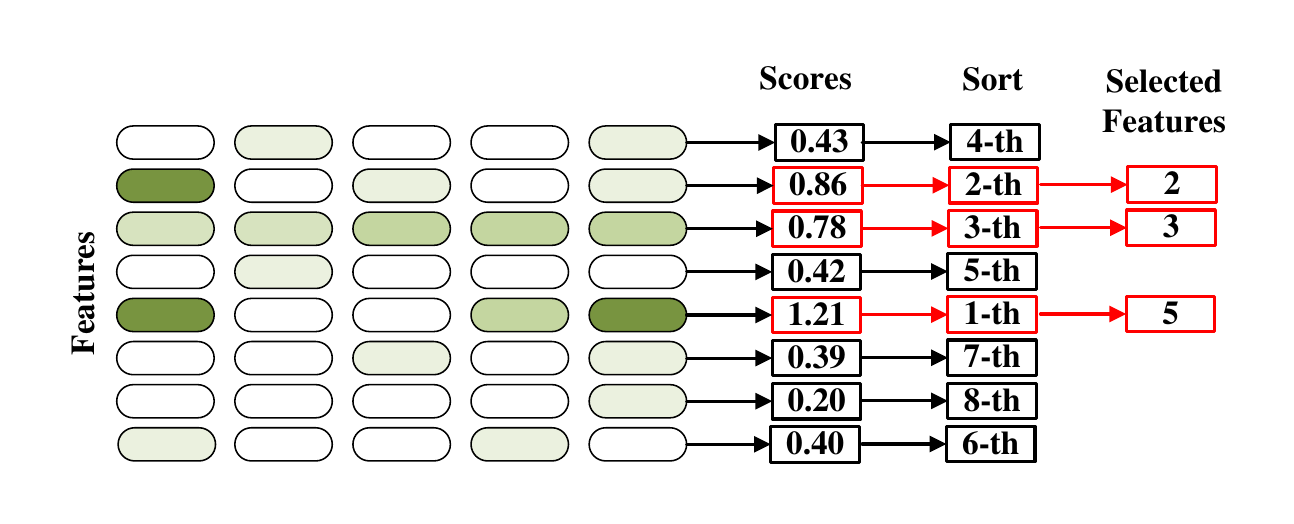}   
		\end{minipage}
	}\hspace{1cm}
	\subfigure[Our proposed BSUFS] 
	{
		\begin{minipage}{11cm}
			\centering      
			\includegraphics[scale=0.43]{figure/fig1.pdf}   
		\end{minipage}
	}
	\caption{Visualization of feature selection results by $(a)$ SPCAFS \cite{li2023sparse} and $(b)$ our proposed  BSUFS. Feature selection is performed by first computing the tranformation matrix $W$, where the dimensions of $W$ correspond to the original data features. Then, by calculating the row-wise norms of $W$, these values are sorted and the top-ranked values are thus selected. Finally, the features associated with these top-ranked values are chosen as the selected features.}
	\label{framework}
\end{figure*}

{However, most works introduced a single row-wise sparse regularization term (e.g., $\ell_{2,1}$-norm) to the transformation matrix in PCA, the element-wise sparsity of this matrix is not considered. Therefore, a bi-sparse optimization model with $\ell_{2,1}$-norm and $\ell_{1}$-norm regularization terms was proposed in \cite{zhu2017double}, where $\ell_{2,1}$-norm is used to select features and $\ell_{1}$-norm is enforced to filter out noise. In fact, the combination of $\ell_{2,1}$-norm and $\ell_{1}$-norm  has been considered in multi-task feature selection and has shown excellent performances \cite{wang2012identifying,wang2013multi}. Although the idea of bi-sparse or double sparse has been applied in image denoising \cite{rubinstein2009double}, compressed sensing \cite{bian2022double}, gene expression \cite{hu2019dstpca}, and radar imaging \cite{9416245}, the involved double sparsity terms are constrained to different variables. Note that $\ell_{2,p}$-norm with $p\in[0,1)$ is a general form of $\ell_{2,p}$-norm in \cite{li2023sparse} and $\ell_{2,0}$-norm \cite{nie2022learning}. \textit{A natural question is whether we can propose a bi-sparse non-convex PCA-based feature selection framework by inheriting the advantages of $\ell_{2,p}$-norm with $p\in[0,1)$ and $\ell_{q}$-norm with $q\in[0,1)$, and verify the effectiveness of $\ell_{q}$-norm for feature selection.}

Motivated by the above observations, we propose a unified bi-sparse UFS method called BSUFS, which is a PCA model with $\ell_{2,p}$-norm and $\ell_{q}$-norm regularization terms, where $p$ and $q$ are within the range of $[0, 1)$. As can be seen from Fig. \ref{framework}, compared with the benchmark SPCAFS which only contains a single $\ell_{2,p}$-norm regularization term with $p\in(0,1)$ on the transformation matrix, adding an additional $\ell_{q}$-norm regularization makes the feature selection results of BSUFS different. Actually, our proposed BSUFS takes SPCAFS and FSPCA as special cases and has higher flexibility than the convex relaxation form in \cite{zhu2017double}. Of course, the choice of $p$ and $q$ may be affected by the specific data structure, which will be discussed in Subsection \ref{parameter}, and it is confirmed that combining the values of 0 and $(0, 1)$ into $[0, 1)$ is of great significance. 
In addition, the idea of bi-sparse optimization is not limited to feature selection, but can  be easily extended to other fields of image processing.

The main contributions of this paper can be summarized in the following three aspects:
\begin{itemize}
 \item We develop a novel BSUFS method that can improve the performance of feature selection by regularizing $\ell_{2,p}$-norm and $\ell_{q}$-norm. The values of $p$ and $q$ are set in the range of $[0,1)$, which has not been considered before.
 \item We design an efficient proximal alternating minimization (PAM) algorithm and all subproblems admit closed-form proximal operators or can be solved by Stiefel manifold optimization. In addition, we provide its computational complexity for every iteration.
 \item We conduct sufficient numerical experiments to evaluate the performance of BSUFS. In particular, we analyze the parameter effects of $p, q\in [0,1)$ and show that, in feature selection, $p$ plays a dominant role, while $q$ serves a complementary role, both of which are indispensable.
\end{itemize}

The structure of this paper is outlined as follows. Section \ref{preliminaries} introduces the notations and related works. Section \ref{method} presents the proposed model, optimization algorithm, and computational complexity. Section \ref{experiment} provides the numerical results, ablation experiments, statistical tests, effects of $p$ and $q$, and discussions. Section \ref{conclusion} concludes this paper.

\section{Preliminaries}\label{preliminaries}

This section first introduces the notations used throughout this paper and then gives some related preliminaries.

\subsection{Notations}
In this paper, matrices are represented by capital letters, vectors by boldface letters, and scalars by lowercase letters. Let $\mathbb{R}^{d}$ and $\mathbb{R}^{m\times n}$ be the sets of all  $d$-dimensional vectors and $m\times n$-dimensional matrices, respectively. For any vector $\mathbf{x}\in\mathbb{R}^{d}$, its $i$-th element is denoted as $x_i$. The $\ell_{2}$-norm of $\mathbf{x}$ is $\|\mathbf{x}\|=(\sum_{i=1}^{d}x_{i}^2)^{1/2}$. 
For any matrix $X\in\mathbb{R}^{m\times n}$, $x_{ij}$ represents its $ij$-th element,  $\mathbf{x}^i$ and $\mathbf{x}_i$ represent its $i$-th row and $i$-th column, respectively. The Frobenius norm of $X$ is $\|X\|_{\textrm{F}} = (\sum_{i=1}^m \sum_{j=1}^n x_{ij}^2)^{1/2}$.  
For $p,q\in (0,1)$, the $\ell_{2,p}$-norm of $X$ is defined as $\|X\|_{2,p}=(\sum_{i=1}^{m}\|\mathbf{x}^i\|^p)^{1/p}$ and the $\ell_{q}$-norm is $\|X\|_{q}=(\sum_{i=1}^m \sum_{j=1}^n |x_{ij}|^q)^{1/q}$. When $p,q=1$, they are exactly $\ell_{2,1}$-norm and $\ell_1$-norm. In addition, $\|X\|_{2,0}$ and $\|X\|_{0}$ count the numbers of non-zero rows and  non-zero elements of $X$, respectively. 
A further notation will be introduced wherever it appears.

\subsection{PCA Basis}

Consider a data matrix $X = (\mathbf{x}_1,\mathbf{x}_2,\ldots,\mathbf{x}_n)\in \mathbb{R}^{d\times n}$  with $\mathbf{x}_i \in \mathbb{R}^{d}$. Denote $W$ as a transformation matrix and $W=(\mathbf{w}_1, \mathbf{w}_2, \ldots, \mathbf{w}_m)\in \mathbb{R}^{d\times m}$, where $\mathbf{w}_i \in \mathbb{R}^{d}$, $\|\mathbf{w}_i\|=1$, and $\mathbf{w}_i^{\top} \mathbf{w}_j=0$ for $i \neq j$. Here, it is supposed that $m<n$.

Mathematically, PCA can be represented in a trace formulation with an orthogonal constraint as
 \begin{equation}\label{pca}
         \begin{aligned}
         \min_{W\in\mathbb{R}^{d\times m}}~&-\textrm{Tr}\left(W^{\top} X X^{\top} W\right)\\
         \textrm{s.t.}~~~~&W^{\top}W=I_m,
         \end{aligned}
 \end{equation}
where the data is assumed to be centralized. For the general case that the data is not centralized, PCA  can be written as 
\begin{equation}\label{pca-1}
        \begin{aligned}
        \min_{W\in\mathbb{R}^{d\times m}}~&-\textrm{Tr}\left(W^{\top} S W\right) \\
        \textrm{s.t.}~~~~&W^{\top}W=I_m,
\end{aligned}
\end{equation}
where $S = XHX^{\top}$ and $H = I_{n}- \frac{1}{n}\mathbf{11}^{\top}$ with $\mathbf{1}\in \mathbb{R}^{n}$ being a vector whose elements are all ones.

Without loss of generality, denote  $\mathbf{w}^{i\top}$ as the transpose of $\mathbf{w}^{i}$, where  $\mathbf{w}^{i}$ is the $i$-th row of $W$.  In feature selection, it can be found that $\mathbf{w}^{i\top}$ represents the transformation vector associated with the $i$-th feature in $X$. 
In details, the matrix $W$ is denoted as
\begin{equation}
        \begin{aligned}
           W =(\mathbf{w}_1, \mathbf{w}_2, \ldots, \mathbf{w}_m)=
            \begin{pmatrix} \mathbf{w}^{1} \\ \mathbf{w}^{2} \\ \vdots \\ \mathbf{w}^{d}\end{pmatrix}\in \mathbb{R}^{d\times m}.
        \end{aligned}
\end{equation}
By transforming the vector 
\begin{equation}
        \begin{aligned}
          \mathbf{x}_i=
            \begin{pmatrix} x_{1,i} \\  x_{2,i} \\ \vdots \\ x_{d,i}\end{pmatrix}\in \mathbb{R}^{d}
        \end{aligned}
\end{equation}
via the transformation matrix $W$, we get the transformation vector $\mathbf{z}_i$ as
\begin{equation}\label{ss}
        \begin{aligned}
            \mathbf{z}_i = W^{\top} \mathbf{x}_i=( \mathbf{w}^{1\top}, ~ \mathbf{w}^{2\top}, ~ \ldots ~, \mathbf{w}^{d\top} )
            \begin{pmatrix} x_{1,i} \\ x_{2,i} \\ \vdots \\ x_{d,i} \end{pmatrix}.
        \end{aligned}
\end{equation}
Accordingly, $\|\mathbf{w}^{i}\|$ can be used to measure the importance of the $i$-th feature \cite{li2013clustering}.

\subsection{Sparse PCA}

In recent years, non-convex optimization has been rapidly developed, which can provide more possibilities than convex optimization \cite{jain2017non,zhou2021sparse,bai2024new}. With this advantage, Li et al. \cite{li2023sparse}  proposed a sparse PCA model called SPCAFS, which is given as the form of
\begin{equation}\label{SPCAFS}
	\begin{aligned}
		\min_{W\in\mathbb{R}^{d\times m}}~&-\textrm{Tr}(W^{\top} SW)+\lambda \|W\|_{2,p}^p \\
		\rm{s.t.}~~~~&W^{\top} W=I_m,
	\end{aligned}	
\end{equation}
where $\lambda \geq0$ is the regularization parameter and $p\in(0,1)$. By introducing $\ell_{2,p}$-norm into the objective function of \eqref{pca-1}, SPCAFS can obtain a row-wise sparse transformation matrix $W$, thereby improving the performance of feature selection. It is also numerically demonstrated that its performance is better than the convex relaxation when $p=1/2$.

FSPCA is another popular UFS method, whose mathematical model is given by
\begin{equation}\label{FSPCA}
	\begin{aligned}
		\min_{W\in\mathbb{R}^{d\times m}}~&-\textrm{Tr}(W^{\top} SW) \\
		\rm{s.t.}~~~~&\|W\|_{2,0}\leq s,~W^{\top} W=I_m,
	\end{aligned}	
\end{equation}
where $s>0$ is the sparsity level. Recall Fig. \ref{framework} and \eqref{ss}, it is seen that $s$ corresponds to the number of selected features.

\section{The Proposed Method} \label{method}

This section first presents our proposed bi-sparse and non-convex UFS model, and then shows an optimization algorithm.

\subsection{Model Construction}

As is demonstrated in \cite{zhu2017double} that, different from these single structured sparse PCA methods, bi-sparse regularized models can better select representative features due to the introduction of $\ell_{2,1}$-norm and $\ell_{1}$-norm on the transformation matrix. Therefore, to inherit the advantages of non-convex optimization and bi-sparse optimization, we propose a sparse PCA variant called BSUFS, which is reformulated as
\begin{equation}\label{raw_problem}
        \begin{aligned}
                \min_{W\in\mathbb{R}^{d\times m}}~ &-\textrm{Tr}(W^{\top} SW)+\lambda_1 \|W\|_{2,p}^{p} +\lambda_2 \|W\|_{q}^{q} \\
                \textrm{s.t.}~~~~&W^{\top} W=I_m,
        \end{aligned}
\end{equation}
where $p,q\in [0,1)$ and $\lambda_{1},\lambda_{2}\geq0$ are the parameters to control the bi-sparse regularization terms.  

Compared with the existing PCA-based UFS methods, the following conclusions can be made:
\begin{itemize}
    \item When $\lambda_2=0$ in \eqref{raw_problem}, our proposed BSUFS unifies SPCAFS \cite{li2023sparse} and the Lagrangian form of FSPCA \cite{nie2022learning}.
     \item When $p$ and $q$ tend to 1 in \eqref{raw_problem}, our proposed BSUFS equals to the model in \cite{zhu2017double}.
\end{itemize}

In summary, BSUFS  can provide flexible sparse solutions by choosing different $p$ and $q$ in the range of $[0,1)$, thus facilitating feature selection.

\subsection{Optimization Algorithm}

Besides the Stiefel manifold constraint $W^{\top} W=I_m$, the introduced $\ell_{2,p}$-norm and $\ell_{q}$-norm in \eqref{raw_problem} makes our proposed BSUFS more difficult. Below, we provide an efficient optimization algorithm based on the proximal alternating minimization (PAM) technique.

By utilizing auxiliary variables $W = V,~W = U$, \eqref{raw_problem} can be rewritten as
\begin{equation}\label{raw_problem-wuv}
        \begin{aligned}
                \underset{W,U,V\in\mathbb{R}^{d\times m}}{\min} ~&-\textrm{Tr}(W^{\top}SW) + \lambda_1 \|V\|_{2,p}^{p} + \lambda_2 \|U\|_{q}^{q} \\
                \textrm{s.t.}~~~~~~& W^{\top}W=I_m,~W = V,~W = U.
        \end{aligned}
\end{equation}

Then, introduce the indicator function
\begin{equation}
\Phi(W)=\left\{
        \begin{aligned}
        &0, ~~~~~~~~W^{\top}W=I_m,\\
        &+\infty, ~~~\text{otherwise},
        \end{aligned}
\right.
\end{equation}
and rewrite \eqref{raw_problem-wuv}  its penalized form
\begin{equation}\label{new_problem}
        \begin{split}
               & \underset{W,U,V\in\mathbb{R}^{d\times m}}{\min} ~ -\textrm{Tr}(W^{\top}SW) +\lambda_1 \|V\|_{2,p}^{p} + \lambda_2 \|U\|_{q}^{q}  \\
                &\qquad \qquad + \frac{\beta_{1}}{2} \|W-U\|^{2}_{\textrm{F}} + \frac{\beta_{2}}{2} \|W-V\|^{2}_{\textrm{F}}  + \Phi(W),
        \end{split}
\end{equation}
where $\beta_{1}, \beta_{2} > 0$ are the penalty parameters. Denote the objective function of \eqref{new_problem} as $f(W,U,V)$. Under the PAM framework, each variable can be updated alternately via the following scheme
\begin{equation}\label{eq:1}
W^{k+1}\in\underset{W\in\mathbb{R}^{d\times m}}{\arg\min} ~ f(W, U^{k},V^{k}) + \frac{\tau_{1}}{2}\|W-W^{k}\|^{2}_{\textrm{F}},
\end{equation}
\begin{equation}\label{eq:2}
U^{k+1}\in\underset{U\in\mathbb{R}^{d\times m}}{\arg\min} ~f(W^{k+1},U, V^{k}) + \frac{\tau_{2}}{2}\|U-U^{k}\|^{2}_{\textrm{F}},
\end{equation}
\begin{equation}\label{eq:3}
V^{k+1}\in\underset{V\in\mathbb{R}^{d\times m}}{\arg\min} ~ f(W^{k+1},U^{k+1}, V) + \frac{\tau_{3}}{2}\|V-V^{k}\|^{2}_{\textrm{F}},
\end{equation}
where $\tau_{1}, \tau_{2}, \tau_{3}  > 0$ and $k$ is the iteration number.	 The overall scheme is presented in Algorithm \ref{algorithm}, and the update rules for $W$, $U$, and $V$ are analyzed as follows.

\begin{algorithm}[t]
    \caption{Proximal alternating minimization (PAM) algorithm for BSUFS} \label{algorithm}
    \textbf{Input:} Data $X \in \mathbb{R}^{d \times n}$, parameters $p$, $q$, $\lambda_1$, $\lambda_2$, $\beta_1$, $\beta_2$,  $\tau_1$, $\tau_2$, $\tau_3$, calculate $H=I_{n}-\frac{1}{n}\mathbf{11}^{\top}$, $S = XHX^{\top}$\\
    \textbf{Output:} $V$\\
    \textbf{Initialize:} $k=0$, $W^k$, $U^k$, $V^k$ \\
    \textbf{While} not converged \textbf{do}
            \begin{algorithmic}[1]
                \STATE  Update $W^{k+1}$ by \eqref{eq:1}
                \STATE  Update $U^{k+1}$ by \eqref{eq:2}
                \STATE  Update $V^{k+1}$ by \eqref{eq:3}
                \STATE  Check convergence: If $$\frac{|f^{k+1} - f^{k}|}{ \max \{ |f^{k}|,1 \} } < 10^{-4}~ \textrm{or} ~k > 500$$ then stop. Otherwise, $k = k + 1$
            \end{algorithmic}
    \textbf{End While}
\end{algorithm}

\subsubsection{Update $W$ by fixing $U$ and $V$}
\eqref{eq:1} can be expressed as 
\begin{equation}\label{update_W}
        \begin{aligned}
               & \min_{W^{\top}W = I_m}~ g(W)=-\textrm{Tr}(W^{\top}SW) + \frac{\beta_{1}}{2} \|W-U^{k}\|^{2}_{\textrm{F}} \\ 
         & \qquad\quad\quad + \frac{\beta_{2}}{2} \|W-V^{k}\|^{2}_{\textrm{F}} + \frac{\tau_{1}}{2}\|W-W^{k}\|^{2}_{\textrm{F}}.
        \end{aligned}
\end{equation}
The Euclidean gradient of $g(W)$ is
\begin{equation}\label{gradient_W}
    \begin{aligned}
        \nabla g(W) = &- 2S W + \beta_{1} (W-U^{k}) \\
        &+\beta_{2} (W-V^{k}) + \tau_{1} (W-W^{k}),
    \end{aligned}
\end{equation}
and the Euclidean Hessian is
\begin{equation}
\begin{aligned}
    \nabla^2 g(W) = -2 I_m \otimes S + (\beta_{1} + \beta_{2} + \tau_{1}) I_{dm},
\end{aligned}
\end{equation}
where $\otimes$ represents the Kronecker product.

Denote
\begin{equation}
\textrm{St}(d, m) = \{W \in \mathbb{R}^{d \times m} \mid  W^{\top}W = I_m\}.
\end{equation}
Then, \eqref{update_W} is a Stiefel manifold optimization problem and can be reformulated as 
\begin{equation}\label{W2}
\begin{aligned}
        \underset{W \in \textrm{St}(d,m)}{\min} ~ & g(W),
\end{aligned}
\end{equation}
where the trust-region Riemannian manifold optimization algorithm \cite{absil2008optimization} can be applied to solve this problem.

\begin{algorithm}[t]
\caption{Trust-region Riemannian manifold optimization algorithm for solving  \eqref{eq:1}} \label{SW}
\textbf{Input:} Data $S \in \mathbb{R}^{d\times d}$, $U^{k} \in \mathbb{R}^{d\times m}$, $V^{k} \in \mathbb{R}^{d\times m}$, parameters $\beta_1$, $\beta_2$, $\tau_1$, $\varepsilon$, $\Delta'>0$, $\rho' \in [0,\frac{1}{4})$\\
\textbf{Output:} $W^{k}$ \\
\textbf{Initialize:} $i=0$,  $W^{k}_{i}\in\textrm{St}(d, m)$,  $\Delta_i \in (0,\Delta')$\\
\textbf{While} not converged \textbf{do}	
\begin{algorithmic}[1]
\STATE Obtain $M_i$ by solving \eqref{trustsub} with $W=W^{k}_{i}$ and $\Delta=\Delta_i$
\STATE Compute  $\rho_i$ from \eqref{rho} with $W=W^{k}_{i}$
\IF{$\rho_i < \frac{1}{4}$}
\STATE $\Delta_{i+1} = \frac{1}{4} \Delta_i$
\ELSIF {$\rho_i > \frac{3}{4}$ and $\|M_i\| = \Delta_i$}
\STATE $\Delta_{i+1} = \min(2\Delta_i, \Delta')$
\ELSE
\STATE $\Delta_{i+1} = \Delta_i$
\ENDIF
\IF {$\rho_i > \rho'$}
\STATE $W^{k}_{i+1} = R_W(M_i)$
\ELSE
\STATE $W^{k}_{i+1} = W^{k}_{i}$
\ENDIF
\STATE  Check convergence: If 
$$\operatorname{grad} g(W^{k}_{i+1}) < 10^{-6}~  \textrm{or} ~i> 100$$ then stop. Otherwise, $i = i +1$
\end{algorithmic}
    \textbf{End While}
\end{algorithm}

To implement this algorithm, there are two basic elements, i.e., the search direction and the trust-region ratio. The search direction relies on the Riemannian gradient and Riemannian Hessian of the objective function $g(W)$ in \eqref{W2}. Specifically, the Riemannian gradient can be  obtained by projecting its Euclidean gradient onto the tangent space of the Stiefel manifold, i.e.,
\begin{equation} \label{grad}
\begin{aligned}
    \operatorname{grad} g(W) & =\mathcal{P}_{W} (\nabla g(W)) \\
    & =\nabla g(W)- W \operatorname{sym}(W^{\top}  \nabla g(W) ),
\end{aligned}
\end{equation}
where $\mathcal{P}_{W} (\nabla g(W))$ is the projection of the Euclidean gradient onto the tangent space of the Stiefel manifold, and $\operatorname{sym}(X)=(X+X^{\top}) / 2$ extracts the symmetric part of a square matrix $X$. 
The Riemannian Hessian can be obtained by projecting its Euclidean Hessian onto the tangent space of the Stiefel manifold, i.e.,
\begin{equation}
\begin{aligned}
    \operatorname{Hess} g(W)[M]&= \mathcal{P}_{W}  (\nabla^2 g(W)[M]\\
    &~~-M\operatorname{sym}(W^{\top}\nabla g(W))\\
    &~~-W\operatorname{sym}(M^{\top}\nabla g(W))\\
    &~~-W\operatorname{sym}(W^{\top}\nabla^2 g(W)[M])),
\end{aligned}
\end{equation}
where $\nabla^2 g(W)[M]$ is the Euclidean Hessian-vector product. Then, the search direction of the trust region algorithm is derived by solving the following problem
\begin{equation}
			\begin{aligned}\label{trustsub}
				\min_{M \in \operatorname{T_W} \textrm{St}(d,m)} ~ &m_W(M) = g(W)+\left\langle \operatorname{grad} g(W),M\right\rangle_{W} \\
				&\quad\quad+\frac{1}{2} \left\langle \operatorname{Hess} g(W)[M],M\right\rangle_{W}\\
				\textrm{s.t.} ~~~~~~& \left\langle M,M\right\rangle_{W} \leq \Delta^{2},
			\end{aligned}
		\end{equation}
		where $\Delta$ is the trust-region radius, and $\operatorname{T_W} \textrm{St}(d,m)$ is the tangent space of the manifold at $W$.

Finally, the trust region ratio is determined by
\begin{equation}
			\begin{aligned} \label{rho}
				\rho = \frac{g(W) - g(R_W(M))}{m_W(0) - m_W(M)},
			\end{aligned}
		\end{equation}
		where $R_W(M)$ is the retraction of $M$ onto the Stiefel manifold.

\subsubsection{Update $U$ by fixing $W$ and $V$}
After updating $W$, \eqref{eq:2} can be solved by
\begin{equation}\label{U1}
        \min_{U\in\mathbb{R}^{d\times m}} ~ \lambda_{2} \|U\|_{q}^{q}+  \frac{\beta_{1}}{2} \|W^{k+1}-U\|^{2}_{\textrm{F}} + \frac{\tau_{2}}{2}\|U-U^{k}\|^{2}_{\textrm{F}}.
\end{equation}
By expanding the norm terms in \eqref{U1}, this subproblem can be rewritten as the form of
\begin{equation}\label{U2}
\begin{split}
        \min_{U\in\mathbb{R}^{d\times m}}~ \lambda_{2} \|U\|_{q}^{q} +\frac{\beta_{1} + \tau_{2}}{2} \|U - Y \|^{2}_{\textrm{F}},
\end{split}
\end{equation}
where
$$Y = \frac{\beta_{1}}{\beta_{1} +\tau_{2}} W^{k+1} + \frac{\tau_{2}}{\beta_{1} +\tau_{2}} U^{k}.$$
By considering each element of the matrix separately, the optimization problem can be reformulated as
\begin{equation}\label{U3}
        \begin{aligned}
                \min_{u_{ij}\in\mathbb{R}} ~ \lambda_{2} |u_{ij}|^{q} + \frac{\beta_{1} + \tau_{2}}{2} (u_{ij} - y_{ij})^{2}.
        \end{aligned}
\end{equation}
The solution of this problem is the proximal operator of the $|\cdot|^{q}$, whose result is reviewed  in the following lemma.
\begin{Lemma}\label{P}
Consider the proximal operator 
\begin{equation}
\begin{aligned} \label {U4}
    \operatorname{Prox}_{\lambda|\cdot|^{q}}(a) & =\underset{x\in\mathbb{R}}{\operatorname{argmin}}~\lambda|x|^q+\frac{1}{2}(x-a)^2\\
    & = \begin{cases}\{0\}, & |a|<\kappa(\lambda, q), \\
    \{0, \operatorname{sgn}(a) c(\lambda, q)\}, & |a|=\kappa(\lambda, q), \\
    \left\{\operatorname{sgn}(a) \varpi_q(|a|)\right\}, & |a|>\kappa(\lambda, q),\end{cases}
\end{aligned}
\end{equation}
where
\begin{equation}
\begin{aligned}
c(\lambda, p)&=(2 \lambda(1-q))^{\frac{1}{2-q}}>0,\\
\kappa(\lambda, q)&=(2-q) \lambda^{\frac{1}{2-q}}(2(1-q))^{\frac{q+1}{q-2}},\\
\varpi_q(a) &\in \{x \mid x-a+\lambda q \operatorname{sgn}(x) x^{q-1}=0, x>0 \}.
\end{aligned}
\end{equation}
More details can be found in \cite{zhou2023revisiting}.
\end{Lemma}

According to Lemma \ref{P}, the solution of \eqref{U3} can be easily obtained by
\begin{equation}\label{U5}
        \begin{aligned}
                u_{ij} \in\operatorname{Prox}_
                {\frac{\lambda_2} {\beta_{1} + \tau_{2}}|\cdot|^{q}}\left(y_{ij}\right).
        \end{aligned}
\end{equation}
Based on \cite{beck2017first}, the proximal operator in \eqref {U4} admits a closed-form when $q=0$. This closed-form also exists when $q=1/2$ and $q=2/3$ as illustrated in \cite{xu2012l_,cao2013fast}. For other choices $q\in(0,1)$, efficient algorithms proposed in \cite{zhou2023revisiting,liu2024efficient} can be considered.

\subsubsection{Update $V$ by fixing $W$ and $U$}
Once $W$ and $U$ have been updated, \eqref{eq:3} can be calculated via
\begin{equation}\label{V1}
        \min_{V\in\mathbb{R}^{d\times m}} ~ \lambda_{1} \|V\|_{2,p}^{p} +  \frac{\beta_{2}}{2} \|W^{k+1}-V\|^{2}_{\textrm{F}} + \frac{\tau_{3}}{2}\|V-V^{k}\|^{2}_{\textrm{F}}.
\end{equation}
It is easy to reformulate \eqref{V1} as
\begin{equation}\label{V2}
\begin{aligned}
         \min_{V\in\mathbb{R}^{d\times m}} ~ \lambda_{1} \|V\|_{2,p}^{p} + \frac{\beta_{2} + \tau_{3}}{2} \|V - Z \|^{2}_{\textrm{F}},
\end{aligned}
\end{equation}
where
$$Z = \frac{\beta_{2}}{\beta_{2} +\tau_{3}} W^{k+1} + \frac{\tau_{3}}{\beta_{2} +\tau_{3}} V^{k}.$$
It can be further decomposed into a series of vector optimization problems as
\begin{equation} \label{V3}
        \begin{aligned}
                \min_{\mathbf{v}^{i}\in\mathbb{R}^{m}} ~ \lambda_1 \|\mathbf{v}^{i}\|^{p} + \frac{\beta_{2} + \tau_{3}}{2} \|\mathbf{v}^{i} - \mathbf{z}^{i} \|^{2},
        \end{aligned}
\end{equation}
where $i\in\{1,2,\cdots,d\}$. The solution of \eqref{V3} is the proximal operator of $\|\cdot\|^{p}$ as in \cite{Zhang2022structured,beck2017first}, which is reviewed in the next lemma.
\begin{Lemma}\label{P2}
Consider the proximal operator 
\begin{equation}
\begin{aligned} \label {U5}
    \operatorname{Prox}_{\lambda\|\cdot\|^{p}}(\mathbf{a}) & =\underset{\mathbf{x}\in\mathbb{R}^{m}}{\operatorname{argmin}}~\lambda\|\mathbf{x}\|^p+\frac{1}{2}\|\mathbf{x}-\mathbf{a}\|^2\\
    & = \begin{cases}
\operatorname{Prox}_{\lambda|\cdot|^{p}}\left(\|\mathbf{a}\| \right) \cdot \frac{\mathbf{a}}{\|\mathbf{a}\|}, & \mathbf{a}\neq \mathbf{0},\\
\{\mathbf{0}\}, & \mathbf{a}=\mathbf{0}.\end{cases}
\end{aligned}
\end{equation}
Here, $\|\mathbf{x}\|^{0}=1$ when $\mathbf{x}\neq \mathbf{0}$, and $\|\mathbf{x}\|^{0}=0$ when $\mathbf{x}=\mathbf{0}$.
\end{Lemma}

From results in Lemma \ref{P2}, the solution of \eqref{V3} is
\begin{equation}\label{V4}
                \mathbf{v}^{i} \in
                 \begin{cases}
            \operatorname{Prox}_{\frac{\lambda_1} {\beta_{2} + \tau_{3}}|\cdot|^{p}}\left(\|\mathbf{z}^{i}\| \right) \cdot \frac{\mathbf{z}^{i}}{\|\mathbf{z}^{i}\|}, & \mathbf{z}^{i}\neq \mathbf{0},\\
            \{\mathbf{0}\},& \mathbf{z}^{i}= \mathbf{0}.
                 \end{cases}
\end{equation}

\begin{Remark}
For Algorithm \ref{algorithm}, it is suggested to update $U^{k+1}$ first for eliminating the interference caused by noise and then update $V^{k+1}$ for selecting discriminative features.
\end{Remark}

\begin{Remark}\label{remark2}
To establish the convergent result of Algorithm \ref{algorithm}, there is a common pattern similar to  \cite{bolte2014proximal}, consisting of a sufficient decent lemma and the Kurdyka-Lojasiewicz (K-L)  property of the objective function. For this algorithm, the main challenge is to obtain a sufficient decent lemma, because of the fact that the update of $W^{k+1}$ relies  on Algorithm \ref{SW}. In \cite{absil2008optimization}, it is proved that Algorithm \ref{SW} converges to the zero gradient point, i.e., the point satisfies  $\operatorname{grad} g(W)=0$. However, to prove the sufficient decent lemma of Algorithm \ref{algorithm}, there is a requirement that  Algorithm \ref{SW} converges to the global optimizer, which is a more strict result. Therefore, it is deferred as a future work to tackle this challenge.
\end{Remark}

\begin{figure*}[h]
    \centering
    \subfigcapskip=-1pt
    \subfigure[Diamond9]{
        \label{a}
        \centering
        \includegraphics[width=3.3cm]{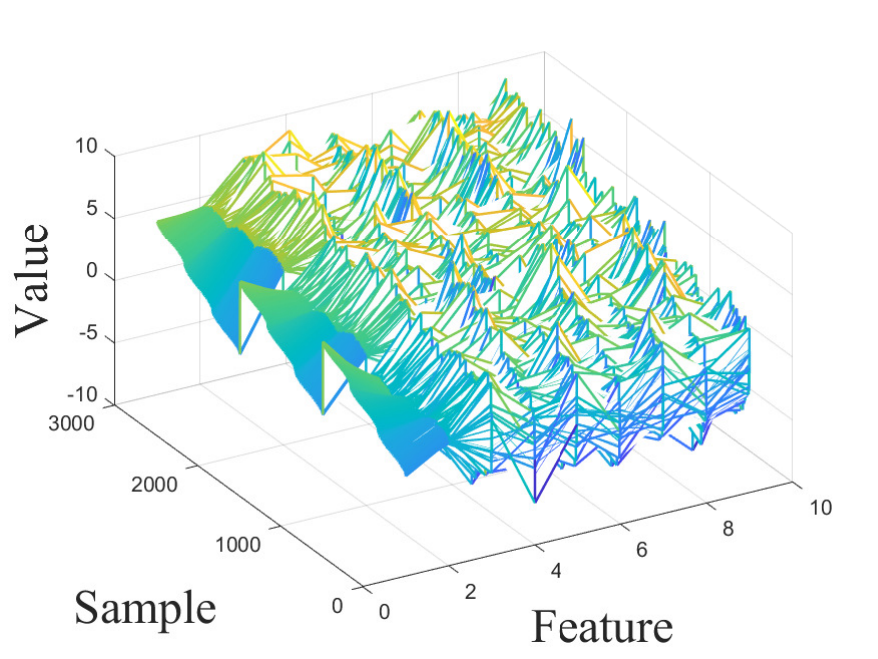}
    }\hspace{-2mm}
    \subfigcapskip=-1pt
    \subfigure[LapScore]{
        \label{b}
        \centering
        \includegraphics[width=3.3cm]{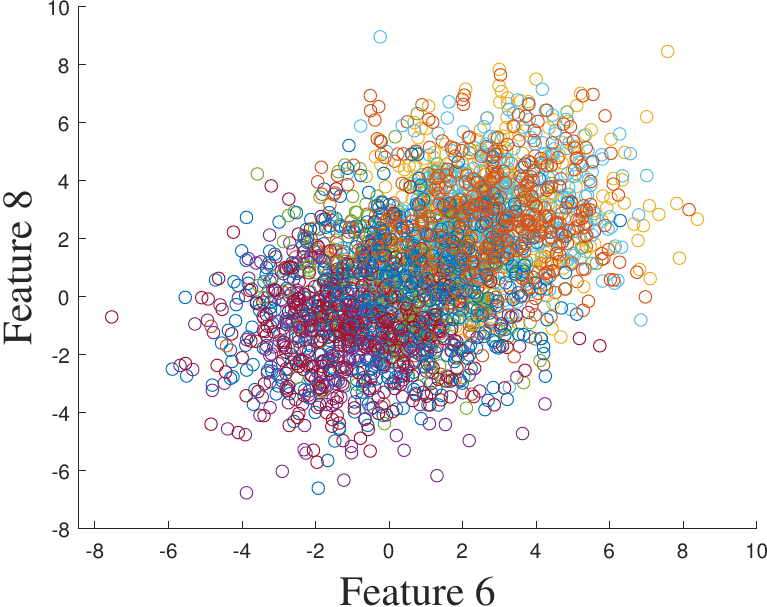}
    }\hspace{-2mm}
    \subfigcapskip=-1pt
    \subfigure[SOGFS]{
        \label{c}
        \centering
        \includegraphics[width=3.3cm]{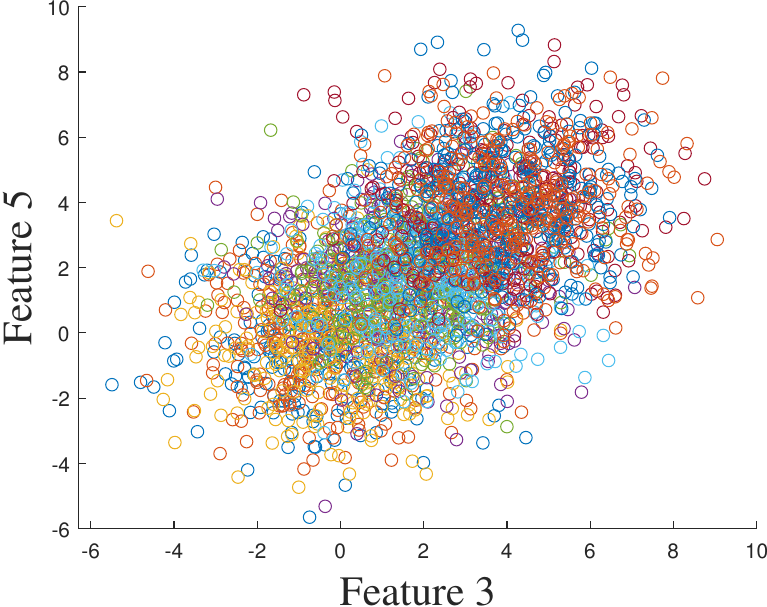}
    }\hspace{-2mm}
    \subfigcapskip=-1pt
    \subfigure[RNE]{
        \label{e}
        \centering
        \includegraphics[width=3.3cm]{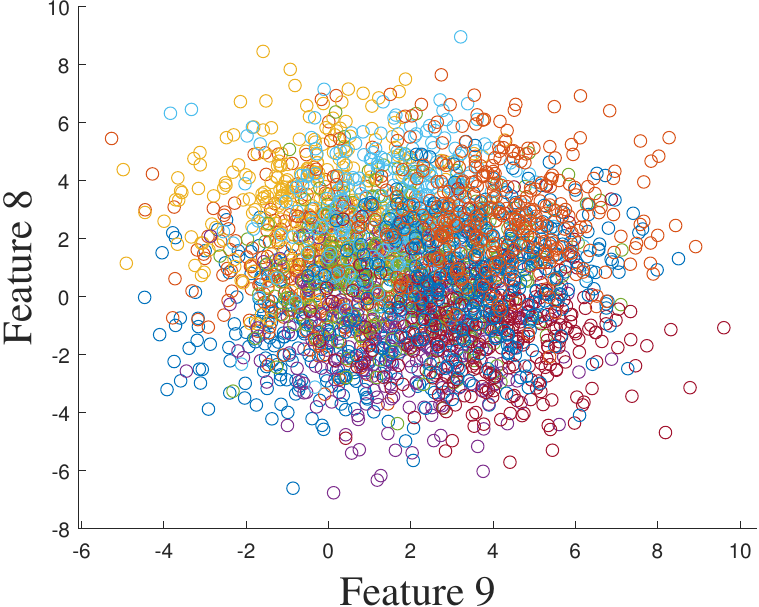}
    }\hspace{-2mm}
    \subfigcapskip=-1pt
    \subfigure[UDFS]{
        \label{f}
        \centering
        \includegraphics[width=3.3cm]{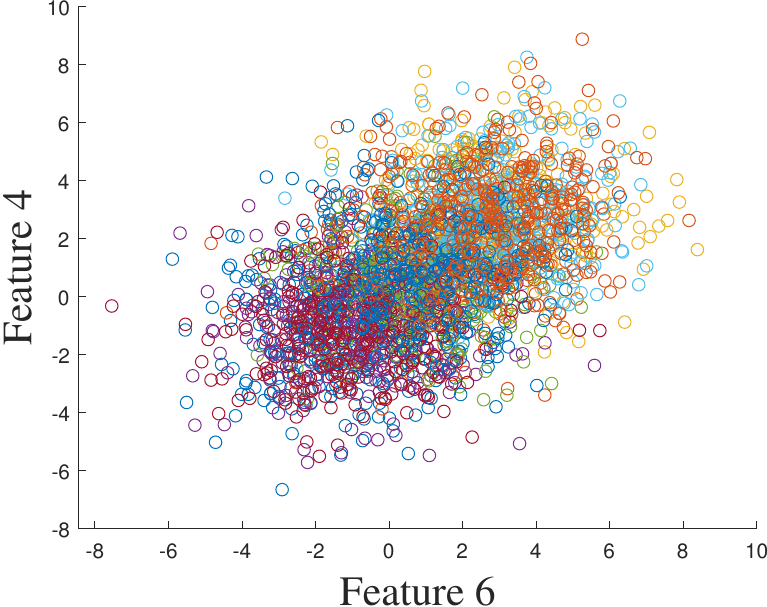}
    }\hspace{-2mm}
     
     \subfigure[SPCAFS]{
        \label{d}
        \centering
        \includegraphics[width=3.3cm]{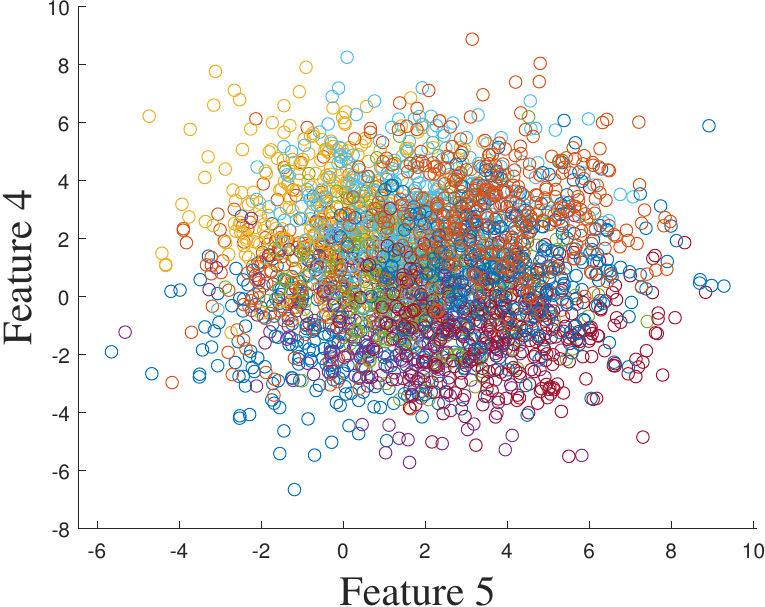}
    }\hspace{-2mm}
    \subfigcapskip=-1pt
    \subfigure[FSPCA]{
        \label{h}
        \centering
        \includegraphics[width=3.3cm]{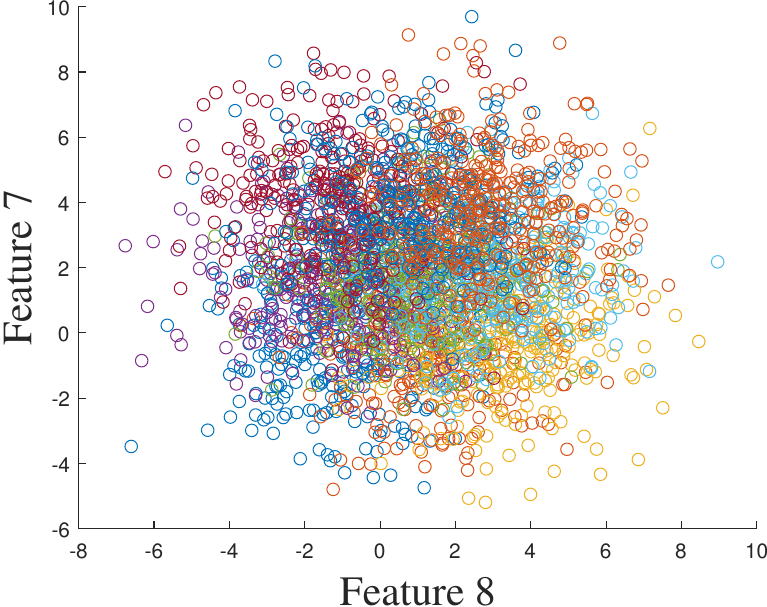}
    }\hspace{-2mm}
            \subfigcapskip=-1pt
    \subfigure[SPCA-PSD]{
        \label{g}
        \centering
        \includegraphics[width=3.3cm]{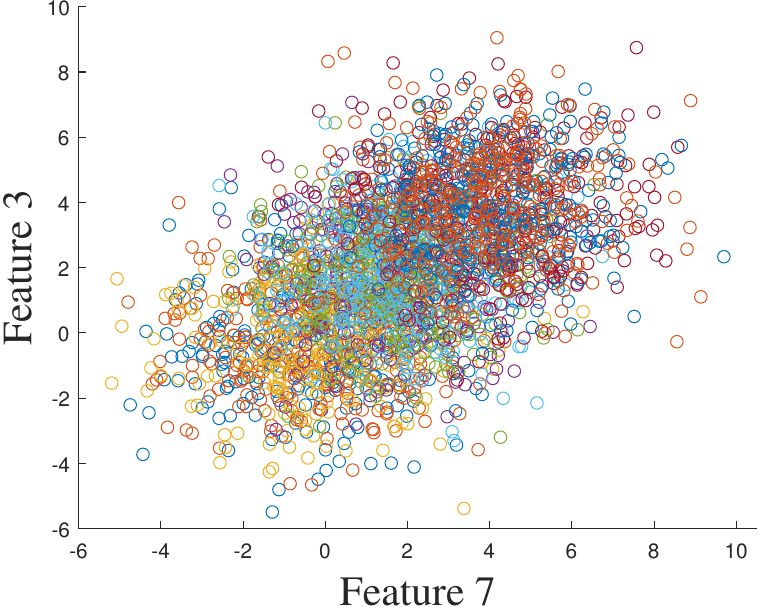}
    }\hspace{-2mm}
    \subfigcapskip=-1pt
    \subfigure[FEN-PCAFS]{
        \label{i}
        \centering
        \includegraphics[width=3.3cm]{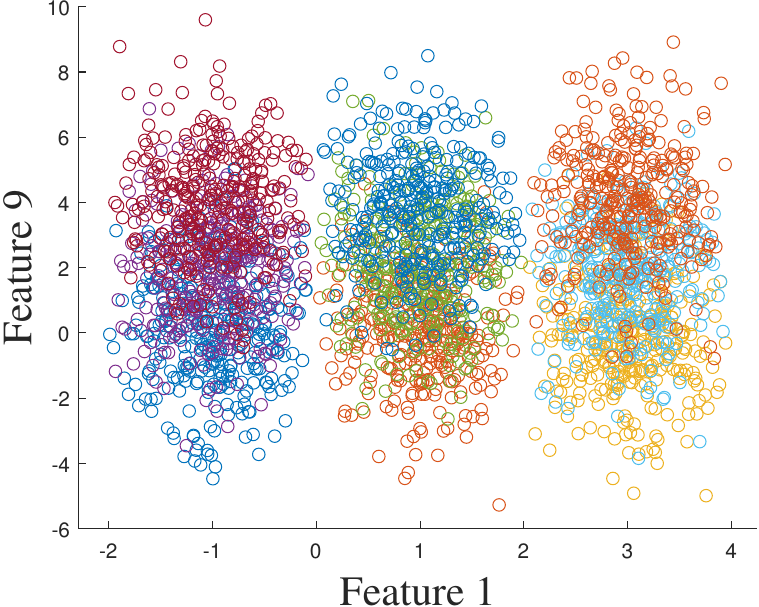}
    }\hspace{-2mm}
    \subfigcapskip=-1pt
    \subfigure[BSUFS]{
        \label{j}
        \centering
        \includegraphics[width=3.3cm]{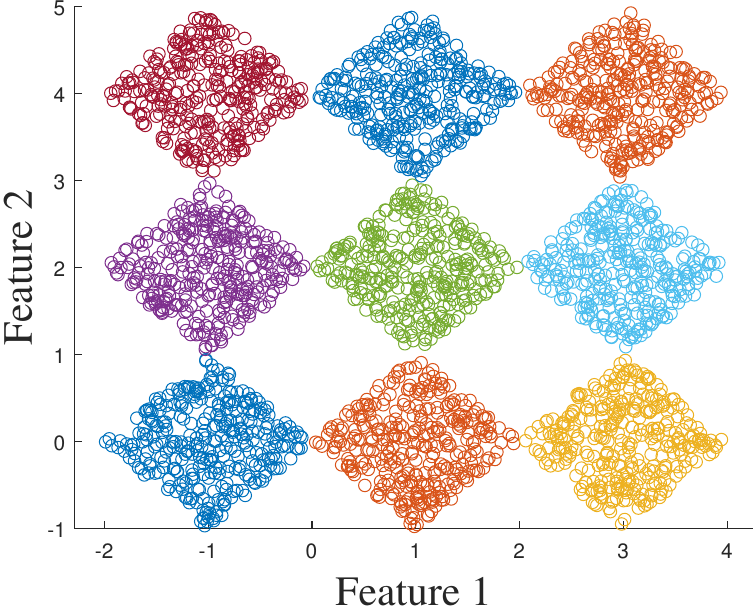}
    }\hspace{-2mm}
    \vspace{-0.1cm}
    \caption{Visual comparisons on the Diamond9 dataset, where (a) is the dataset distribution and (b)-(j) are the feature selection results.}
    \label{diamond9}
\end{figure*}

\subsection{Computational Complexity Analysis}\label{complexity}
At the end of this section, we provide the computational complexity of Algorithm \ref{algorithm} (including Algorithm \ref{SW}), which can be decomposed into the following four aspects.
\begin{itemize}
\item When initializing Algorithm \ref{algorithm}, there require $O(n^2)$ and $O(dn^2)$ to compute $H$ and $S$, respectively, which leads to the computational complexity of initialization is $O(dn^2)$.
\item When updating $W$ by Algorithm \ref{SW}, the computational complexity is that of solving \eqref{trustsub} and \eqref{rho}. As pointed out in \cite{absil2008optimization}, \eqref{trustsub} is approximated solved by the truncated (Steihaug-Toint) conjugate-gradient method with caching, where the computational cost relies on computations of the Hessian-vector product and inner product. Here, the computational complexity of the Hessian-vector product is $O(d^2m+dm^2)$ and the inner product is $O(dm^2)$. For \eqref{rho}, the computational complexity is $O(d^2m+dm^2)$. By summing up these results, the computational complexity of each iteration in Algorithm \ref{SW} is $O(d^2m+dm^2)$.
\item When updating $U$ and $V$, the computational complexity only relies on the proximal operators, which are proved as closed-form functions in our numerical studies, and their computational complexity is $O(dm)$.
\item The convergent check is based on the loss function $f$, which has a computational complexity of $O(d^2m)$.
\end{itemize}
Overall, the  computational complexity for every iteration of  Algorithm \ref{algorithm} is $O((\kappa+1)d^2m+\kappa dm^2+dm)$, where  $\kappa$ is the iteration number of Algorithm \ref{SW}.

\begin{table}[t]    \renewcommand\arraystretch{1.3}
    \caption{Statistical information of selected datasets.}\label{data}
    \centering
    \begin{tabular}{|c|c|c|c|c|}
        \hline
        \textbf{Type}&\textbf{Datasets} & \textbf{Features} &\textbf{Samples} &\textbf{Classes} \\
        \hline\hline
        \multirow{2}*{\textbf{Synthetic}}
         &Diamond9 & 9 & 3000  & 9\\  
        \cline{2-5}
        &Dartboard1  & 9& 1000 & 4 \\
        \hline
        \multirow{8}*{\textbf{Real-world}}&
         COIL20 & 1024 & 1440 & 20 \\
        \cline{2-5}
         &Isolet & 617 & 1560 & 26 \\
        \cline{2-5}
        &USPS & 256 & 1000 & 10 \\
        \cline{2-5}
        &umist & 644 & 575 & 20 \\
        \cline{2-5} 
        &GLIOMA & 4434 & 50 & 4 \\
        \cline{2-5}
        &pie & 1024 & 1166  & 53 \\
        \cline{2-5}
        &LUNG & 325 & 73 & 7 \\
        \cline{2-5}
        &MSTAR & 1024 & 2425 & 10 \\
        \hline
    \end{tabular}
\end{table}

\section{Numerical Studies}\label{experiment}
To illustrate the effectiveness of the proposed BSUFS, this section presents extensive comparisons with several benchmark UFS methods, including graph-based methods as LapScore \cite{he2005laplacian}, SOGFS \cite{nie2016unsupervised}, RNE \cite{liu2020robust}, and UDFS \cite{yang2011}, and PCA-based methods as SPCAFS \cite{li2023sparse}, FSPCA \cite{tian2020learning}, SPCA-PSD\cite{zheng2023fast}, and FEN-PCAFS \cite{gao2024principal}. Note that LapScore, SOGFS, RNE, and UDFS are directly implemented by the AutoUFSTool toolbox\footnote{https://github.com/farhadabedinzadeh/AutoUFSTool}, while
SPCAFS\footnote{https://github.com/quiter2005/algorithm}, FSPCA\footnote{https://github.com/tianlai09/FSPCA}, SPCA-PSD\footnote{https://github.com/zjj20212035/SPCA-PSD}, and FEN-PCAFS\footnote{https://github.com/gaoyl-group/FEN-PCAFS} are from the authors' GitHub repositories.

To verify whether BSUFS is better than these UFS methods and test every component of BSUFS, this section is organized as follows. Subsection \ref{setup} describes the datasets, parameter settings, and evaluation metrics. 
Subsection \ref{experiments} and Subsection \ref{experiments-real} present numerical experiments on synthetic and real-world datasets, respectively.  
Subsection \ref{ablation} provides ablation experiments. 
Subsection \ref{stat} shows the statistical test results. Subsection \ref{parameter} analyzes the effects of $p$ and $q$.  Subsection \ref{discussion} gives more discussions.

\begin{figure*}[t]
    \centering
    \subfigcapskip=-1pt
    \subfigure[Dartboard1]{
        \label{a}
        \centering
        \includegraphics[width=3.3cm]{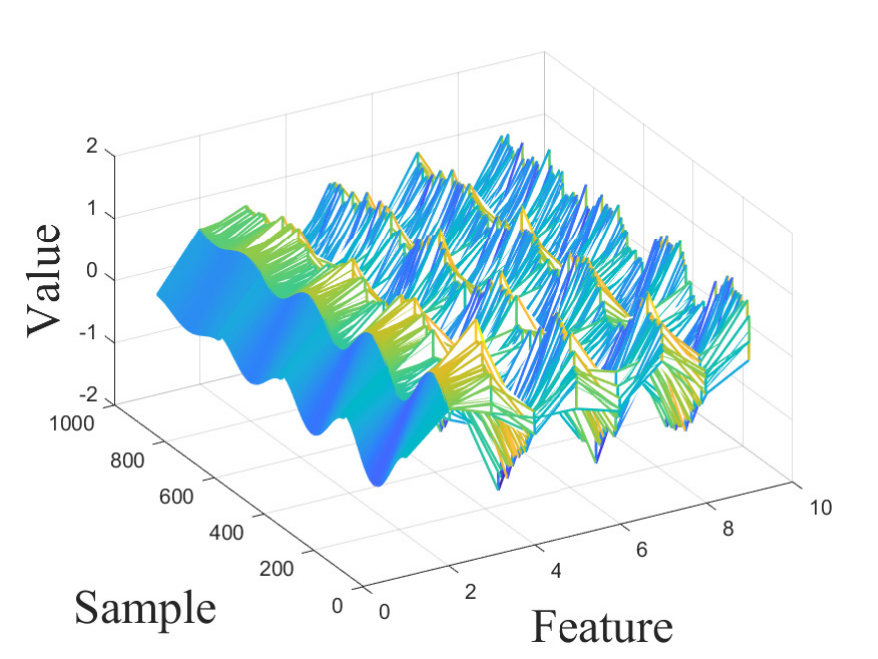}
    }\hspace{-2mm}
    \subfigcapskip=-1pt
    \subfigure[LapScore]{
        \label{b}
        \centering
        \includegraphics[width=3.3cm]{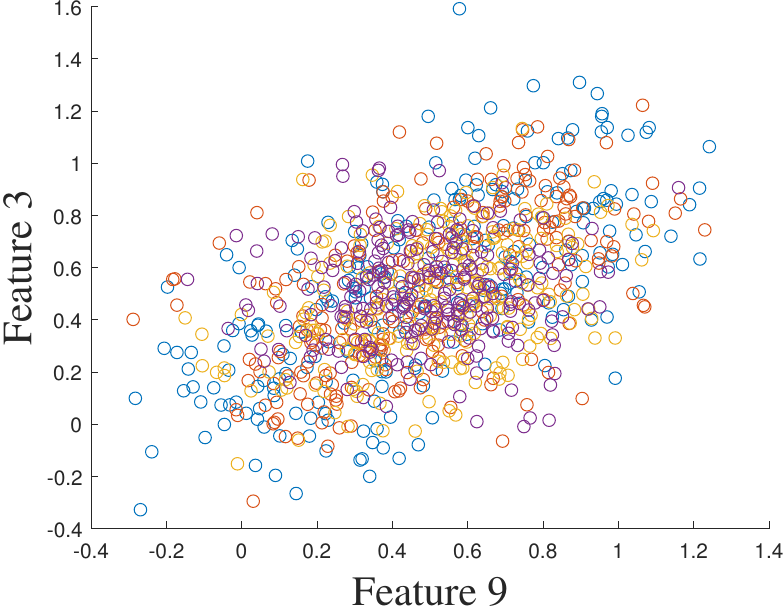}
    }\hspace{-2mm}
    \subfigcapskip=-1pt
    \subfigure[SOGFS]{
        \label{c}
        \centering
        \includegraphics[width=3.3cm]{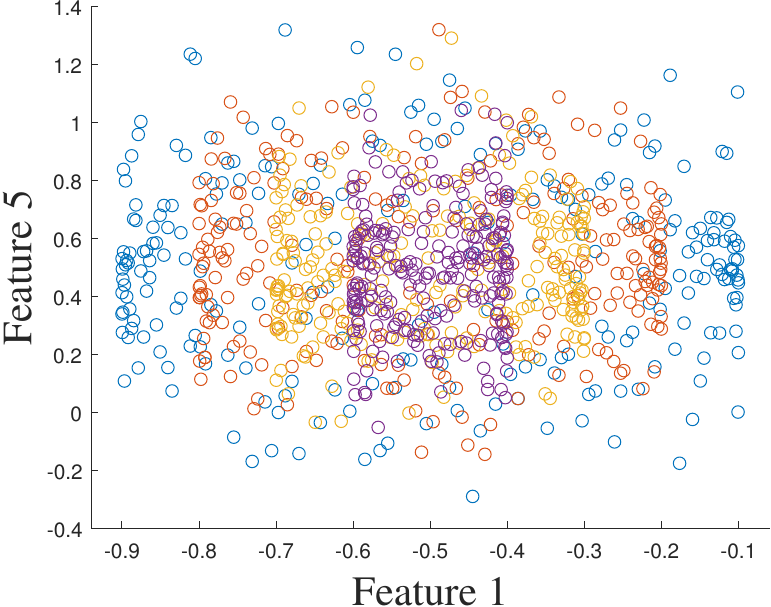}
    }\hspace{-2mm}
    \subfigcapskip=-1pt
    \subfigure[RNE]{
        \label{e}
        \centering
        \includegraphics[width=3.3cm]{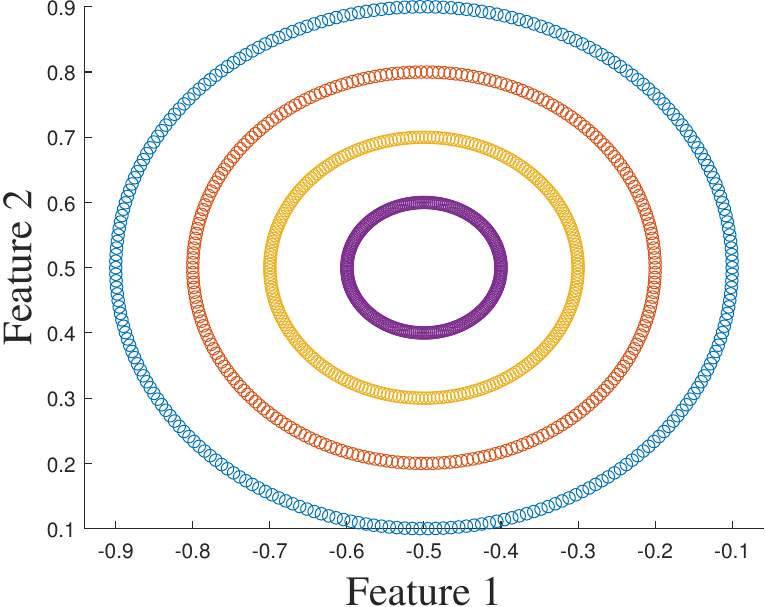}
    }\hspace{-2mm}
    \subfigcapskip=-1pt
    \subfigure[UDFS]{
        \label{f}
        \centering
        \includegraphics[width=3.3cm]{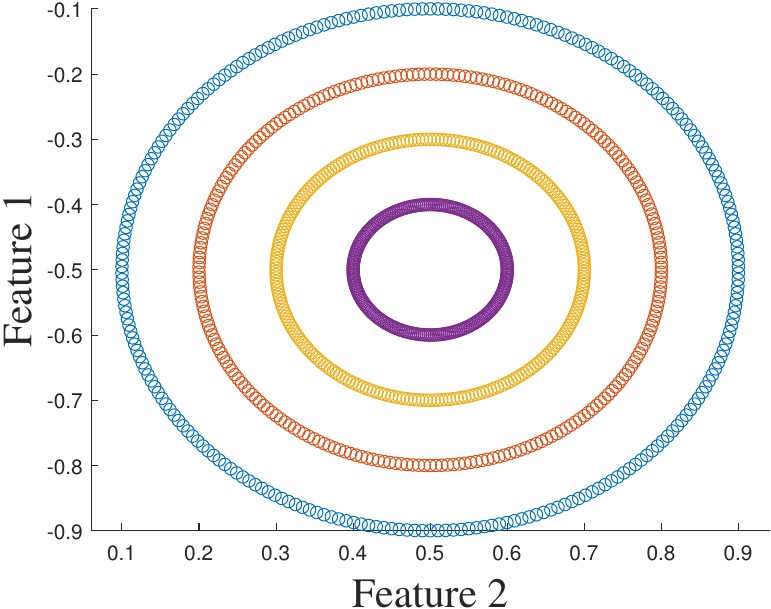}
    }\hspace{-2mm}

    \subfigure[SPCAFS]{
        \label{d}
        \centering
        \includegraphics[width=3.3cm]{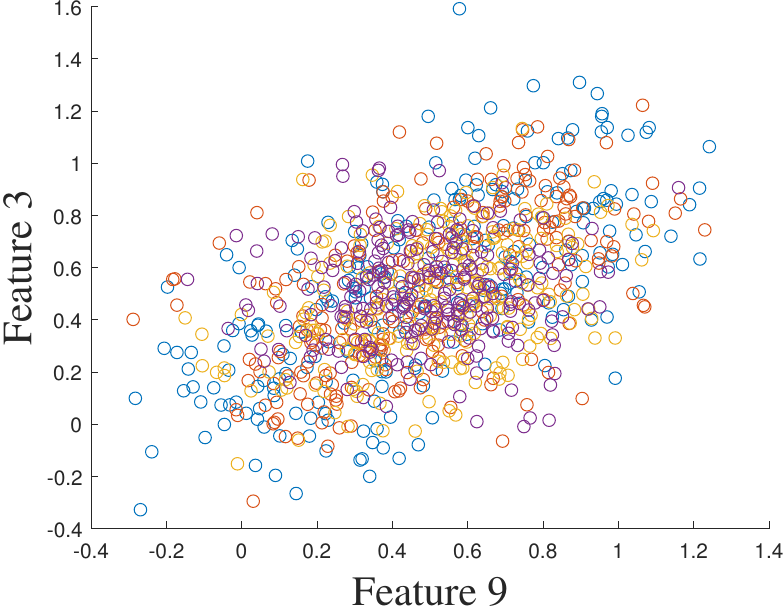}
    }\hspace{-2mm}
    \subfigcapskip=-1pt
    \subfigure[FSPCA]{
        \label{h}
        \centering
        \includegraphics[width=3.3cm]{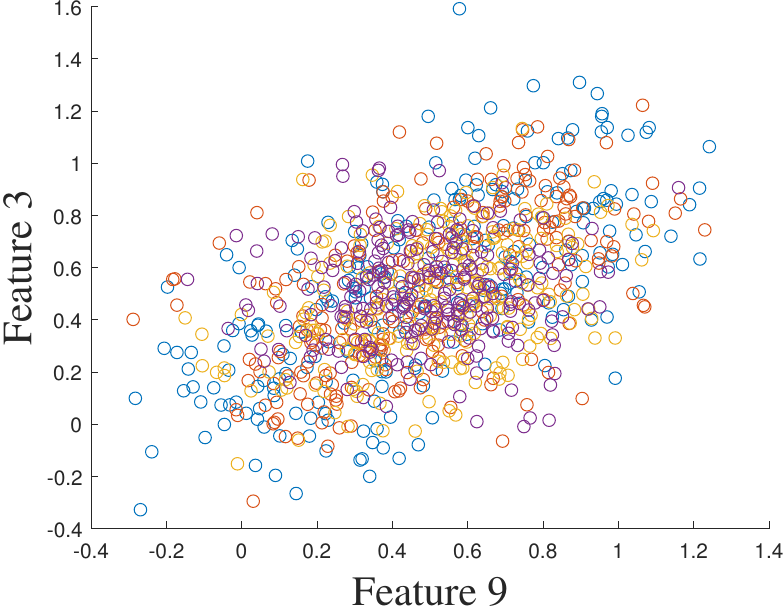}
    }\hspace{-2mm}
     \subfigcapskip=-1pt
    \subfigure[SPCA-PSD]{
        \label{g}
        \centering
        \includegraphics[width=3.3cm]{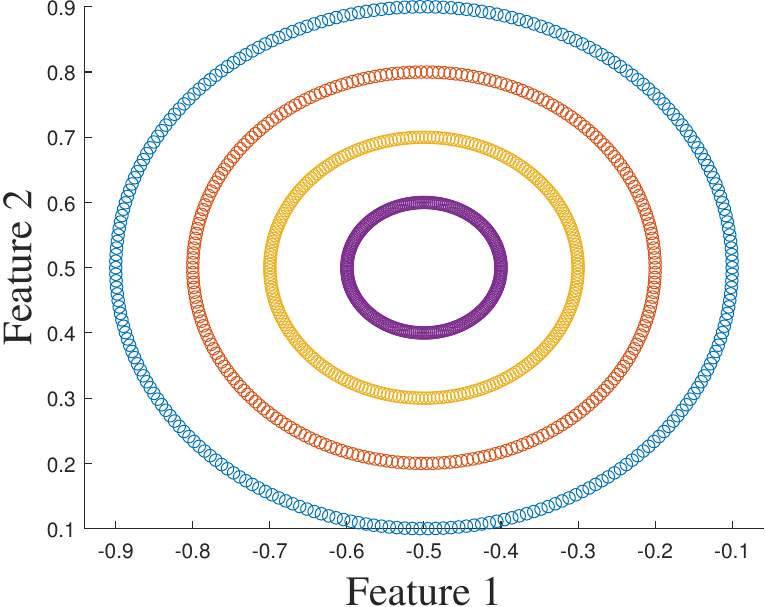}
    }\hspace{-2mm}
    \subfigcapskip=-1pt
    \subfigure[FEN-PCAFS]{
        \label{i}
        \centering
        \includegraphics[width=3.3cm]{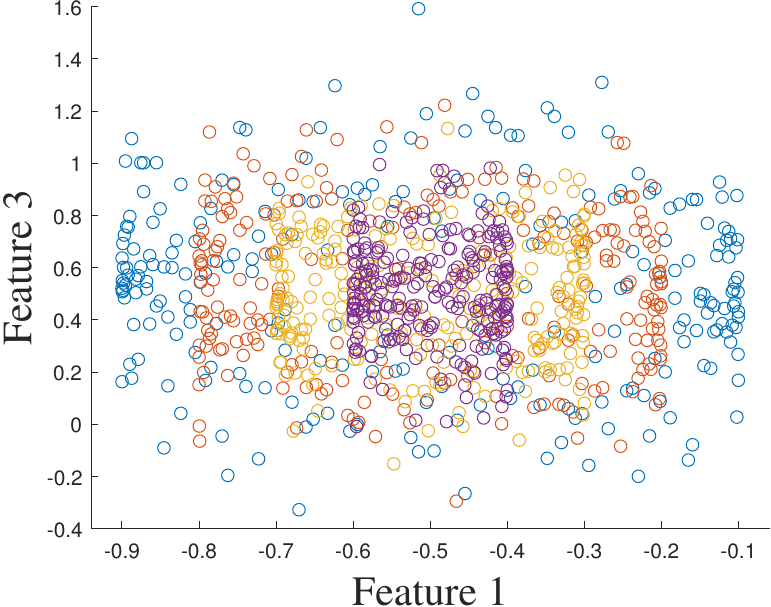}
    }\hspace{-2mm}
    \subfigcapskip=-1pt
    \subfigure[BSUFS]{
        \label{j}
        \centering
        \includegraphics[width=3.3cm]{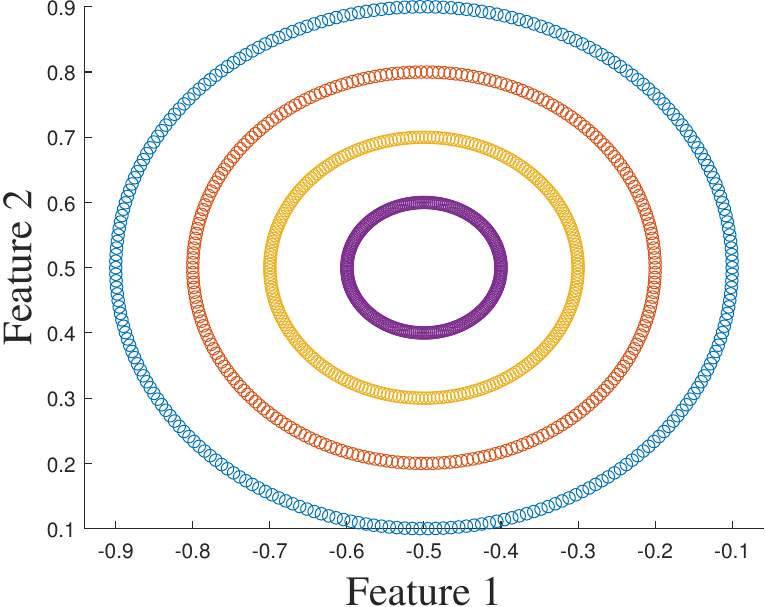}
    }\hspace{-2mm}
    \vspace{-0.05cm}
    \caption{Visual comparisons on the Dartboard1 dataset, where (a) is the dataset distribution and (b)-(j) are the feature selection results.}
     \label{dartboard1}
\end{figure*}

\begin{figure*}[!ht]
   \centering
   \subfigcapskip=-1pt
   \subfigure[Dartboard1]{
       \label{a}
       \centering
       \includegraphics[width=3.3cm]{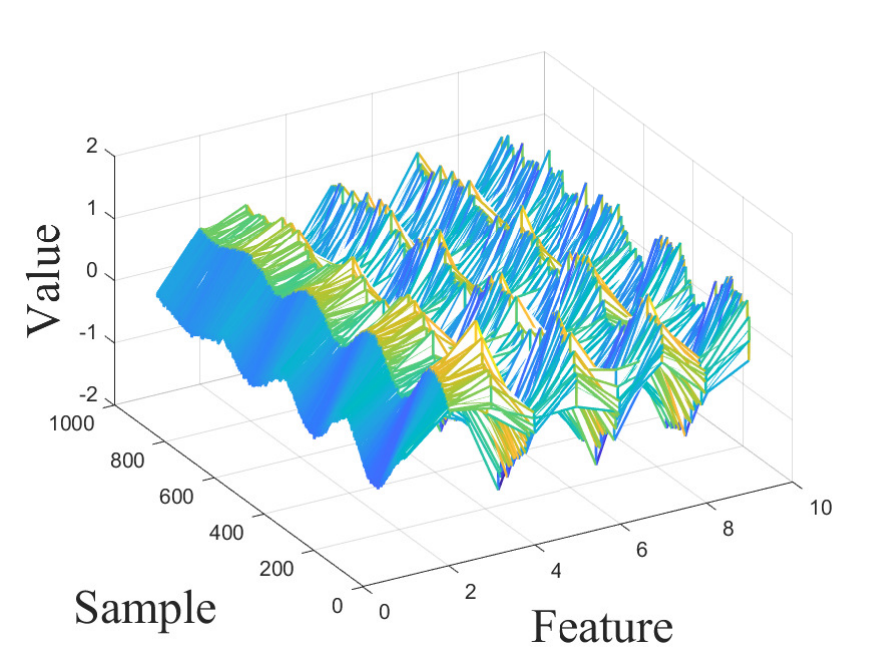}
   }\hspace{-2mm}
   \subfigcapskip=-1pt
   \subfigure[LapScore]{
       \label{b}
       \centering
       \includegraphics[width=3.3cm]{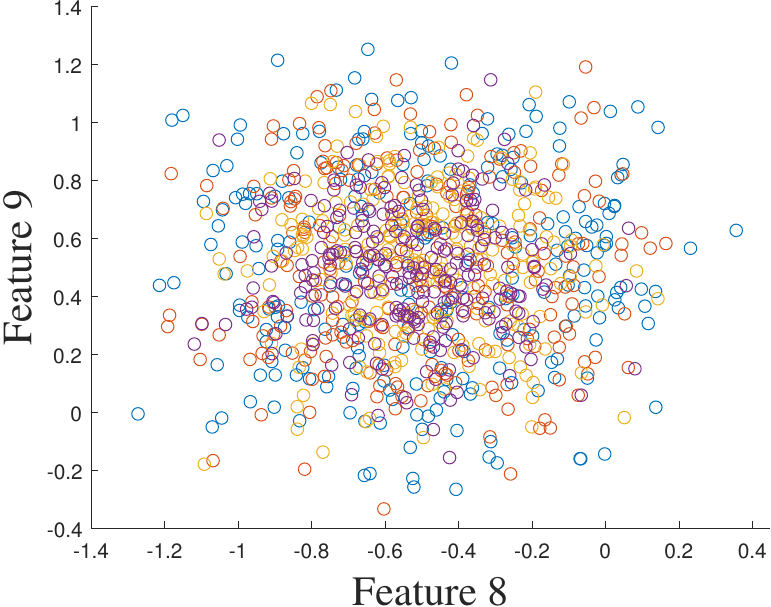}
   }\hspace{-2mm}
   \subfigcapskip=-1pt
   \subfigure[SOGFS]{
       \label{c}
       \centering
       \includegraphics[width=3.3cm]{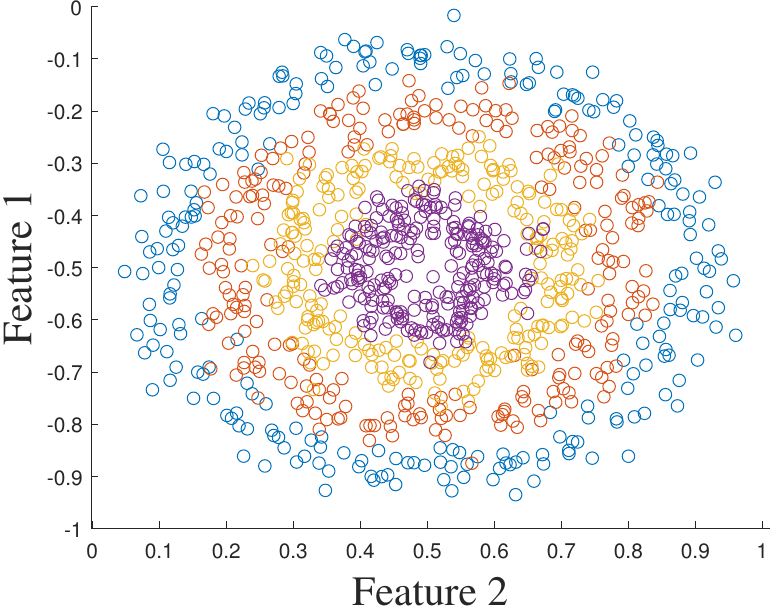}
   }\hspace{-2mm}
   \subfigcapskip=-1pt
   \subfigure[RNE]{
       \label{e}
       \centering
       \includegraphics[width=3.3cm]{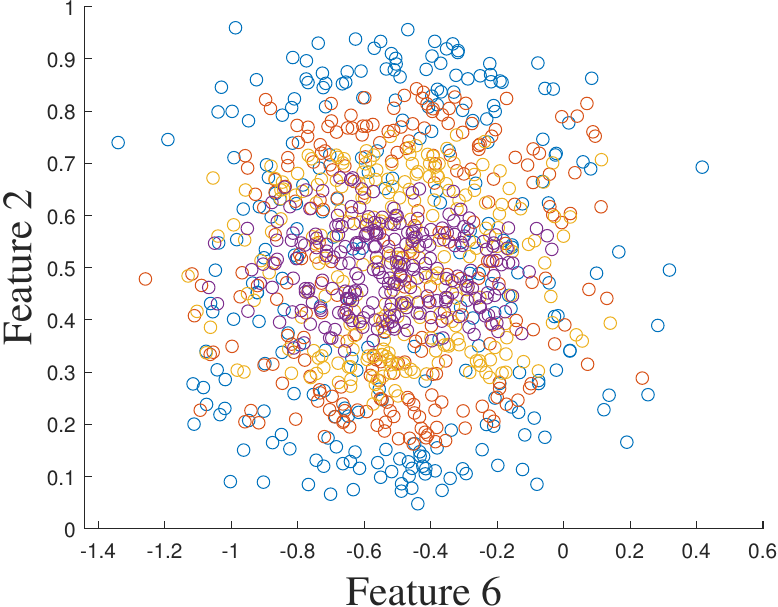}
   }\hspace{-2mm}
   \subfigcapskip=-1pt
   \subfigure[UDFS]{
       \label{f}
       \centering
       \includegraphics[width=3.3cm]{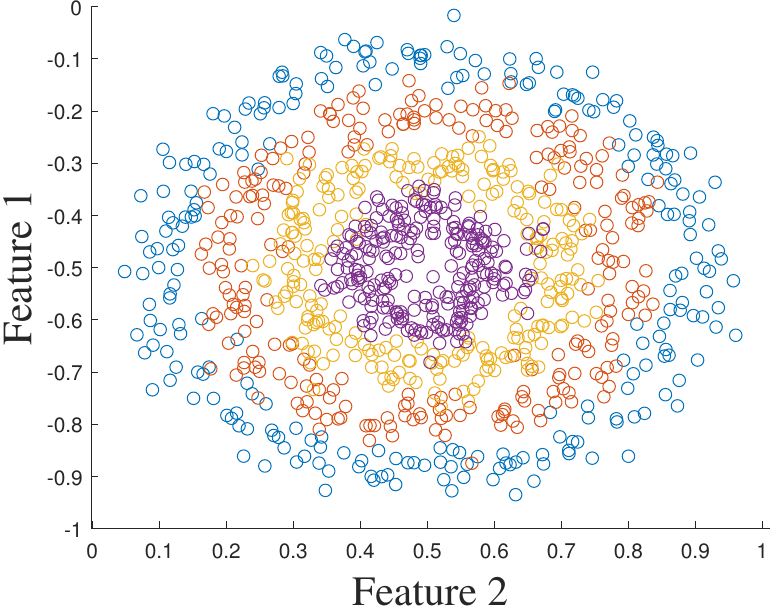}
   }\hspace{-2mm}

   \subfigure[SPCAFS]{
       \label{d}
       \centering
       \includegraphics[width=3.3cm]{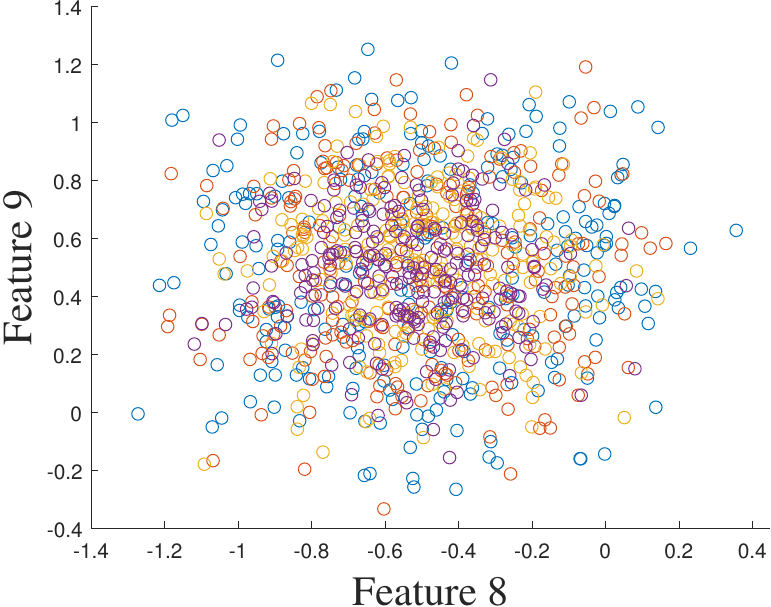}
   }\hspace{-2mm}
   \subfigcapskip=-1pt
   \subfigure[FSPCA]{
       \label{h}
       \centering
       \includegraphics[width=3.3cm]{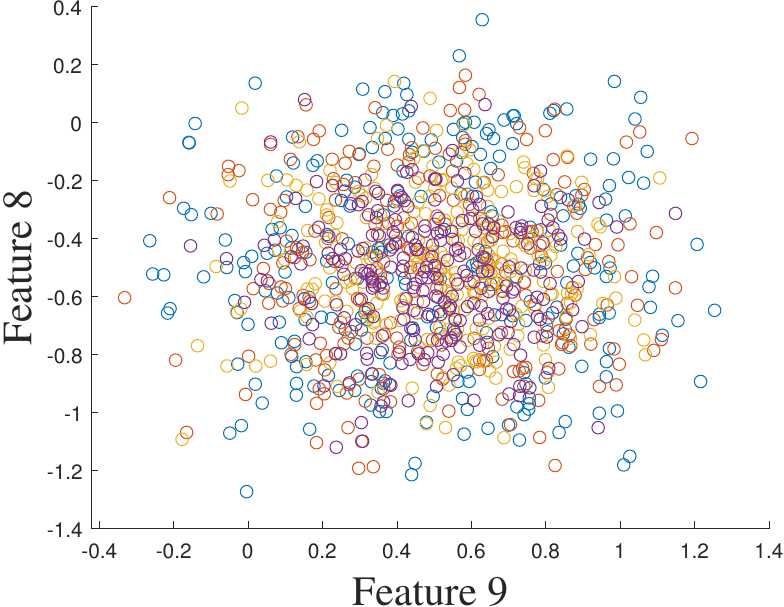}
   }\hspace{-2mm}
    \subfigcapskip=-1pt
   \subfigure[SPCA-PSD]{
       \label{g}
       \centering
       \includegraphics[width=3.3cm]{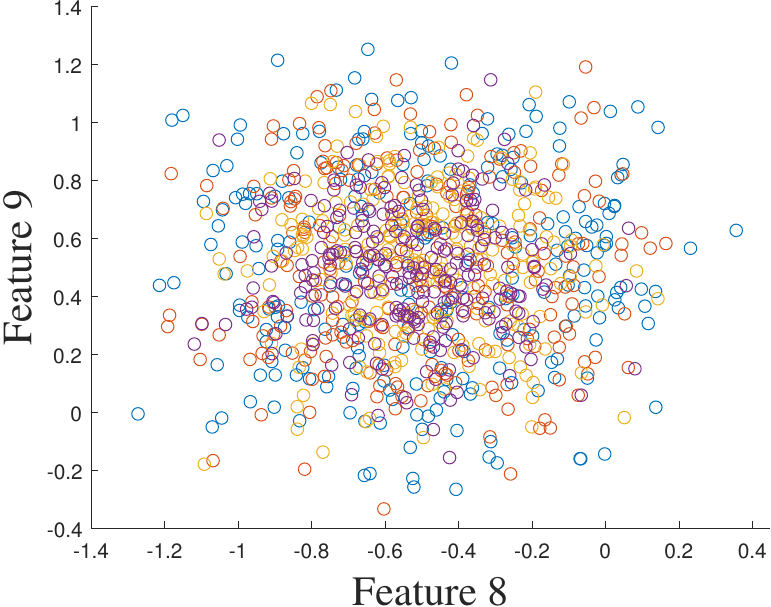}
   }\hspace{-2mm}
   \subfigcapskip=-1pt
   \subfigure[FEN-PCAFS]{
       \label{i}
       \centering
       \includegraphics[width=3.3cm]{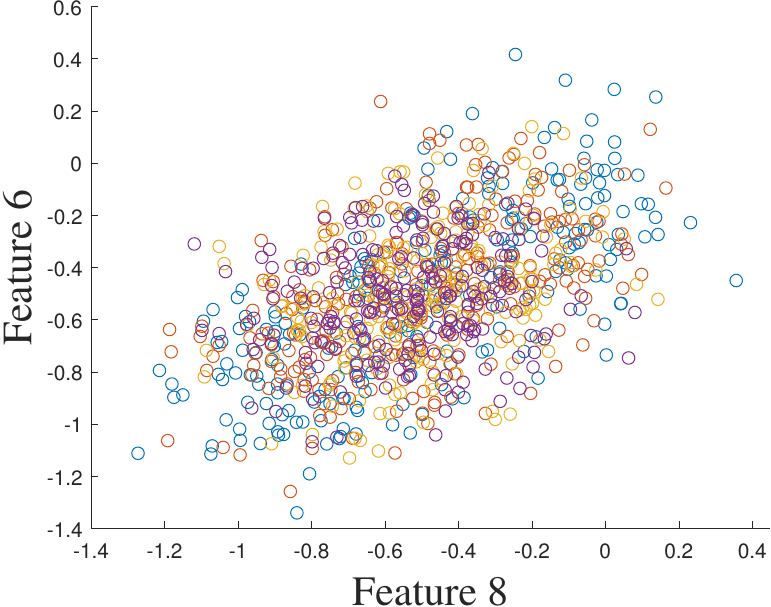}
   }\hspace{-2mm}
   \subfigcapskip=-1pt
   \subfigure[BSUFS]{
       \label{j}
       \centering
       \includegraphics[width=3.3cm]{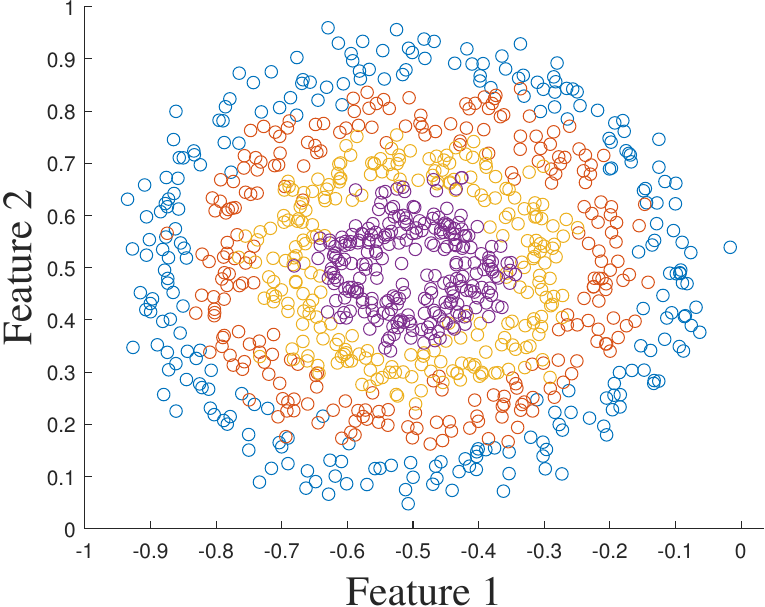}
   }\hspace{-2mm}
   \vspace{-0.05cm}
   \caption{Visual comparisons on the Dartboard1 dataset corrupted by 0.01 Gaussian noise, where (a) is the dataset distribution and (b)-(j) are the feature selection results.}
    \label{dartboard1-gaussian}
\end{figure*}

\begin{figure*}[!ht]
   \centering
   \subfigcapskip=-1pt
   \subfigure[Dartboard1]{
       \label{a}
       \centering
       \includegraphics[width=3.3cm]{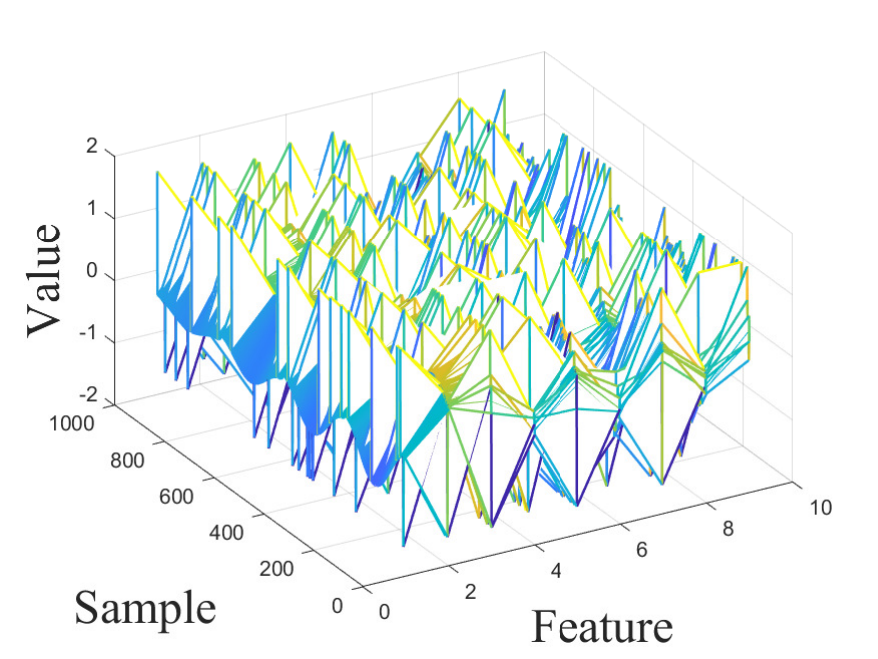}
   }\hspace{-2mm}
   \subfigcapskip=-1pt
   \subfigure[LapScore]{
       \label{b}
       \centering
       \includegraphics[width=3.3cm]{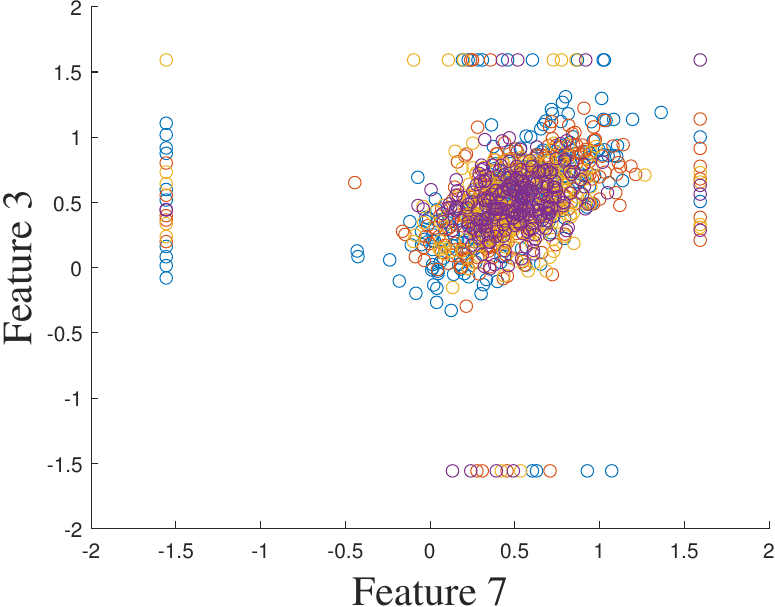}
   }\hspace{-2mm}
   \subfigcapskip=-1pt
   \subfigure[SOGFS]{
       \label{c}
       \centering
       \includegraphics[width=3.3cm]{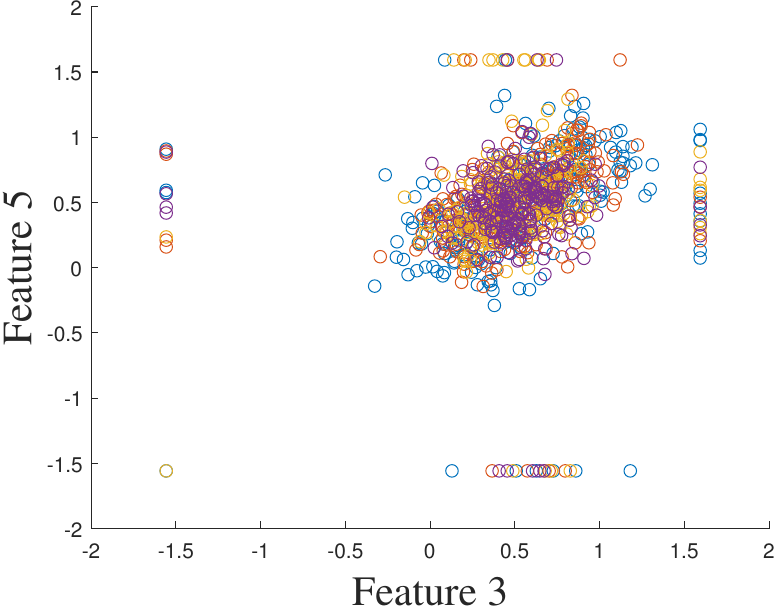}
   }\hspace{-2mm}
   \subfigcapskip=-1pt
   \subfigure[RNE]{
       \label{e}
       \centering
       \includegraphics[width=3.3cm]{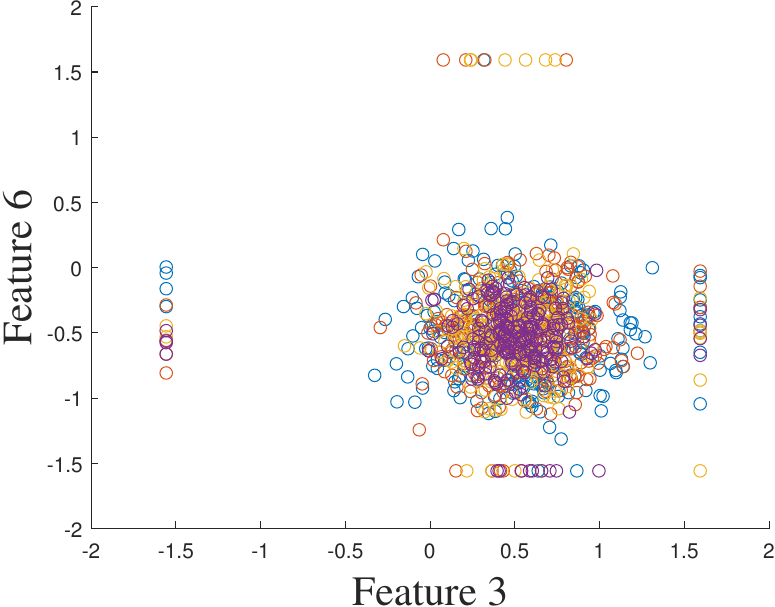}
   }\hspace{-2mm}
   \subfigcapskip=-1pt
   \subfigure[UDFS]{
       \label{f}
       \centering
       \includegraphics[width=3.3cm]{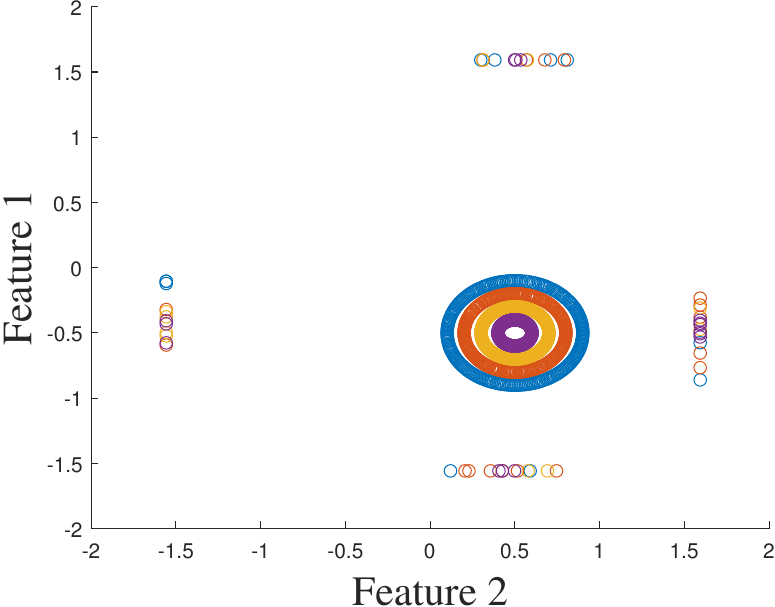}
   }\hspace{-2mm}

   \subfigure[SPCAFS]{
       \label{d}
       \centering
       \includegraphics[width=3.3cm]{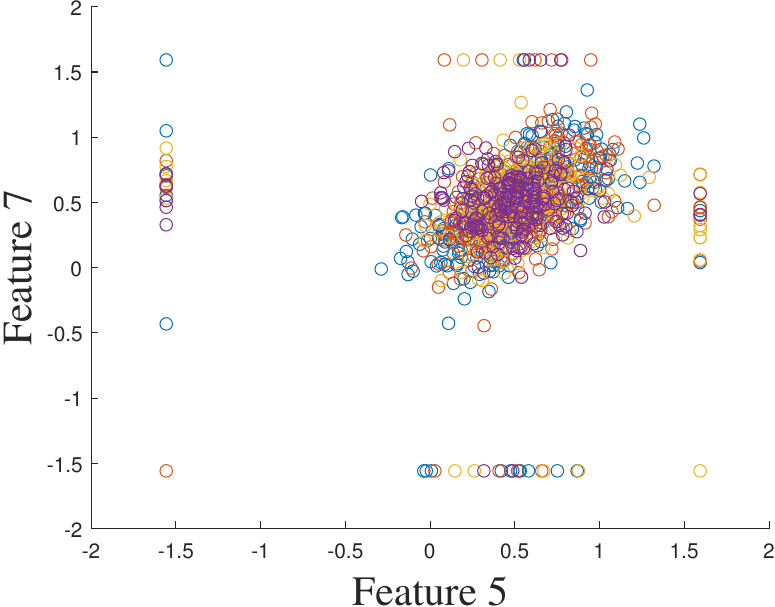}
   }\hspace{-2mm}
   \subfigcapskip=-1pt
   \subfigure[FSPCA]{
       \label{h}
       \centering
       \includegraphics[width=3.3cm]{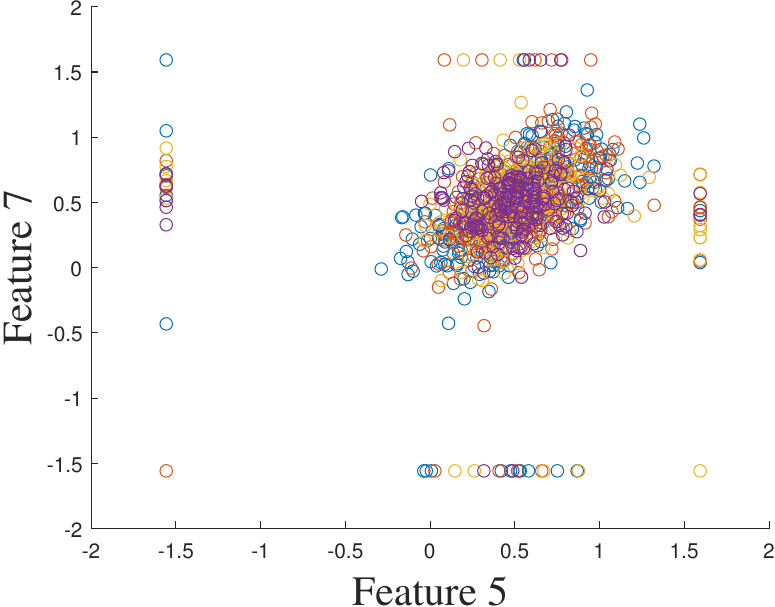}
   }\hspace{-2mm}
    \subfigcapskip=-1pt
   \subfigure[SPCA-PSD]{
       \label{g}
       \centering
       \includegraphics[width=3.3cm]{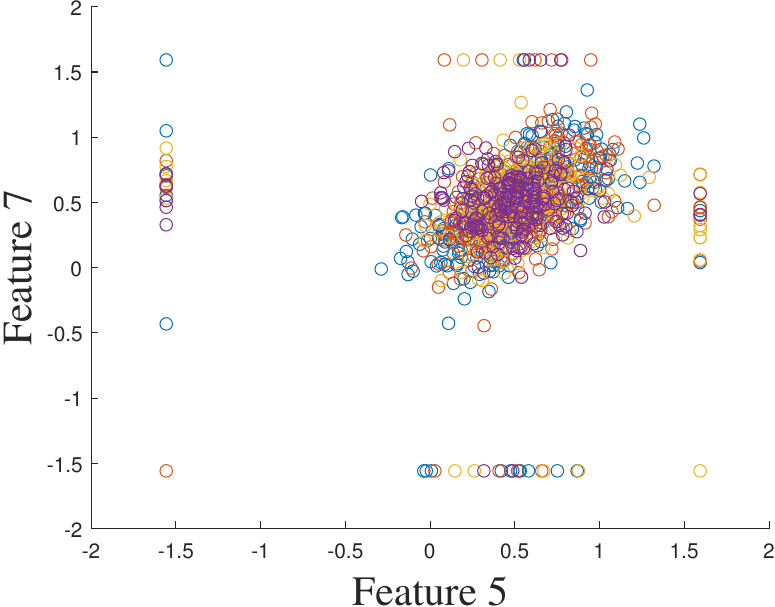}
   }\hspace{-2mm}
   \subfigcapskip=-1pt
   \subfigure[FEN-PCAFS]{
       \label{i}
       \centering
       \includegraphics[width=3.3cm]{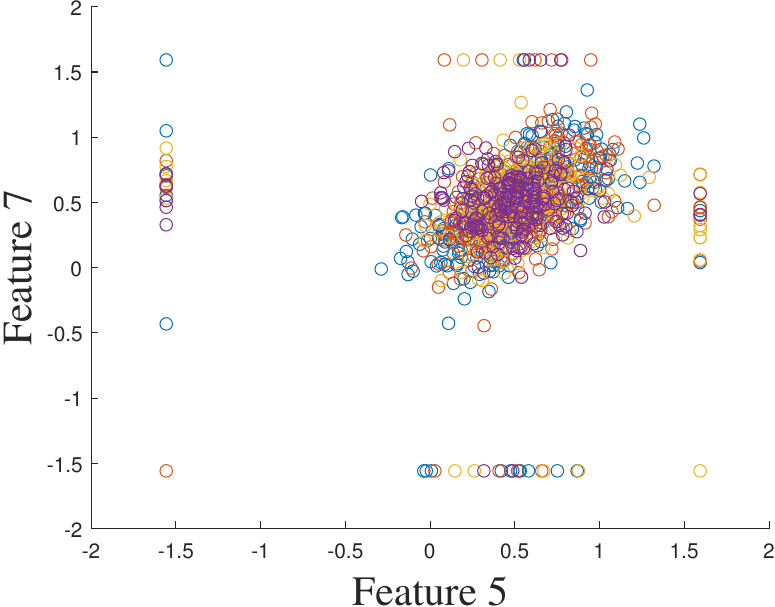}
   }\hspace{-2mm}
   \subfigcapskip=-1pt
   \subfigure[BSUFS]{
       \label{j}
       \centering
       \includegraphics[width=3.3cm]{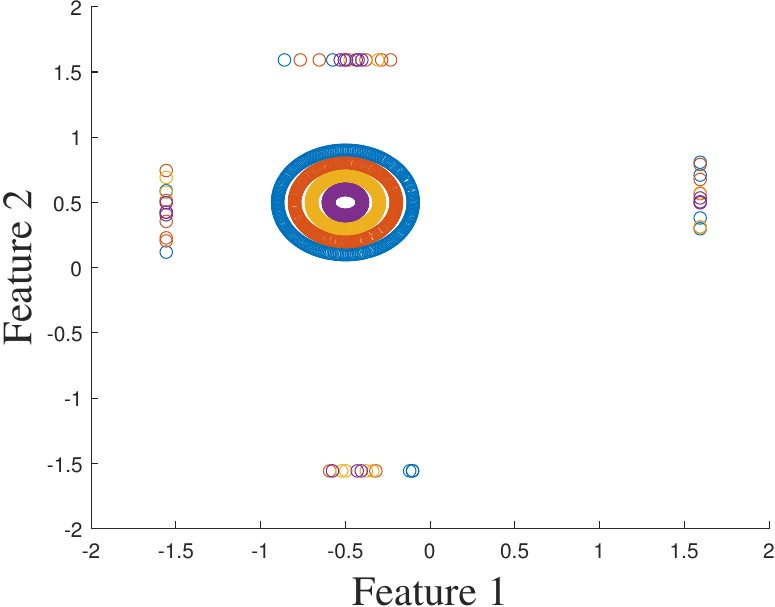}
   }\hspace{-2mm}
   \vspace{-0.05cm}
   \caption{Visual comparisons on the Dartboard1 dataset corrupted by 0.03 salt-and-pepper noise, where (a) is the dataset distribution and (b)-(j) are the feature selection results.}
    \label{dartboard1-sap}
\end{figure*}



\subsection{Experimental Setup}\label{setup}
\subsubsection{Dataset Description}

There are two synthetic datasets and eight real-world datasets that are used to validate the performance of our proposed BSUFS. For more details about these datasets, please refer to Table \ref{data}. 

The two synthetic datasets\footnote{https://github.com/milaan9/Clustering-Datasets} are generated by assigning specific distributions to the first two features, while the remaining seven features are filled with Gaussian noise. 
The eight real-world datasets encompass a diverse range of domains, such as Isolet\footnote{https://jundongl.github.io/scikit-feature/datasets.html \label{web-data}} for spoken letter recognition, MSTAR\_SOC\_CNN\footnote{https://github.com/zjj20212035/SPCA-PSD} (referred to as MSTAR) for deep learning, GLIOMA\textsuperscript{\ref{web-data}} and lung\_discrete\textsuperscript{\ref{web-data}} (referred to as LUNG) for biological information, COIL20\textsuperscript{\ref{web-data}}, USPS\textsuperscript{\ref{web-data}}, pie\footnote{https://data.nvision2.eecs.yorku.ca/PIE\_dataset/}, and umist\footnote{https://github.com/saining/PPSL/blob/master/Platform/Data/UMIST\\/UMIST.mat} for image processing.

\begin{figure*}[t]
    \centering
    \subfigcapskip=-1pt
    \subfigure[COIL20]{
        \label{a}
        \centering
        \includegraphics[width=4.3cm]{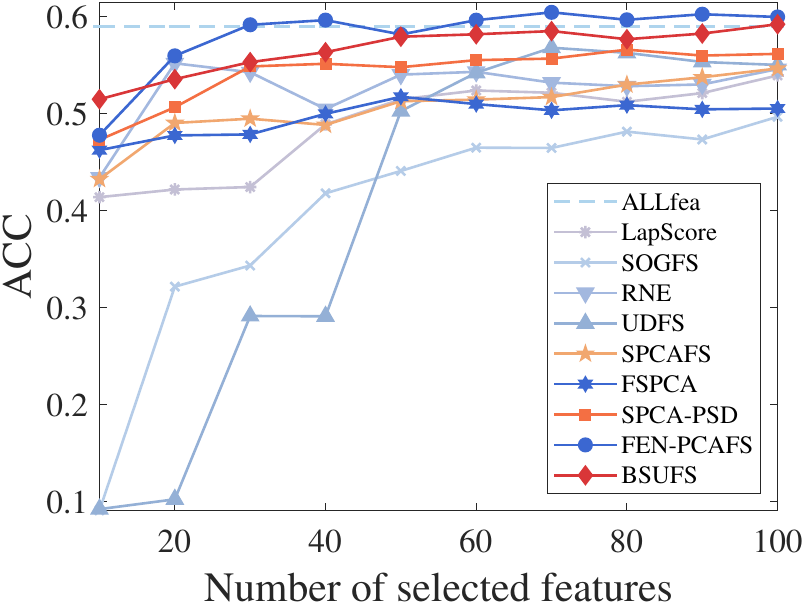}
    }\hspace{-2mm}
    \subfigcapskip=-1pt
    \subfigure[Isolet]{
        \label{b}
        \centering
        \includegraphics[width=4.3cm]{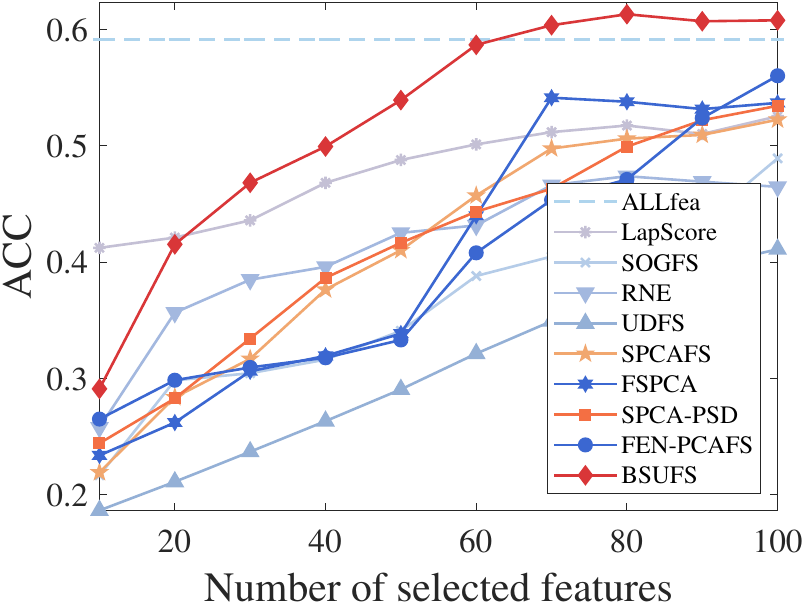}
    }\hspace{-2mm}
    \subfigcapskip=-1pt
    \subfigure[USPS]{
        \label{c}
        \centering
        \includegraphics[width=4.3cm]{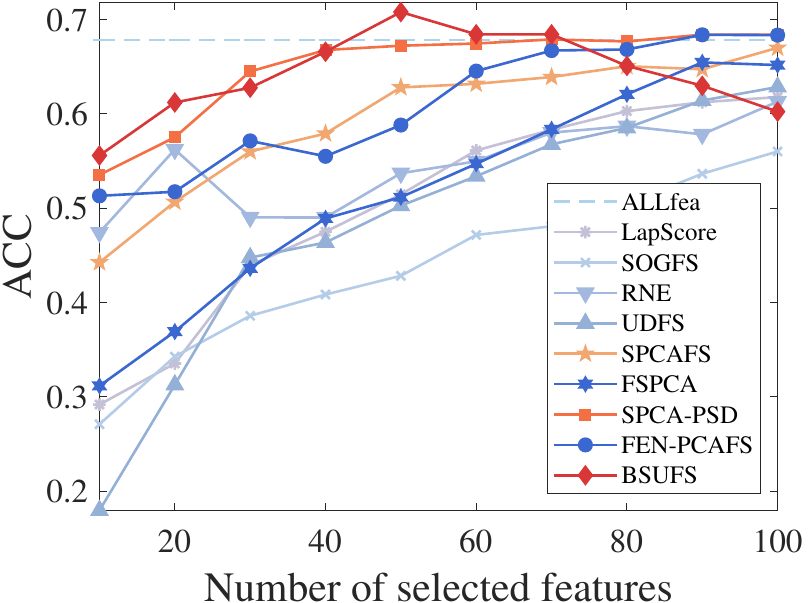}
    }\hspace{-2mm}
    \subfigcapskip=-1pt
    \subfigure[umist]{
        \label{d}
        \centering
        \includegraphics[width=4.3cm]{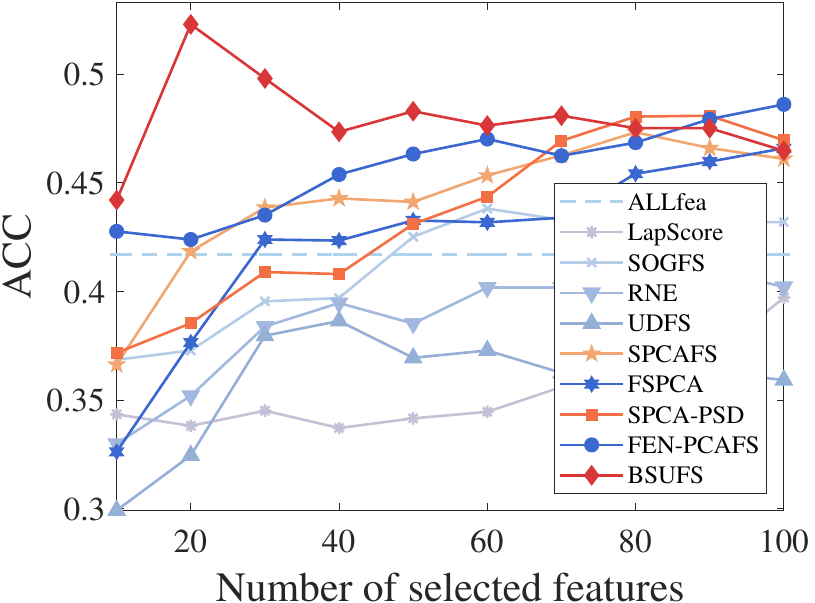}
    }\hspace{-2mm}

    \subfigure[GLIOMA]{
        \label{e}
        \centering
        \includegraphics[width=4.3cm]{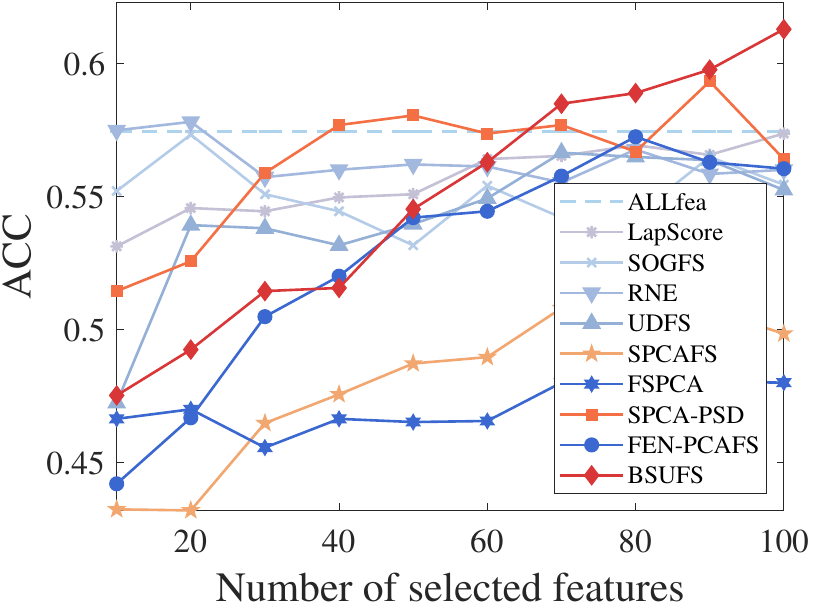}
    }\hspace{-2mm}
    \subfigcapskip=-1pt
    \subfigure[pie]{
        \label{f}
        \centering
        \includegraphics[width=4.3cm]{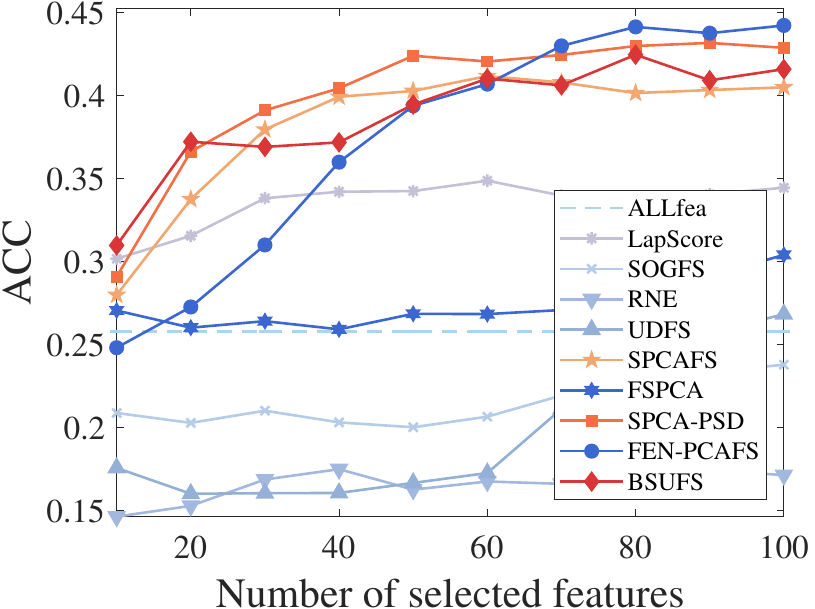}
    }\hspace{-2mm}
    \subfigcapskip=-1pt
    \subfigure[LUNG]{
        \label{g}
        \centering
        \includegraphics[width=4.3cm]{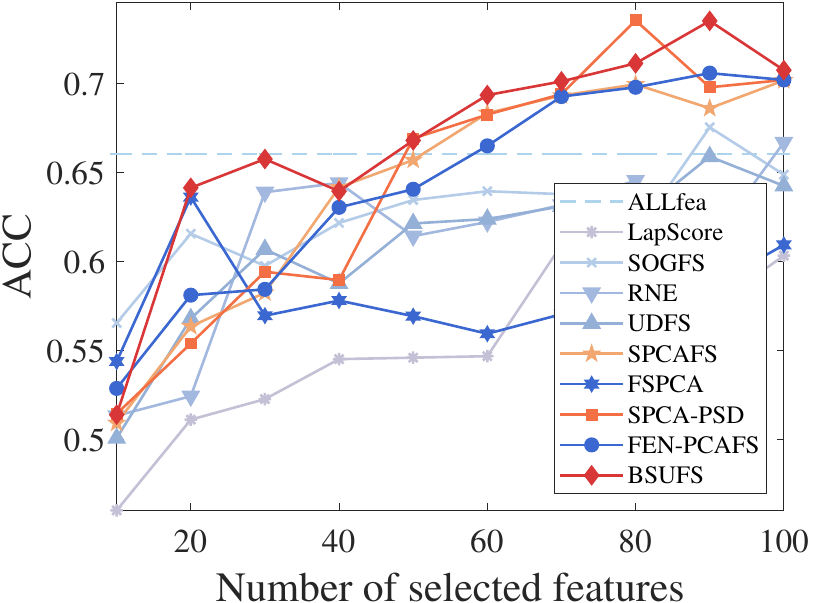}
    }\hspace{-2mm}
    \subfigcapskip=-1pt
    \subfigure[MSTAR]{
        \label{h}
        \centering
        \includegraphics[width=4.3cm]{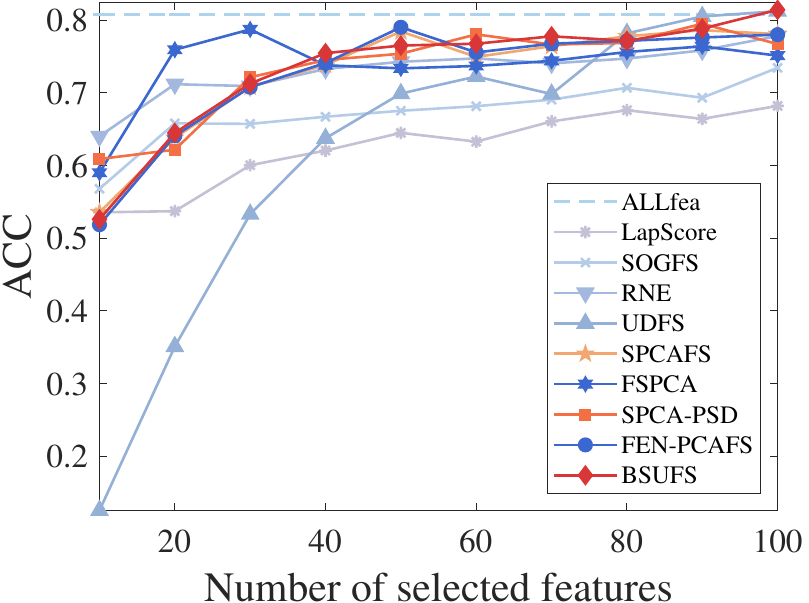}
    }\hspace{-2mm}
    \vspace{-0.1cm}
\caption{Visual comparisons of the ACC metric under different real-world datasets with different number of selected features.} \label{plot-acc}
\end{figure*}

\begin{figure*}[t]
\centering
\subfigcapskip=-1pt
\subfigure[COIL20]{
    \label{a}
    \centering
    \includegraphics[width=4.3cm]{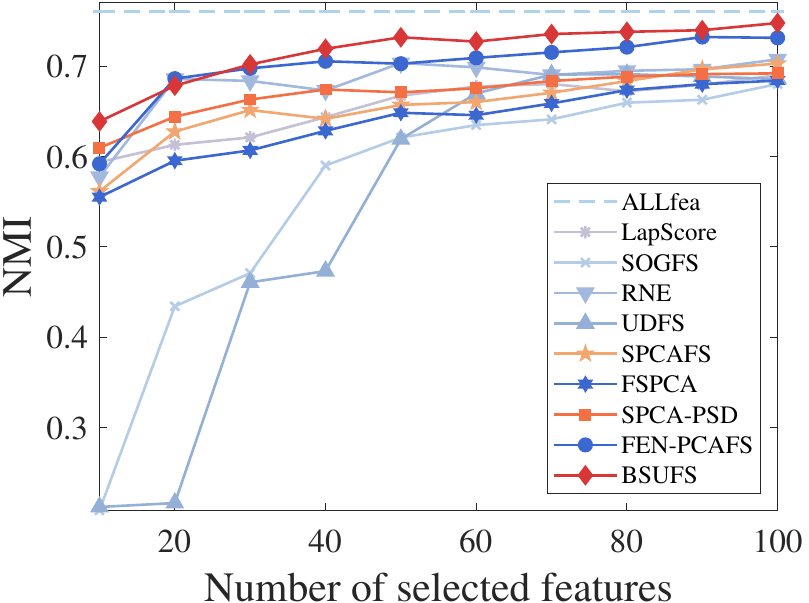}
}\hspace{-2mm}
\subfigcapskip=-1pt
\subfigure[Isolet]{
    \label{b}
    \centering
    \includegraphics[width=4.3cm]{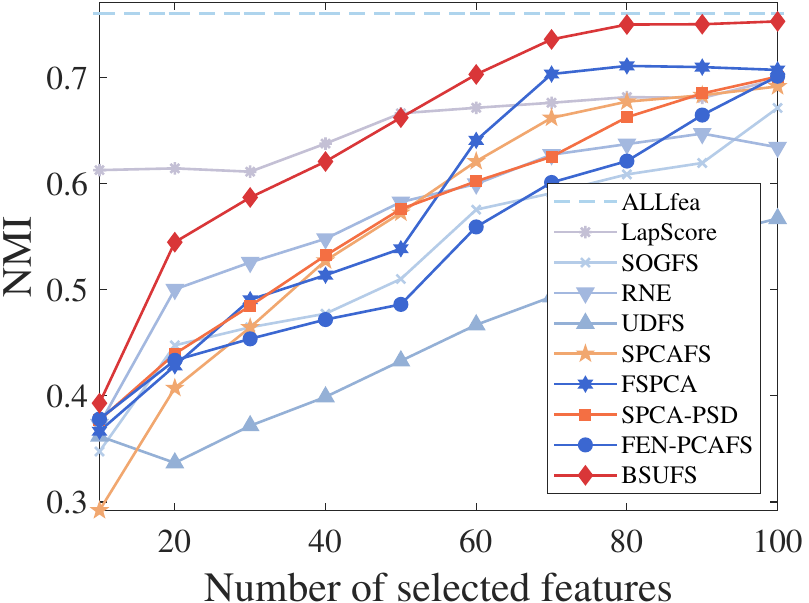}
}\hspace{-2mm}
\subfigcapskip=-1pt
\subfigure[USPS]{
    \label{c}
    \centering
    \includegraphics[width=4.3cm]{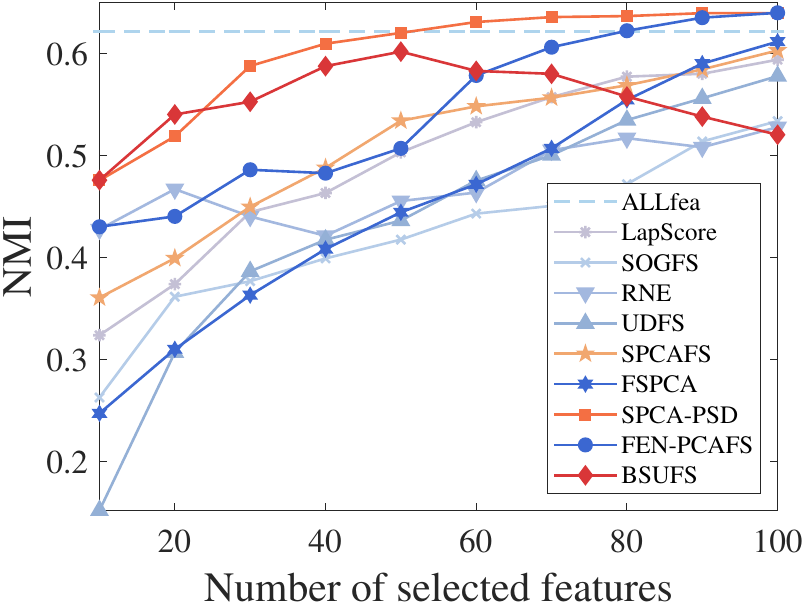}
}\hspace{-2mm}
\subfigcapskip=-1pt
\subfigure[umist]{
    \label{d}
    \centering
    \includegraphics[width=4.3cm]{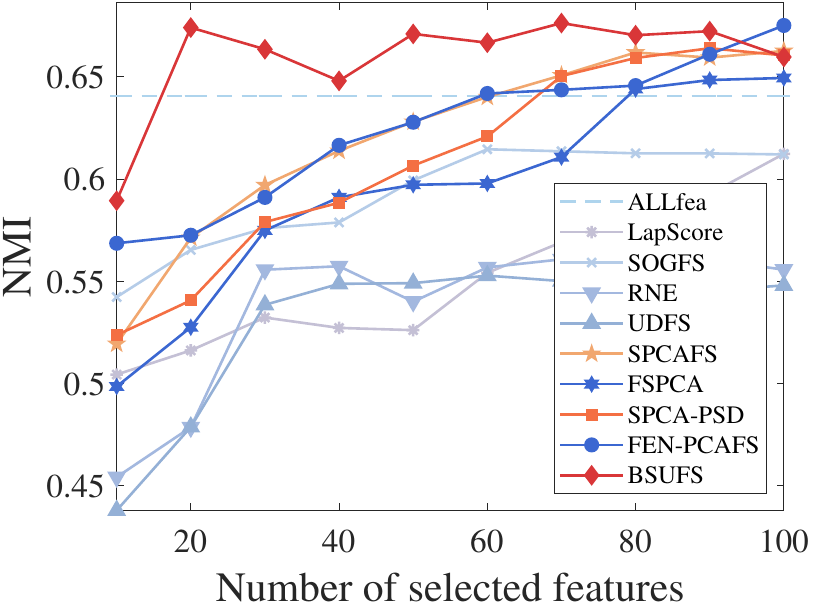}
}\hspace{-2mm}

\subfigure[GLIOMA]{
    \label{e}
    \centering
    \includegraphics[width=4.3cm]{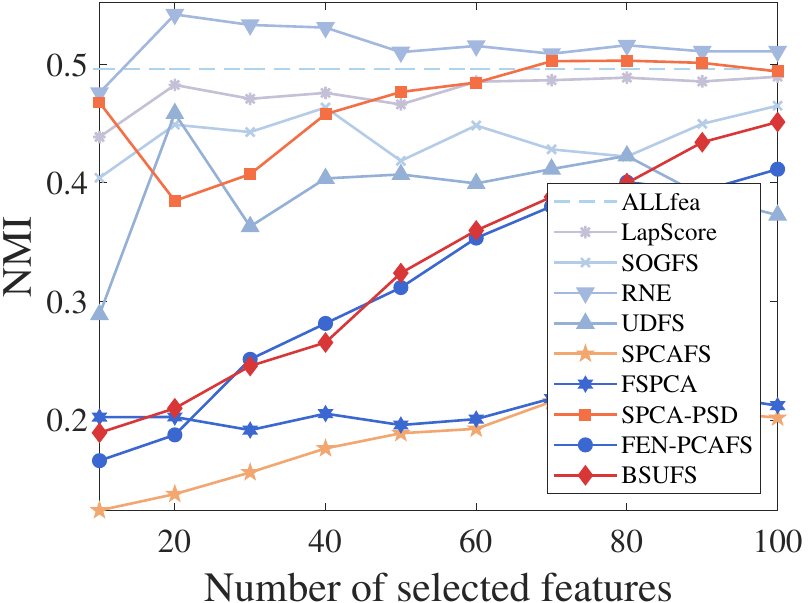}
}\hspace{-2mm}
\subfigcapskip=-1pt
\subfigure[pie]{
    \label{f}
    \centering
    \includegraphics[width=4.3cm]{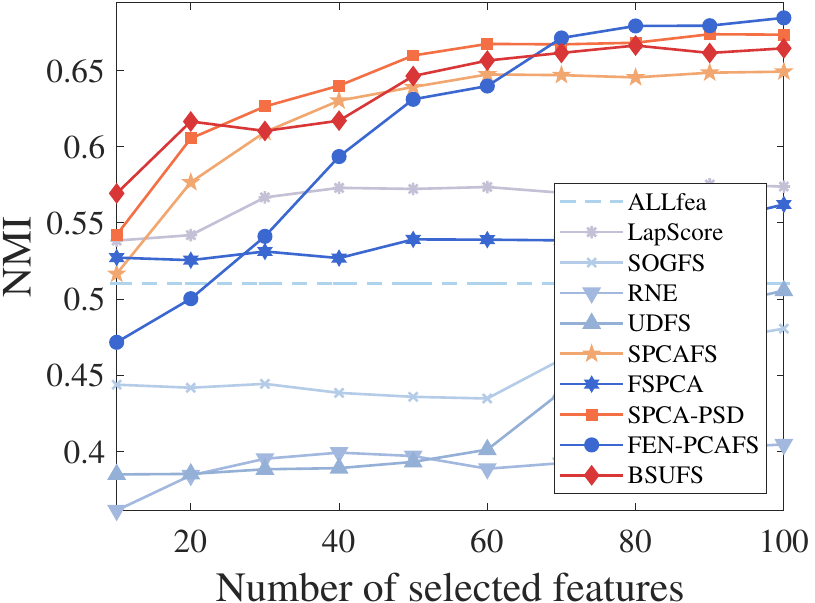}
}\hspace{-2mm}
\subfigcapskip=-1pt
\subfigure[LUNG]{
    \label{g}
    \centering
    \includegraphics[width=4.3cm]{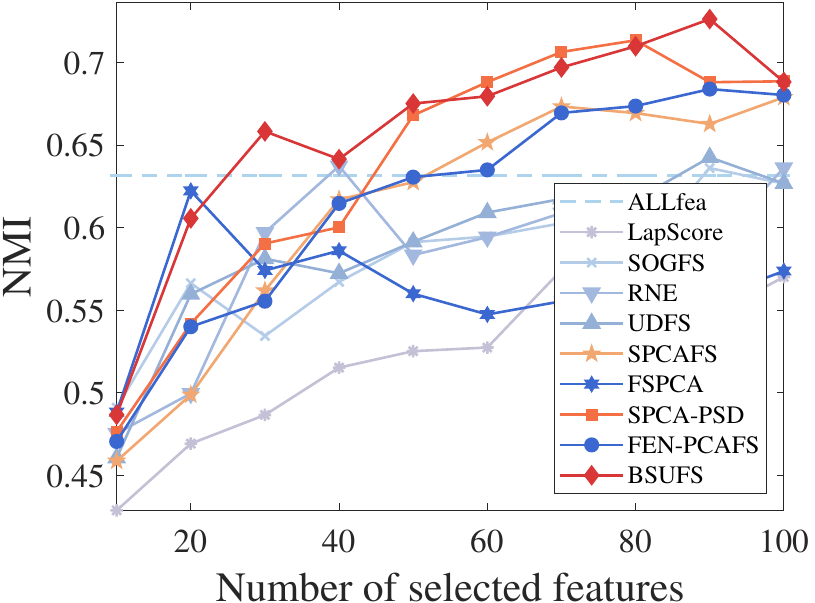}
}\hspace{-2mm}
\subfigcapskip=-1pt
\subfigure[MSTAR]{
    \label{h}
    \centering
    \includegraphics[width=4.3cm]{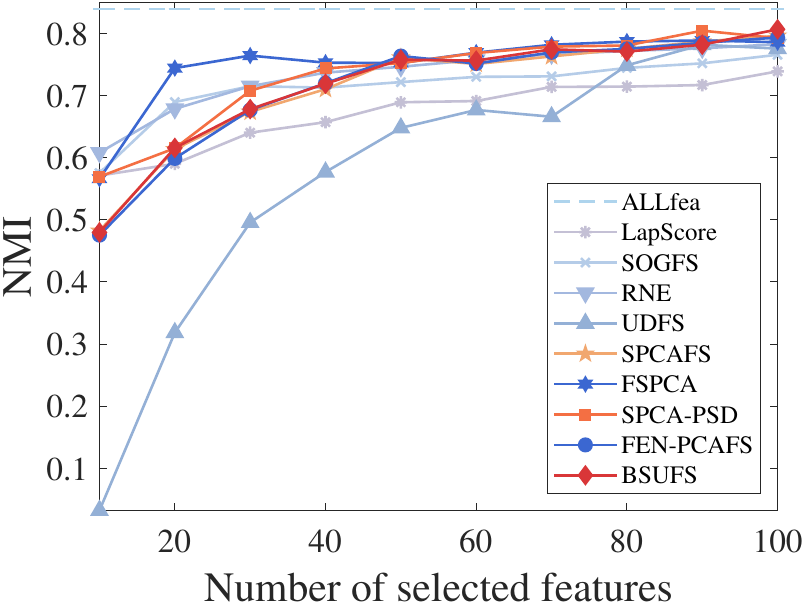}
}\hspace{-2mm}
\vspace{-0.1cm}
\caption{Visual comparisons of the NMI metric under different real-world datasets with different number of selected features.}
\label{plot-nmi}
\end{figure*}
\subsubsection{Parameter Settings}\label{parameter setting}

For LapScore, SOGFS, and RNE, the value of $k$-neighbors is chosen as 5.
For SOGFS, SPCAFS, SPCA-PSD, FEN-PCAFS, and BSUFS, their regularization parameters are selected from the candidate set $\{10^{-6},10^{-4},10^{-2},10^{0},10^{2},10^{4},10^{6}\}$. As suggested in \cite{xu2012l_,zhou2023revisiting}, the values of $p$ and $q$ for BSUFS are selected from $\{0,1/2,2/3\}$. Although other values in $[0, 1)$ are also possible in principle, they are not considered in this paper because there is no closed-form proximal operator and the corresponding solution must be obtained by iterative calculations. For other parameters, the default values or the best parameters provided by the authors are used.

\subsubsection{Evaluation Metrics}

Two key metrics are applied to evaluate these compared methods, including clustering accuracy (ACC) and normalized mutual information (NMI). Here, ACC is defined as
\begin{equation}
\begin{aligned}
    \textrm{ACC} = \frac{1}{n} \sum_{i=1}^{n} \delta(y_{i}, c_{i}),
\end{aligned}
\end{equation}
where $n$ is the number of samples, $y_{i}$ is the true label of the $i$-th sample, and $c_{i}$ is the cluster label of the $i$-th sample. The function $\delta(y_{i}, c_{i})$ is the Kronecker delta function, which equals 1 if $y_{i} = c_{i}$ and 0 otherwise.
Besides, NMI is defined as
\begin{equation}
\begin{aligned}
    \textrm{NMI} = \frac{I(\mathbf{y}, \mathbf{c})}{\sqrt{H(\mathbf{y})H(\mathbf{c})}},
\end{aligned}
\end{equation}
where $\mathbf{y}=(y_1,y_2,\cdots,y_n)\in\mathbb{R}^n$, $\mathbf{c}=(c_1,c_2,\cdots,c_n)\in\mathbb{R}^n$, $I(\mathbf{y}, \mathbf{c})$ is the mutual information between the true label vector $\mathbf{y}$ and the cluster label $\mathbf{c}$, $H(\mathbf{y})$ and $H(\mathbf{c})$ are the entropies of the true label and the cluster label, respectively.

Following a similar way in \cite{li2023sparse}, we select the number of features across all datasets in increments of 10, ranging from 10 to 100. To be fair and counteract the variability introduced by different initial conditions, 50 repetitions of the $k$-means clustering algorithm are conducted. This means that the final results are reported as the mean and standard deviation of the 50 repetitions.

\begin{table*}[h]
    \renewcommand\arraystretch{1.2}
    \caption{Average ACC (mean \% $\pm$ std  \%) and number of selected features  (in brackets) of $K$-means clustering. The top two values are marked as   \textcolor[rgb]{1.00,0.00,0.00}{\textbf{red}} and \textcolor[rgb]{0.00,0.00,1.00}{\textbf{blue}}.}\label{acc}
    \centering%
    \setlength{\tabcolsep}{1.2mm}{
    \begin{tabular}{|c|c|c|c|c|c|c|c|c|c|c|}%
    \hline
     \textbf{Datasets}   & \textbf{ALLfea} & \textbf{LapScore} & \textbf{SOGFS} & \textbf{RNE} & \textbf{UDFS} & \textbf{SPCAFS} & \textbf{FSPCA} & \textbf{SPCA-PSD} & \textbf{FEN-PCAFS} & \textbf{BSUFS} \\ \hline\hline
        \multirow{2}{*}{\textbf{COIL20 }} & 58.97$\pm$4.99 & 53.91$\pm$3.61 & 56.77$\pm$3.09 & 49.66$\pm$3.63 & 55.16$\pm$3.35 & 51.71$\pm$3.05 & 54.63$\pm$3.64 & 56.57$\pm$4.08 & \textcolor[rgb]{1.00,0.00,0.00}{\textbf{60.41$\pm$4.41}} & \textcolor[rgb]{0.00,0.00,1.00}{\textbf{59.18$\pm$3.49}} \\
               & (10) & (100) & (70) & (100) & (20) & (50) & (100) & (80) & \textcolor[rgb]{1.00,0.00,0.00} {(70)} & \textcolor[rgb]{0.00,0.00,1.00}  {(100)} \\ \hline
        \multirow{2}{*}{\textbf{Isolet}} & 59.18$\pm$3.19 & 52.55$\pm$2.83 & 41.11$\pm$1.71 & 48.93$\pm$2.69 & 47.39$\pm$2.91 & 54.15$\pm$2.69 & 52.26$\pm$2.81 & 53.45$\pm$2.82 & \textcolor[rgb]{0.00,0.00,1.00}{\textbf{56.04$\pm$3.50}} & \textcolor[rgb]{1.00,0.00,0.00}{\textbf{61.34$\pm$3.33}} \\
               & (10) & (100) & (100) & (100) & (80) & (70) & (100) & (100) & \textcolor[rgb]{0.00,0.00,1.00} {(100)} & \textcolor[rgb]{1.00,0.00,0.00}  {(80)} \\ \hline
        \multirow{2}{*}{\textbf{USPS}} & 67.79$\pm$4.96 & 61.76$\pm$4.52 & 62.83$\pm$3.79 & 56.00$\pm$3.48 & 61.28$\pm$3.46 & 65.43$\pm$4.90 & 66.98$\pm$3.92 & \textcolor[rgb]{0.00,0.00,1.00}{\textbf{68.38$\pm$3.85}} & 68.36$\pm$4.62 & \textcolor[rgb]{1.00,0.00,0.00}{\textbf{70.77$\pm$3.73}} \\
               & (10) & (100) & (100) & (100) & (100) & (90) & (100) & \textcolor[rgb]{0.00,0.00,1.00}  {(100)} & (90) & \textcolor[rgb]{1.00,0.00,0.00} {(50)} \\ \hline
        \multirow{2}{*}{\textbf{umist}} & 41.68$\pm$2.46 & 39.71$\pm$3.28 & 38.64$\pm$1.61 & 43.81$\pm$2.98 & 41.01$\pm$2.25 & 46.58$\pm$2.34 & 47.32$\pm$3.48 & 48.08$\pm$3.06 & \textcolor[rgb]{0.00,0.00,1.00}{\textbf{48.61$\pm$3.23}} & \textcolor[rgb]{1.00,0.00,0.00}{\textbf{52.29$\pm$3.61}} \\
               & (10) & (100) & (40) & (60) & (90) & (100) & (80) & (90) & \textcolor[rgb]{0.00,0.00,1.00}  {(100)} & \textcolor[rgb]{1.00,0.00,0.00} {(20)} \\ \hline
        \multirow{2}{*}{\textbf{GLIOMA }} & 57.44$\pm$6.40 & 57.36$\pm$3.60 & 56.64$\pm$6.47 & 57.32$\pm$6.47 & 57.80$\pm$2.98 & 48.04$\pm$5.26 & 52.08$\pm$3.64 & \textcolor[rgb]{0.00,0.00,1.00}{\textbf{59.32$\pm$6.27}} & 57.24$\pm$8.16 & \textcolor[rgb]{1.00,0.00,0.00}{\textbf{61.28$\pm$9.01}} \\
               & (10) & (100) & (70) & (20) & (20) & (90) & (80) & \textcolor[rgb]{0.00,0.00,1.00}  {(90)} & (80) & \textcolor[rgb]{1.00,0.00,0.00} {(100)} \\ \hline
        \multirow{2}{*}{\textbf{pie}} & 25.79$\pm$1.39 & 34.86$\pm$1.43 & 26.82$\pm$1.32 & 23.78$\pm$1.19 & 17.49$\pm$0.76 & 30.39$\pm$1.43 & 41.16$\pm$2.46 & \textcolor[rgb]{0.00,0.00,1.00}{\textbf{43.16$\pm$2.38}} & \textcolor[rgb]{1.00,0.00,0.00}{\textbf{44.21$\pm$2.03}} & 42.45$\pm$1.74 \\
               & (10) & (60) & (100) & (100) & (40) & (100) & (60) & \textcolor[rgb]{0.00,0.00,1.00}  {(90)} & \textcolor[rgb]{1.00,0.00,0.00} {(100)} & (80) \\ \hline
        \multirow{2}{*}{\textbf{LUNG}} & 66.03$\pm$7.23 & 60.93$\pm$8.02 & 65.89$\pm$7.43 & 67.53$\pm$7.73 & 66.68$\pm$8.32 & 63.62$\pm$5.45 & 70.16$\pm$7.71 & \textcolor[rgb]{1.00,0.00,0.00}{\textbf{73.53$\pm$8.91}} & 70.58$\pm$6.88 & \textcolor[rgb]{0.00,0.00,1.00}{\textbf{73.51$\pm$6.80}} \\
               & (10) & (70) & (90) & (90) & (100) & (20) & (100) & \textcolor[rgb]{1.00,0.00,0.00} {(80)} & (90) & \textcolor[rgb]{0.00,0.00,1.00}  {(90)} \\ \hline
        \multirow{2}{*}{\textbf{MSTAR}} & 80.81$\pm$8.76 & 68.21$\pm$4.57 & \textcolor[rgb]{0.00,0.00,1.00}{\textbf{81.25$\pm$7.48}} & 73.46$\pm$5.61 & 77.82$\pm$6.16 & 78.74$\pm$5.20 & 78.63$\pm$8.68 & 79.53$\pm$6.75 & 79.03$\pm$6.02 & \textcolor[rgb]{1.00,0.00,0.00}{\textbf{81.43$\pm$6.89}} \\
               & (10) & (100) & \textcolor[rgb]{0.00,0.00,1.00}  {(100)} & (100) & (100) & (30) & (90) & (90) & (50) & \textcolor[rgb]{1.00,0.00,0.00} {(100)} \\ \hline \hline
        \multirow{1}{*}{\textbf{Average}} & 57.21$\pm$4.92 & 53.66$\pm$3.98 & 53.74$\pm$4.11 & 52.56$\pm$4.22 & 53.08$\pm$3.77 & 54.83$\pm$3.79 & 57.90$\pm$4.54 & 60.25$\pm$4.76 & \textcolor[rgb]{0.00,0.00,1.00}{\textbf{60.56$\pm$4.86}} & \textcolor[rgb]{1.00,0.00,0.00}{\textbf{62.78$\pm$4.83}}  \\ \hline  
    \end{tabular}
 } 
\end{table*}

\begin{table*}[h]
    \renewcommand\arraystretch{1.2}
    \caption{Average NMI (mean \% $\pm$ std \%) and number of selected features  (in brackets) of $K$-means clustering. The top two values are marked as   \textcolor[rgb]{1.00,0.00,0.00}{\textbf{red}} and \textcolor[rgb]{0.00,0.00,1.00}{\textbf{blue}}.} \label{nmi}
    \centering%
    \setlength{\tabcolsep}{1.2mm}{
    \begin{tabular}{|c|c|c|c|c|c|c|c|c|c|c|}%
    \hline
     \textbf{Datasets}   & \textbf{ALLfea} & \textbf{LapScore} & \textbf{SOGFS} & \textbf{RNE} & \textbf{UDFS} & \textbf{SPCAFS} & \textbf{FSPCA} & \textbf{SPCA-PSD} & \textbf{FEN-PCAFS} & \textbf{BSUFS} \\ \hline\hline
        \multirow{2}{*}{\textbf{COIL20 }} & 76.04$\pm$1.69 & 69.01$\pm$1.53 & 69.12$\pm$1.17 & 68.03$\pm$1.59 & 70.76$\pm$2.07 & 68.41$\pm$1.60 & 70.29$\pm$1.31 & 69.21$\pm$1.41 & \textcolor[rgb]{0.00,0.00,1.00}{\textbf{73.23$\pm$1.31}}& \textcolor[rgb]{1.00,0.00,0.00}{\textbf{74.78$\pm$1.79}} \\
               & (10) & (100) & (80) & (100) & (100) & (100) & (100) & (100) & \textcolor[rgb]{0.00,0.00,1.00}  {(90)} & \textcolor[rgb]{1.00,0.00,0.00} {(100)} \\ \hline
        \multirow{2}{*}{\textbf{Isolet}} & 76.09$\pm$1.77 & 69.86$\pm$1.26 & 56.73$\pm$1.05 & 67.15$\pm$1.45 & 64.74$\pm$1.28 & \textcolor[rgb]{0.00,0.00,1.00}{\textbf{71.12$\pm$1.11}}& 69.18$\pm$1.33 & 70.11$\pm$1.11 & 70.14$\pm$1.56 & \textcolor[rgb]{1.00,0.00,0.00}{\textbf{75.32$\pm$1.22}} \\
               & (10) & (100) & (100) & (100) & (90) & \textcolor[rgb]{0.00,0.00,1.00}  {(80)} & (100) & (100) & (100) & \textcolor[rgb]{1.00,0.00,0.00} {(100)} \\ \hline
        \multirow{2}{*}{\textbf{USPS}} & 62.11$\pm$2.24 & 59.37$\pm$1.98 & 57.76$\pm$2.02 & 53.36$\pm$1.83 & 52.77$\pm$2.01 & 61.14$\pm$1.87 & 60.28$\pm$2.17 & \textcolor[rgb]{0.00,0.00,1.00}{\textbf{63.93$\pm$2.06}}& \textcolor[rgb]{1.00,0.00,0.00}{\textbf{63.96$\pm$2.24}} & 60.16$\pm$1.68 \\
               & (10) & (100) & (100) & (100) & (100) & (100) & (100) & \textcolor[rgb]{0.00,0.00,1.00}  {(90)} & \textcolor[rgb]{1.00,0.00,0.00} {(100)} & (50) \\ \hline
        \multirow{2}{*}{\textbf{umist}} & 64.07$\pm$1.76 & 61.23$\pm$2.15 & 55.43$\pm$1.50 & 61.46$\pm$2.03 & 56.08$\pm$1.80 & 64.94$\pm$1.65 & 66.26$\pm$1.74 & 66.39$\pm$1.93 & \textcolor[rgb]{0.00,0.00,1.00}{\textbf{67.51$\pm$1.92}}& \textcolor[rgb]{1.00,0.00,0.00}{\textbf{67.62$\pm$1.91}} \\
               & (10) & (100) & (80) & (60) & (70) & (100) & (100) & (90) & \textcolor[rgb]{0.00,0.00,1.00}  {(100)} & \textcolor[rgb]{1.00,0.00,0.00} {(70)} \\ \hline
        \multirow{2}{*}{\textbf{GLIOMA }} & 49.59$\pm$6.76 & 48.96$\pm$3.59 & 45.86$\pm$8.08 & 46.51$\pm$9.11 & \textcolor[rgb]{1.00,0.00,0.00}{\textbf{54.21$\pm$2.23}} & 22.17$\pm$5.17 & 22.01$\pm$4.88 & \textcolor[rgb]{0.00,0.00,1.00}{\textbf{50.31$\pm$6.65}}& 41.16$\pm$7.66 & 45.14$\pm$8.66 \\
               & (10) & (100) & (20) & (100) & \textcolor[rgb]{1.00,0.00,0.00} {(20)} & (90) & (80) & \textcolor[rgb]{0.00,0.00,1.00}  {(80)} & (100) & (100) \\ \hline
        \multirow{2}{*}{\textbf{pie}} & 51.01$\pm$1.02 & 57.53$\pm$0.73 & 50.55$\pm$1.03 & 48.05$\pm$0.76 & 40.45$\pm$0.79 & 56.21$\pm$0.90 & 64.94$\pm$1.30 & \textcolor[rgb]{0.00,0.00,1.00}{\textbf{67.40$\pm$1.21}}& \textcolor[rgb]{1.00,0.00,0.00}{\textbf{68.47$\pm$1.15}} & 66.66$\pm$1.14 \\
               & (10) & (90) & (100) & (100) & (100) & (100) & (100) & \textcolor[rgb]{0.00,0.00,1.00}  {(90)} & \textcolor[rgb]{1.00,0.00,0.00} {(100)} & (80) \\ \hline
        \multirow{2}{*}{\textbf{LUNG}} & 63.18$\pm$5.48 & 57.44$\pm$6.44 & 64.27$\pm$5.35 & 63.62$\pm$5.41 & 63.74$\pm$5.30 & 62.23$\pm$4.80 & 67.91$\pm$6.23 & \textcolor[rgb]{0.00,0.00,1.00}{\textbf{71.36$\pm$6.71}}& 68.40$\pm$5.34 & \textcolor[rgb]{1.00,0.00,0.00}{\textbf{72.64$\pm$4.69}} \\
               & (10) & (70) & (90) & (90) & (40) & (20) & (100) & \textcolor[rgb]{0.00,0.00,1.00}  {(80)} & (90) & \textcolor[rgb]{1.00,0.00,0.00} {(90)} \\ \hline
        \multirow{2}{*}{\textbf{MSTAR}} & 83.96$\pm$3.14 & 73.90$\pm$1.62 & 78.18$\pm$3.64 & 76.56$\pm$1.54 & 78.26$\pm$2.51 & 78.87$\pm$2.52 & 79.62$\pm$2.30 & \textcolor[rgb]{0.00,0.00,1.00}{\textbf{80.44$\pm$2.04}}& 79.34$\pm$3.27 & \textcolor[rgb]{1.00,0.00,0.00}{\textbf{80.66$\pm$2.68}} \\
               & (10) & (100) & (90) & (100) & (100) & (90) & (100) & \textcolor[rgb]{0.00,0.00,1.00}  {(90)} & (100) & \textcolor[rgb]{1.00,0.00,0.00} {(100)} \\ \hline \hline
        \multirow{1}{*}{\textbf{Average}} & 65.76$\pm$2.98 & 62.16$\pm$2.41 & 59.74$\pm$2.98 & 60.59$\pm$2.96 & 60.13$\pm$2.25 & 60.64$\pm$2.45 & 62.56$\pm$2.66 & \textcolor[rgb]{0.00,0.00,1.00}{\textbf{67.39$\pm$2.89}}& 66.53$\pm$3.06 & \textcolor[rgb]{1.00,0.00,0.00}{\textbf{67.87$\pm$2.97}}  \\ \hline
    \end{tabular}
 } 
\end{table*}

\subsection{Synthetic Results}\label{experiments}

In this experiment, various UFS methods are conducted on two synthetic datasets, i.e., Diamond9 and Dartboard1. Each UFS method is applied to obtain scores for all nine features and the top two features are selected. After that, the selected features are visualized in scatter diagrams alongside the samples.

Fig. \ref{diamond9} shows the feature selection results of the Diamond9 dataset, where (a) is the dataset distribution and (b)-(j) are the feature selection results. It is demonstrated that BSUFS is the only method that can successfully identify the two most discriminative features.

For the Dartboard1 dataset, Fig. \ref{dartboard1} shows that RNE, UDFS, SPCA-PSD, and BSUFS are all capable of identifying the appropriate features, compared to the other methods. In Fig. \ref{dartboard1-gaussian}, by adding  Gaussian noise with mean 0 and standard deviation 0.01 on this dataset, UDFS, SPCA-PSD, and BSUFS identify the right features, but RNE fails. As for Fig. \ref{dartboard1-sap}, the 0.03 salt-and-pepper noise is added that is 3\% of the pixels are  affected by this noise. From this figure, only UDFS and BSUFS are capable of selecting the right features. Therefore, it can be concluded that UDFS and BSUFS perform better and are more robust on the Dartboard1 dataset. Besides, between  graph-based methods and PCA-based methods, it is difficult to say which ones perform better.

To sum up, our proposed BSUFS consistently selects discriminative features and performs robustly on different noises, while the other methods fail to select the correct features in some cases. All those results suggest the effectiveness of our proposed BSUFS on synthetic datasets.

\begin{table}[t]
    \renewcommand\arraystretch{1.2} 
    \caption{ACC (\%) comparisons for four cases in ablation experiments. The top two values are marked as   \textcolor[rgb]{1.00,0.00,0.00}{\textbf{red}} and \textcolor[rgb]{0.00,0.00,1.00}{\textbf{blue}}.}\label{a-acc} 
    \centering
    \setlength{\tabcolsep}{1.5mm} 
    \begin{tabular}{|c|c|c|c|c|}
      \hline
     ~{\textbf{Datasets}}~ & ~\textbf{Case I}~ & ~{\textbf{Case II}}~ & ~{\textbf{Case III}}~ &~{\textbf{Case IV}}~ \\ 
      \hline \hline
      \textbf{COIL20} & 54.09   &  57.33   & \textcolor[rgb]{0.00,0.00,1.00} {\textbf{58.76}}  & \textcolor[rgb]{1.00,0.00,0.00}{\textbf{59.18}}  \\ \hline
      \textbf{Isolet} & 51.77   &  58.78   & \textcolor[rgb]{0.00,0.00,1.00} {\textbf{56.19}}  & \textcolor[rgb]{1.00,0.00,0.00}{\textbf{61.34}} \\ \hline
      \textbf{USPS}   & 67.06   &  66.42   & \textcolor[rgb]{0.00,0.00,1.00} {\textbf{68.11}}  & \textcolor[rgb]{1.00,0.00,0.00}{\textbf{70.77}}  \\ \hline
      \textbf{umist}  & 47.16   &  47.57   & \textcolor[rgb]{0.00,0.00,1.00} {\textbf{49.23}}  & \textcolor[rgb]{1.00,0.00,0.00}{\textbf{52.29}}  \\ \hline
      \textbf{GLIOMA} & 49.76   &  58.40   & \textcolor[rgb]{0.00,0.00,1.00} {\textbf{60.12}}  & \textcolor[rgb]{1.00,0.00,0.00}{\textbf{61.28}}  \\ \hline
      \textbf{pie}    & 40.98   &  41.05   & \textcolor[rgb]{0.00,0.00,1.00} {\textbf{41.15}}  & \textcolor[rgb]{1.00,0.00,0.00}{\textbf{42.45}}  \\ \hline
      \textbf{LUNG}   & 71.34   &  71.51   & \textcolor[rgb]{0.00,0.00,1.00} {\textbf{72.33}}  & \textcolor[rgb]{1.00,0.00,0.00}{\textbf{73.51}}  \\ \hline
      \textbf{MSTAR}  & 79.25   &  74.67   & \textcolor[rgb]{0.00,0.00,1.00} {\textbf{80.08}}  & \textcolor[rgb]{1.00,0.00,0.00}{\textbf{81.43}}  \\
      \hline
    \end{tabular}
\end{table}

\begin{table}[t]
    \renewcommand\arraystretch{1.2} 
    \caption{NMI  (\%) comparisons for four cases in ablation experiments. The top two values are marked as   \textcolor[rgb]{1.00,0.00,0.00}{\textbf{red}} and \textcolor[rgb]{0.00,0.00,1.00}{\textbf{blue}}.}\label{ff-nmi} 
    \centering
    \setlength{\tabcolsep}{1.5mm} 
    \begin{tabular}{|c|c|c|c|c|}
      \hline
     ~{\textbf{Datasets}}~ & ~\textbf{Case I}~ & ~{\textbf{Case II}}~ & ~{\textbf{Case III}}~ &~{\textbf{Case IV}}~ \\ 
      \hline \hline
      \textbf{COIL20} & 69.94   & 72.12    & \textcolor[rgb]{0.00,0.00,1.00} {\textbf{74.57}}  & \textcolor[rgb]{1.00,0.00,0.00}{\textbf{74.78}}  \\ \hline
      \textbf{Isolet} & 66.84   & 73.30    & \textcolor[rgb]{0.00,0.00,1.00} {\textbf{72.73}}  & \textcolor[rgb]{1.00,0.00,0.00}{\textbf{75.32}}  \\ \hline
      \textbf{USPS}   & \textcolor[rgb]{0.00,0.00,1.00}{\textbf{58.86}}   & 35.66    & \textcolor[rgb]{1.00,0.00,0.00} {\textbf{61.14}}  & 60.16  \\ \hline
      \textbf{umist}  & 66.48   & 67.77    & \textcolor[rgb]{1.00,0.00,0.00} {\textbf{69.45}}  & \textcolor[rgb]{0.00,0.00,1.00}{\textbf{67.62}}  \\ \hline
      \textbf{GLIOMA} & 20.64   & \textcolor[rgb]{1.00,0.00,0.00} {\textbf{52.11}}    & 43.13  & \textcolor[rgb]{0.00,0.00,1.00}{\textbf{45.14}}  \\ \hline
      \textbf{pie}    & 65.02   & 65.13    & \textcolor[rgb]{0.00,0.00,1.00} {\textbf{65.23}}  & \textcolor[rgb]{1.00,0.00,0.00}{\textbf{66.66}}  \\ \hline
      \textbf{LUNG}   & 69.31   &69.17    & \textcolor[rgb]{0.00,0.00,1.00} {\textbf{71.94}}  & \textcolor[rgb]{1.00,0.00,0.00}{\textbf{72.64}}  \\ \hline
      \textbf{MSTAR}  & 79.92   & 73.14    & \textcolor[rgb]{0.00,0.00,1.00} {\textbf{79.97}}  & \textcolor[rgb]{1.00,0.00,0.00}{\textbf{80.66}}  \\
      \hline
    \end{tabular}
\end{table}

\begin{table}[t]
\centering
\caption{Visual comparisons of selected image samples from the pie dataset with the corresponding ACC (\%) and NMI (\%). \\The top two  values are marked as   \textcolor[rgb]{1.00,0.00,0.00}{\textbf{red}} and \textcolor[rgb]{0.00,0.00,1.00}{\textbf{blue}}.}\label{graph}
\begin{tblr}{
  cells={valign=m,halign=c},
  colspec={p{1cm}|p{8mm}|p{8mm}|p{8mm}|p{8mm}|p{7mm}|p{7mm}}, 
  hlines,
  vlines
}
  \textbf{Methods} & \SetCell[c=4]{c} \textbf{Samples} &&&& \textbf{ACC} & \textbf{NMI} \\
  \hline
\textbf{Case I} & \includegraphics[scale=0.1,valign=c]{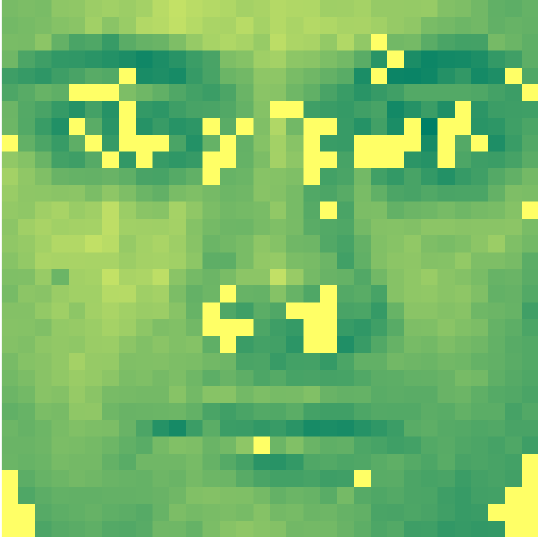} & \includegraphics[scale=0.1,valign=c]{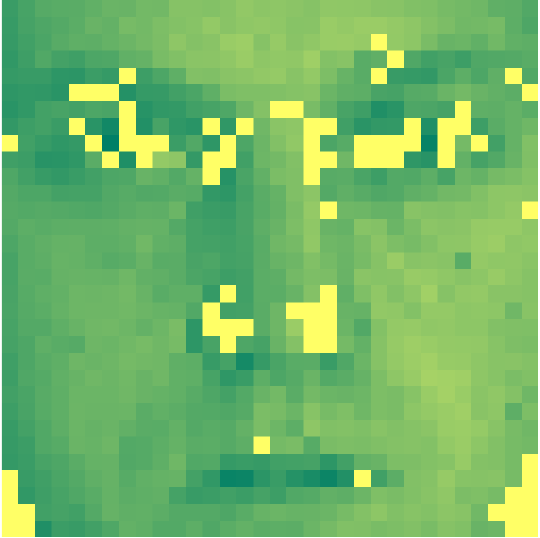} & \includegraphics[scale=0.1,valign=c]{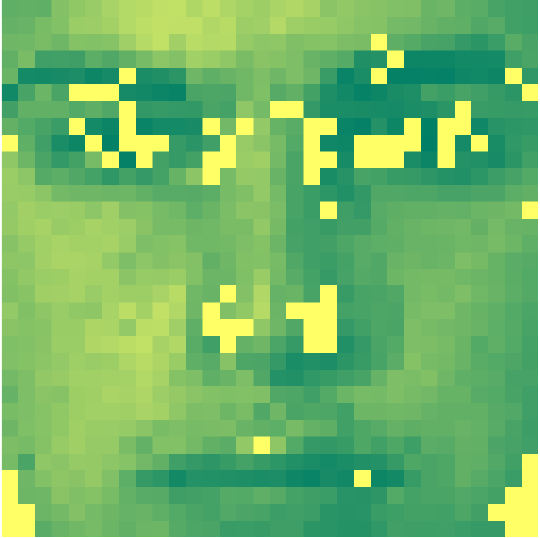} & \includegraphics[scale=0.1,valign=c]{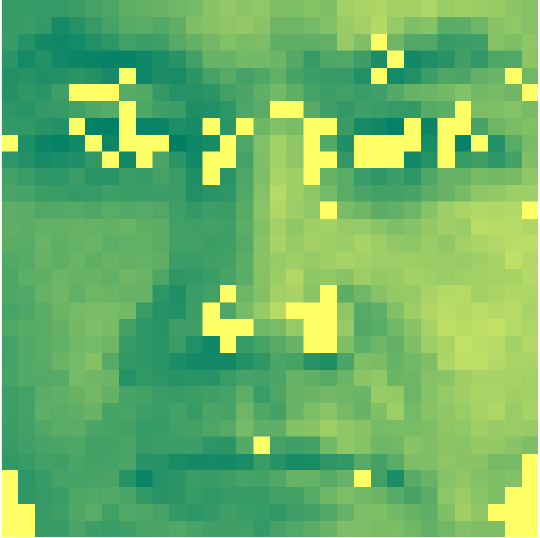} & 40.98 & 65.02\\
\textbf{Case II} & \includegraphics[scale=0.1,valign=c]{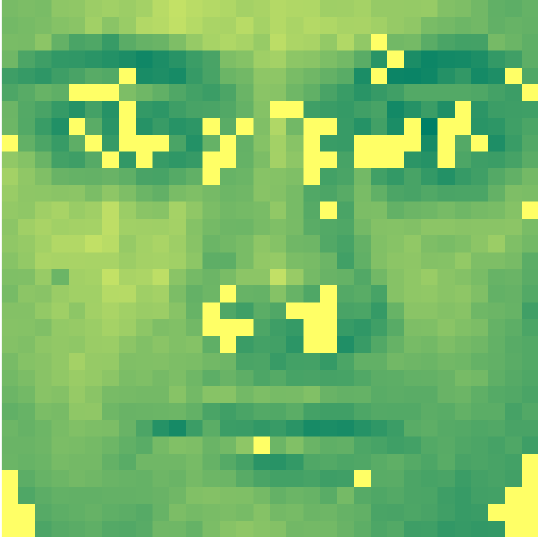} & \includegraphics[scale=0.1,valign=c]{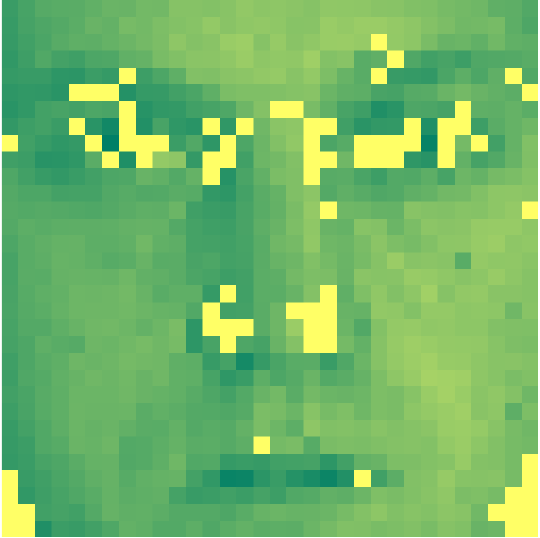} & \includegraphics[scale=0.1,valign=c]{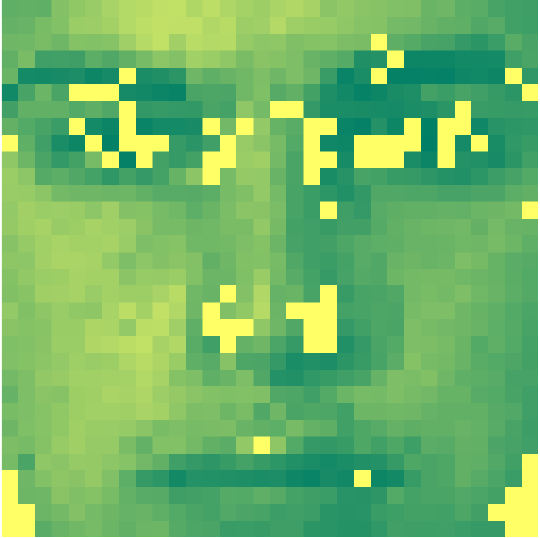} & \includegraphics[scale=0.1,valign=c]{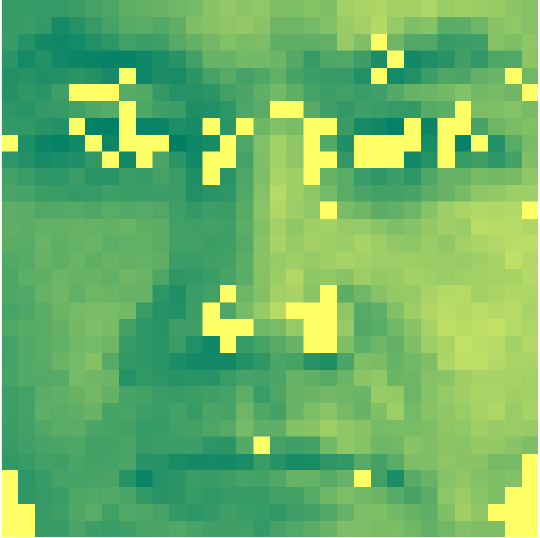} & 41.05 & 65.13 \\
\textbf{Case III} & \includegraphics[scale=0.1,valign=c]{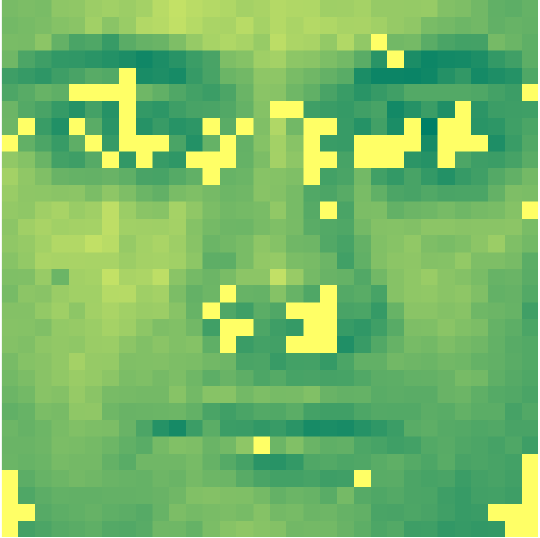} & \includegraphics[scale=0.1,valign=c]{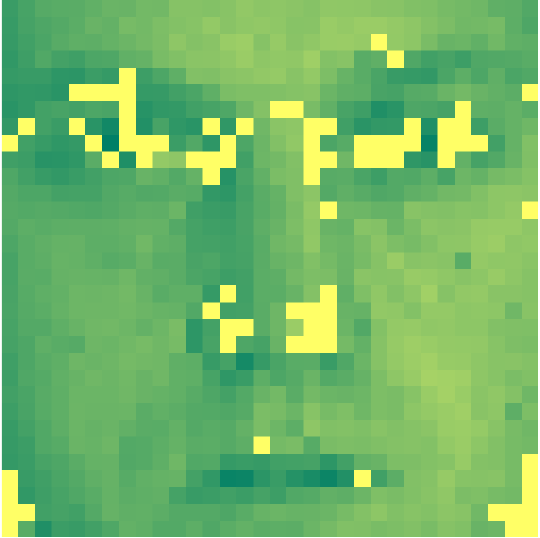} & \includegraphics[scale=0.1,valign=c]{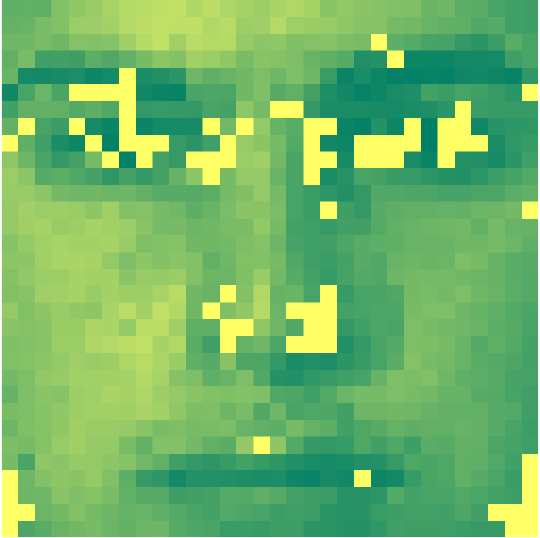} & \includegraphics[scale=0.1,valign=c]{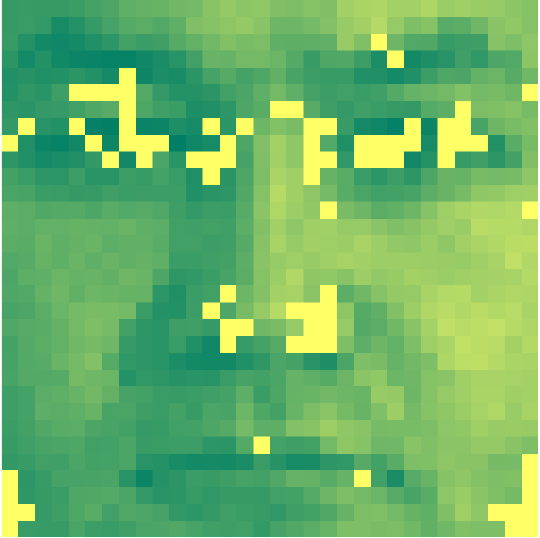} & \textcolor[rgb]{0.00,0.00,1.00}{\textbf{41.15}} & \textcolor[rgb]{0.00,0.00,1.00}{\textbf{65.23}}\\
\textbf{Case IV} & \includegraphics[scale=0.1,valign=c]{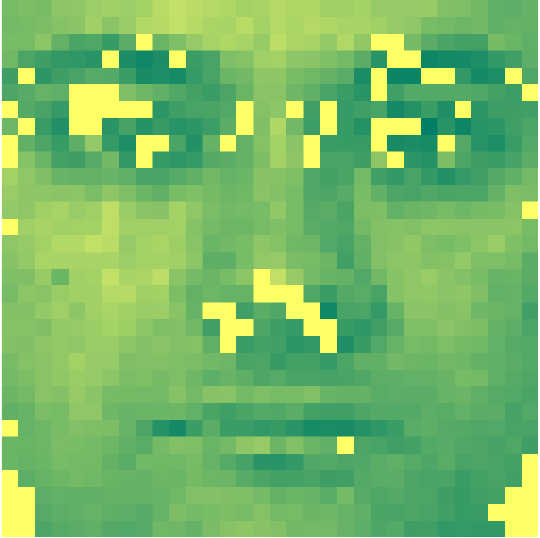} & \includegraphics[scale=0.1,valign=c]{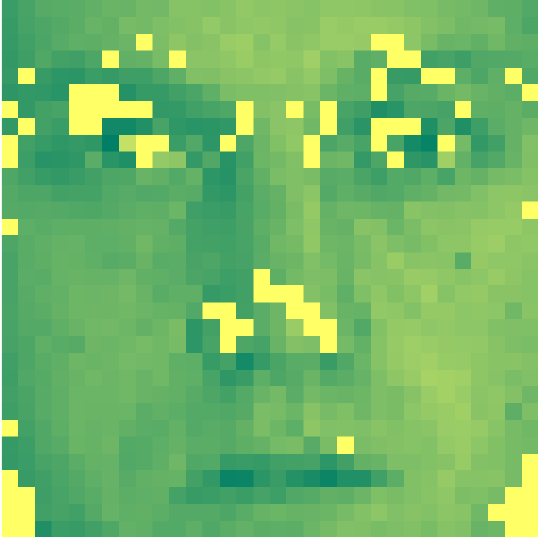} & \includegraphics[scale=0.1,valign=c]{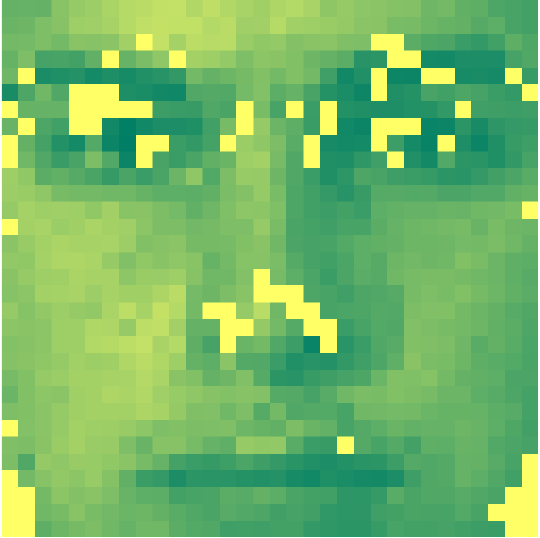} & \includegraphics[scale=0.1,valign=c]{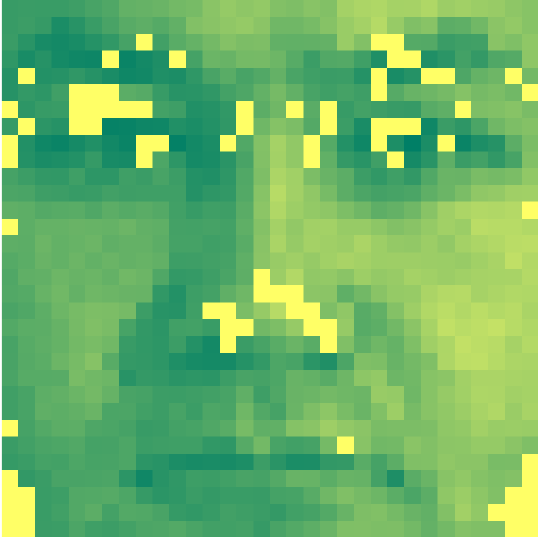} & \textcolor[rgb]{1.00,0.00,0.00}{\textbf{42.45}} & \textcolor[rgb]{1.00,0.00,0.00}{\textbf{66.66}} \\
\end{tblr}
\end{table}

\subsection{Real-world Results}\label{experiments-real}

This section presents numerical results on eight real-world datasets. Here is another compared method called ALLfea, which stands for all features used for clustering and is the gold standard for comparison. 
Fig. \ref{plot-acc} and Fig. \ref{plot-nmi} show the average values of ACC and NMI for 50 repetitions  under different selected feature numbers. 
Table \ref{acc} and Table \ref{nmi} summarize the detailed results, with the top two values marked in \textcolor[rgb]{1.00,0.00,0.00}{\textbf{red}} and \textcolor[rgb]{0.00,0.00,1.00}{\textbf{blue}} except for ALLfea.
To be specific, we set the number of selected features from 10 to 100 and report the best result with the number of features shown in brackets.

\textbf{For the ACC metric, our proposed BSUFS  obtains good results on most datasets, even performs better than the latest SPCA-PSD and FEN-PCAFS.} In Fig. \ref{plot-acc}, almost all BSUFS lines are  higher than other lines under different numbers of selected features, which implies that BSUFS performs better than others in terms of ACC. 
In Table \ref{acc}, BSUFS has the largest average ACC value on these real-world datasets, followed by FEN-PCAFS and SPCA-PSD.
In addition, compared with graph-based methods (such as RNE and UDFS),  PCA-based methods (such as SPCAFS and SPCA-PSD) show good performance. 
Due to the introduction of bi-sparse regularization, BSUFS achieves at least an average improvement of 7.95\% than SPCAFS. Of course, in some cases, BSUFS is slightly inferior to SPCA-PSD, we believe that low-rank priors can also improve the performance of UFS.
For the Isolet dataset, the ACC value is significantly improved. The reason may be that the fine-grained noise and speaker variability in this dataset can be effectively handled by the additional $\ell_q$-norm, thereby achieving more effective features.

\textbf{For the corresponding NMI metric, similar results can be obtained, that is, the proposed BSUFS generally outperforms other competitors.} Note that the values of NMI here adopt these parameters corresponding to the best ACC, which may cause NMI to be not such good in some cases. Although not all BSUFS lines in Fig. \ref{plot-nmi} are higher than others, they are still the most lines.  As shown in Table \ref{nmi}, on average, BSUFS improves the NMI value by at least 0.48\% compared with other methods. For the GLIOMA dataset with a smaller sample size and a larger number of features, it achieves a higher ACC value but a relatively lower NMI value. This may be because this dataset has only four categories, which easily leads to an imbalance in ACC and NMI.

Overall, our proposed BSUFS outperforms compared methods on many real-world datasets by achieving higher ACC and NMI values. On the other hand, achieving higher accuracy with a smaller number of selected features makes BSUFS a very practical method for real-world applications.

\begin{figure*}[t]
\centering
\hspace{-0.1cm}
\subfigcapskip=-1pt
\subfigure[USPS (Case I)]{
    \label{a}
    \centering
    \includegraphics[width=4cm]{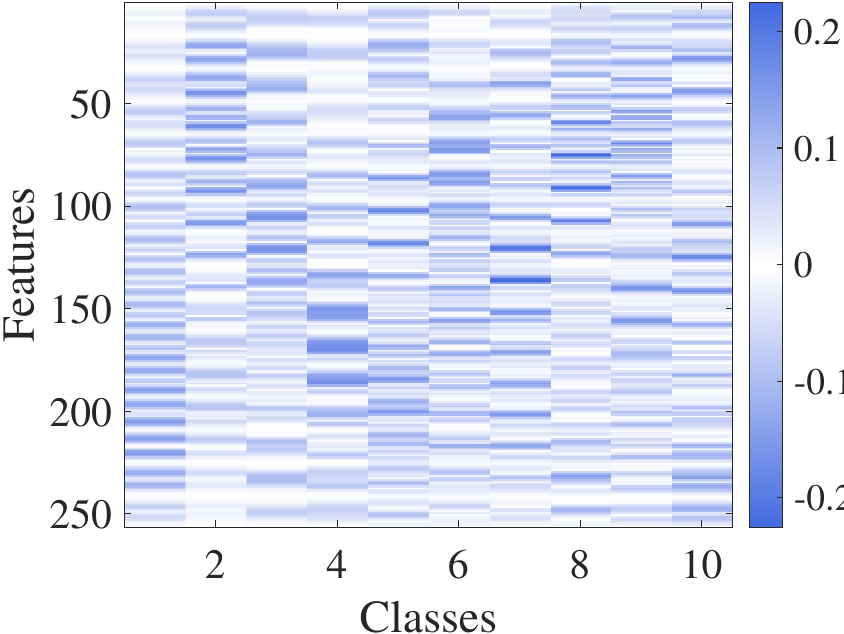}
}\hspace{-0mm}
\subfigcapskip=-1pt
\subfigure[USPS (Case II)]{
    \label{b}
    \centering
    \includegraphics[width=4cm]{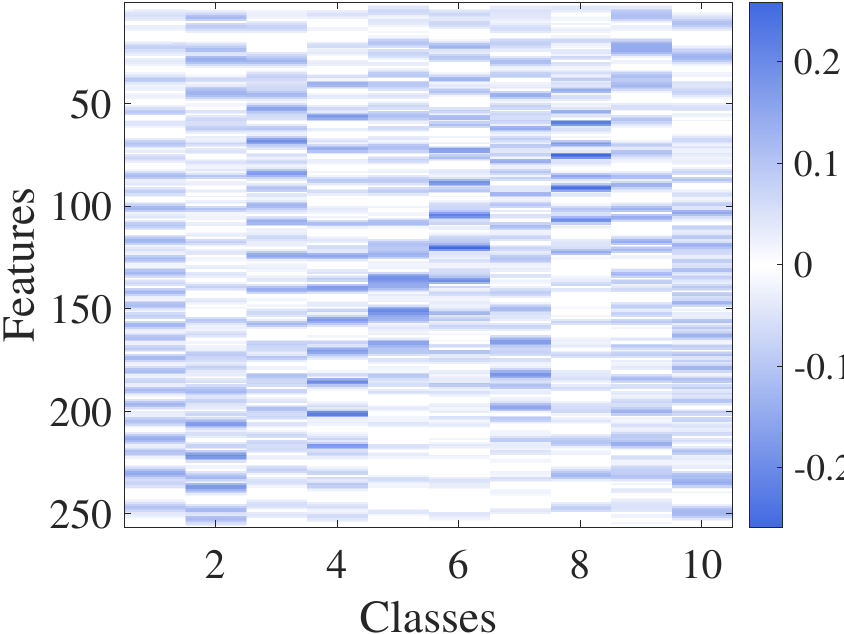}
}\hspace{-0mm}
\subfigcapskip=-1pt
\subfigure[USPS (Case III)]{
    \label{c}
    \centering
    \includegraphics[width=4cm]{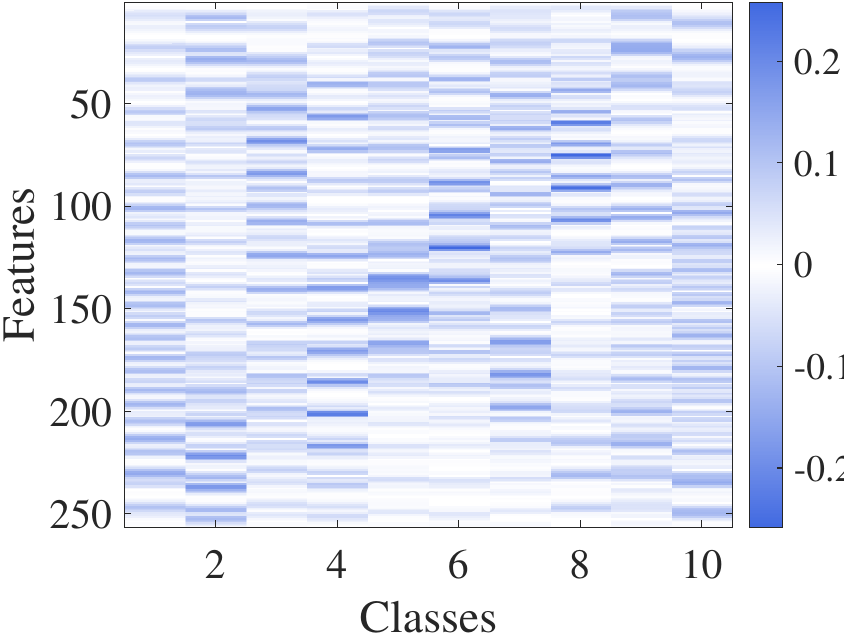}
}\hspace{-0mm}
\subfigcapskip=-1pt
\subfigure[USPS (Case IV)]{
    \label{d}
    \centering
    \includegraphics[width=4cm]{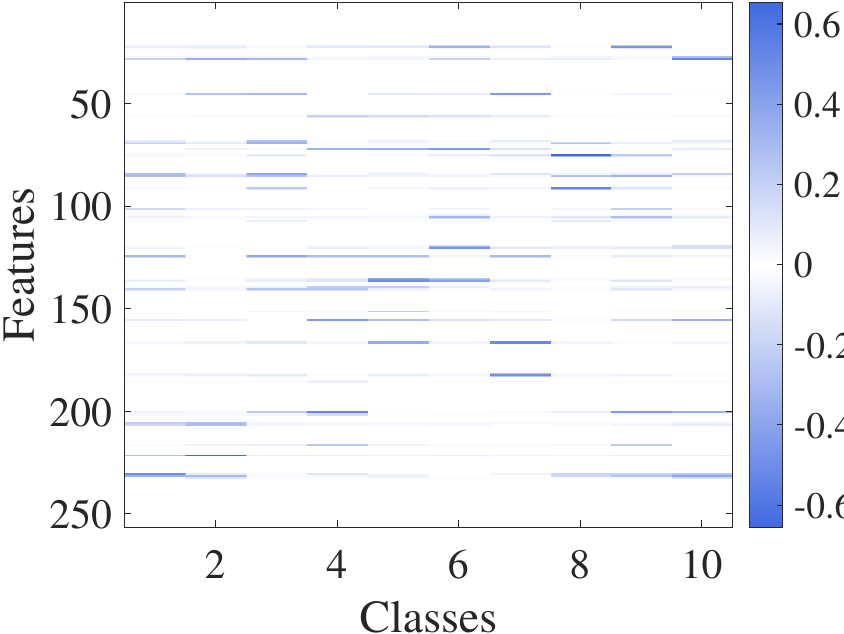}
}\hspace{-0mm}
\subfigcapskip=-1pt
\subfigure[umist (Case I)]{
    \label{a}
    \centering
    \includegraphics[width=4cm]{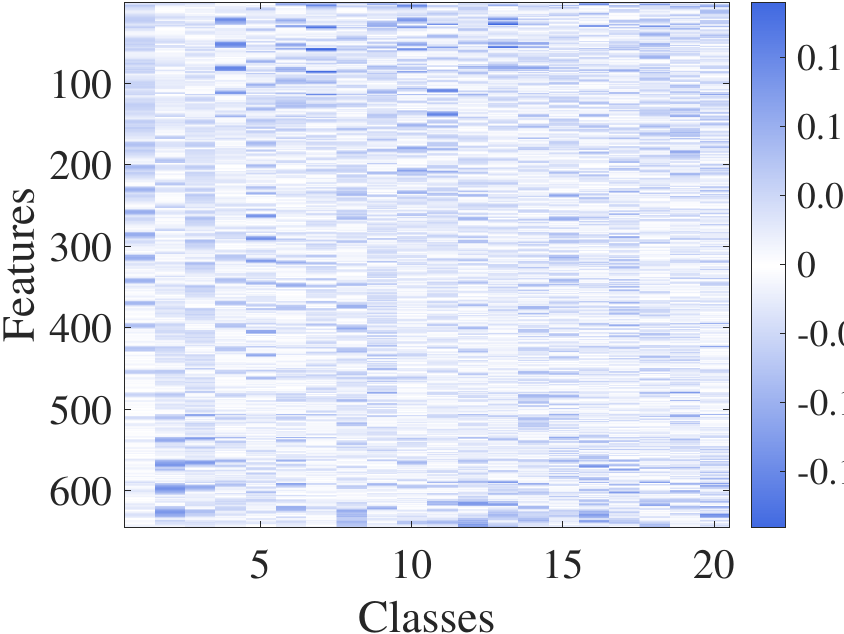}
}\hspace{-0mm}
\subfigcapskip=-1pt
\subfigure[umist (Case II)]{
    \label{b}
    \centering
    \includegraphics[width=4cm]{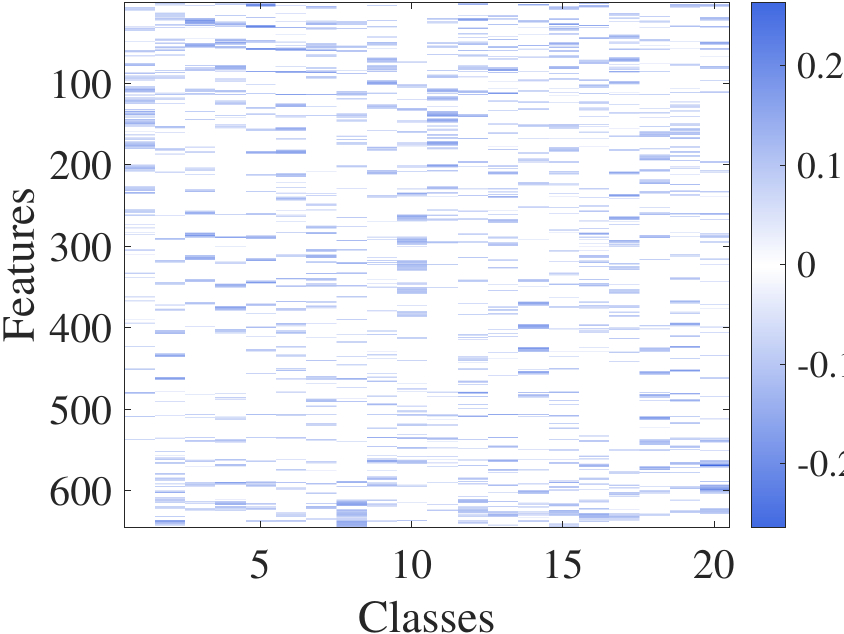}
}\hspace{-0mm}
\subfigcapskip=-1pt
\subfigure[umist (Case III)]{
    \label{c}
    \centering
    \includegraphics[width=4cm]{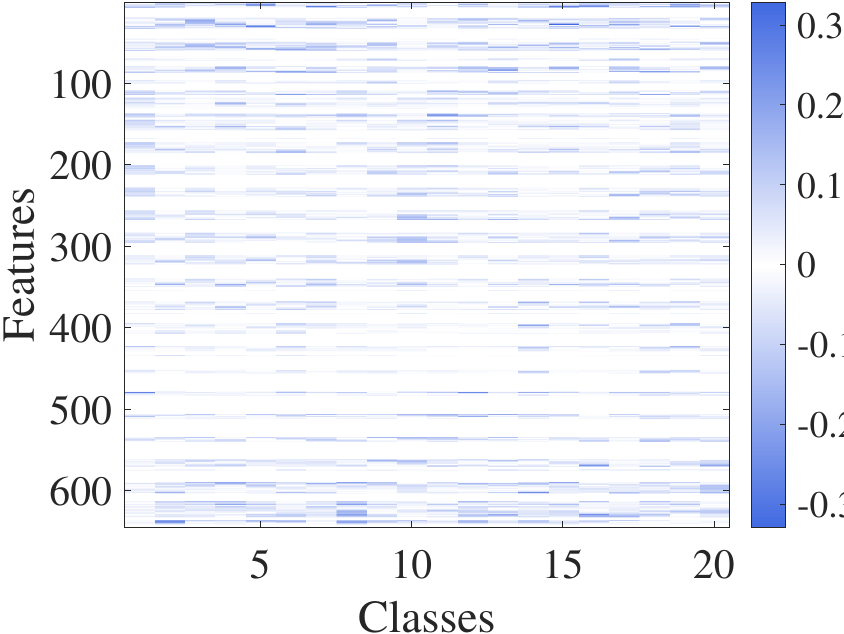}
}\hspace{-0mm}
\subfigcapskip=-1pt
\subfigure[umist (Case IV)]{
    \label{d}
    \centering
    \includegraphics[width=4cm]{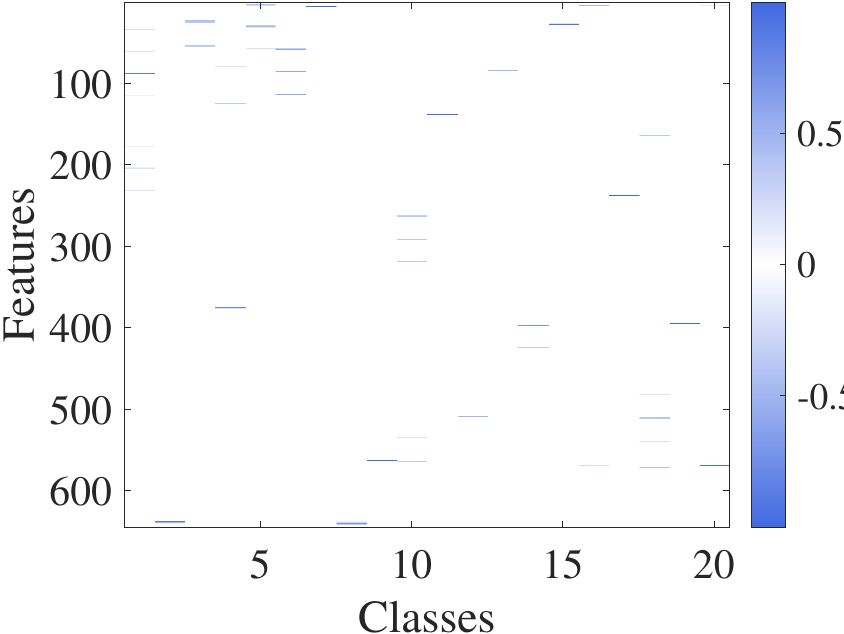}
}\hspace{-0mm}
\caption{Sparse visualization of the transformation matrix, where (a)-(d) are the results on the USPS dataset and (e)-(h) are the results on the umist dataset.}
\label{heat}
\end{figure*}

\subsection{Ablation Experiments} \label{ablation}

To investigate the effect of bi-sparse regularization terms in BSUFS, this section conducts ablation experiments on 
\begin{itemize}
  \item Case I: BSUFS without $\ell_{2,p}$-norm and $\ell_{q}$-norm, i.e.,
\begin{equation}
        \begin{aligned}
                \min_{W\in\mathbb{R}^{d\times m}}~ &-\textrm{Tr}(W^{\top} SW)\\
                \textrm{s.t.}~~~~&W^{\top} W=I_m.
        \end{aligned}
\end{equation}
  \item Case II: BSUFS without $\ell_{2,p}$-norm, i.e.,
  \begin{equation}
        \begin{aligned}
                \min_{W\in\mathbb{R}^{d\times m}}~ &-\textrm{Tr}(W^{\top} SW)+\lambda_2 \|W\|_{q}^{q} \\
                \textrm{s.t.}~~~~&W^{\top} W=I_m.
        \end{aligned}
\end{equation}
  \item Case III: BSUFS without $\ell_{q}$-norm, i.e.,
  \begin{equation}
        \begin{aligned}
                \min_{W\in\mathbb{R}^{d\times m}}~ &-\textrm{Tr}(W^{\top} SW)+\lambda_1 \|W\|_{2,p}^{p}  \\
                \textrm{s.t.}~~~~&W^{\top} W=I_m.
        \end{aligned}
\end{equation}
  \item Case IV: BSUFS.
\end{itemize}
Note that $p,q\in [0,1)$ for all cases.

\subsubsection{Clustering Efficacy}

Table \ref{a-acc} and Table \ref{ff-nmi} list  the clustering performance in terms of ACC and NMI, respectively. These results indicate that Case IV consistently achieves the top performance. In particular, on the USPS and GLIOMA datasets, the ACC values of Case IV increase from 67.06\% to 70.77\% and from 49.76\% to 61.28\% compared to Case I, respectively. In addition, compared with Case III, almost all ACC and NMI results of Case IV are improved, which shows that the introduction of $\ell_q$-norm to Case III is meaningful for feature selection.

\subsubsection{Feature Visualization}

Table \ref{graph} presents visual comparisons of feature selection results of Cases I to IV on the pie dataset. In this study, we set the number of selected features to 80 and randomly select 4 image samples to highlight the effectiveness of feature selection. Although all methods can select the basic facial features (eyes, mouth, nose, lips), Case IV selects more additional features, such as eyebrows. This diverse selection helps to maintain a more complete geometric structure of the face, potentially making better use of smaller regions of an image and reducing redundant features. Consequently, Case IV achieves higher ACC and NMI values, which reflects its superior effectiveness.

\begin{figure}[t]
    \centering
    \includegraphics[width=0.50\textwidth]{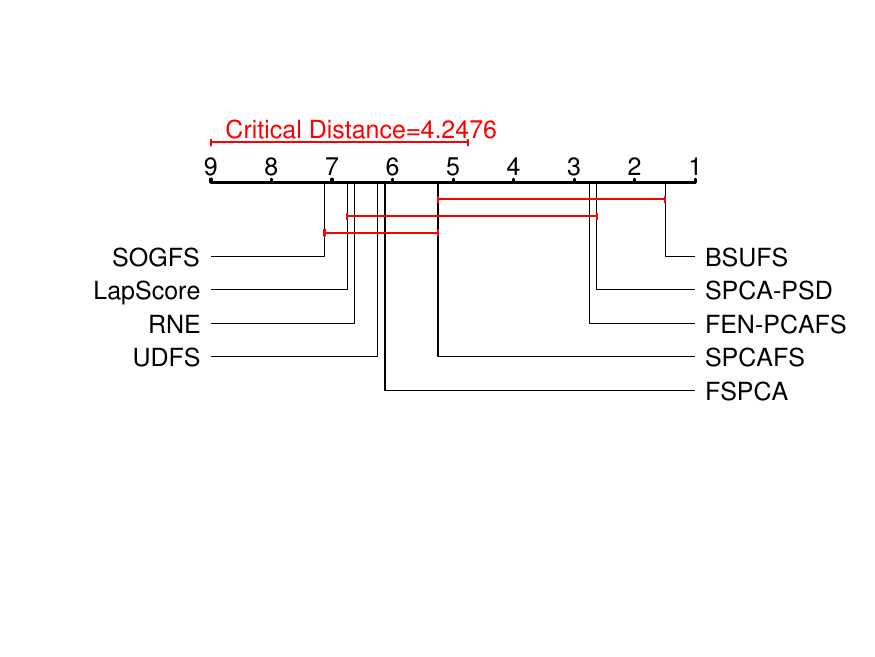}
    \vspace{-32mm}
    \caption{Post-hoc Nemenyi test in terms of ACC.}\label{hacc}
    \label{ft}
\end{figure}

\begin{figure}[t]
    \centering
    \includegraphics[width=0.50\textwidth]{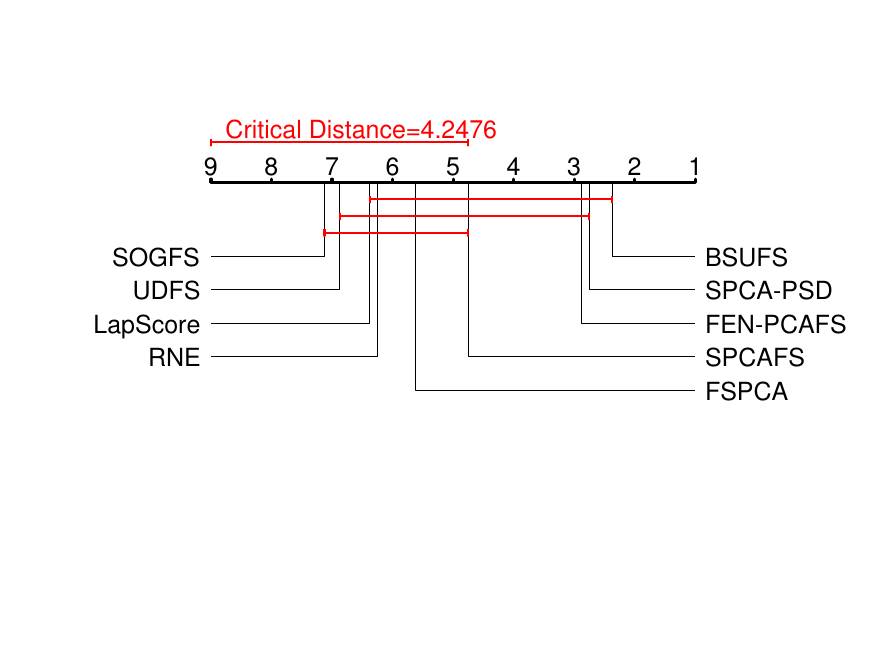}
    \vspace{-32mm}
    \caption{Post-hoc Nemenyi test in terms of NMI.}\label{hnmi}
    \label{ft}
\end{figure}

\subsubsection{Sparse Analysis}
Fig. \ref{heat} shows the sparse visualization of the transformation matrix $W$ on the USPS and umist datasets. Since both USPS and umist are image datasets, they usually contain a lot of noise. Obviously, Case IV combines the basic sparsity of Case II and the row sparsity of Case III to achieve a more sparse $W$, that is, more focused on effective features. This is because the introduction of $\ell_q$-norm regularization can eliminate noise in the datasets, thereby affecting feature selection and clustering.

\textbf{In summary, by comparing Cases I, II, III with Case IV in different measurements, the introduced bi-sparse term improves the performance of PCA in feature selection, i.e., our proposed BSUFS is promising.}

\begin{figure*}[t]
\centering
\hspace{-0.5cm}
\subfigcapskip=-3pt
\subfigure[COIL20]{
    \label{a}
    \centering
    \includegraphics[width=4cm]{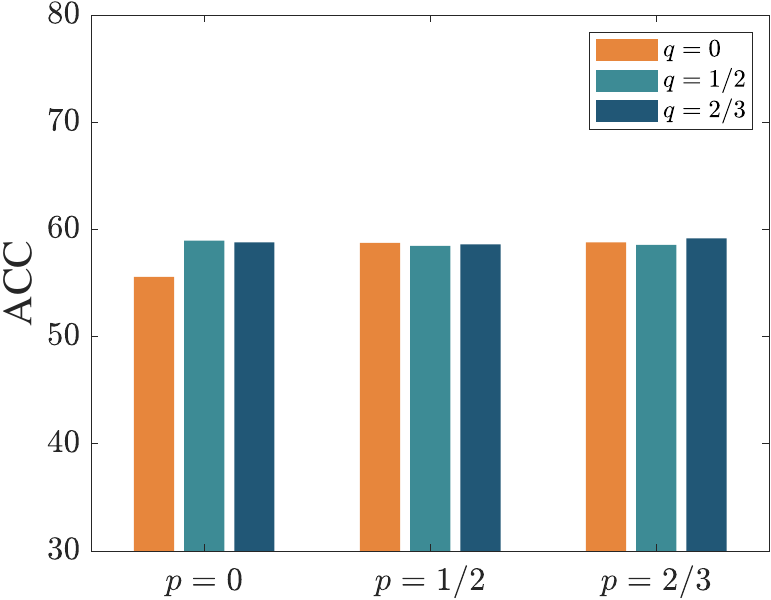}
}\hspace{-0mm}
\subfigcapskip=-1pt
\subfigure[Isolet]{
    \label{b}
    \centering
    \includegraphics[width=4cm]{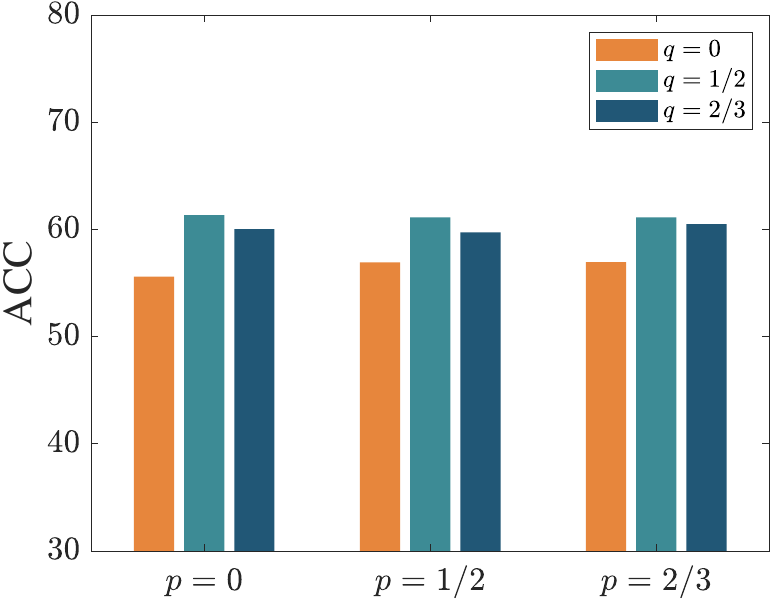}
}\hspace{-0mm}
\subfigcapskip=-1pt
\subfigure[USPS]{
    \label{c}
    \centering
    \includegraphics[width=4cm]{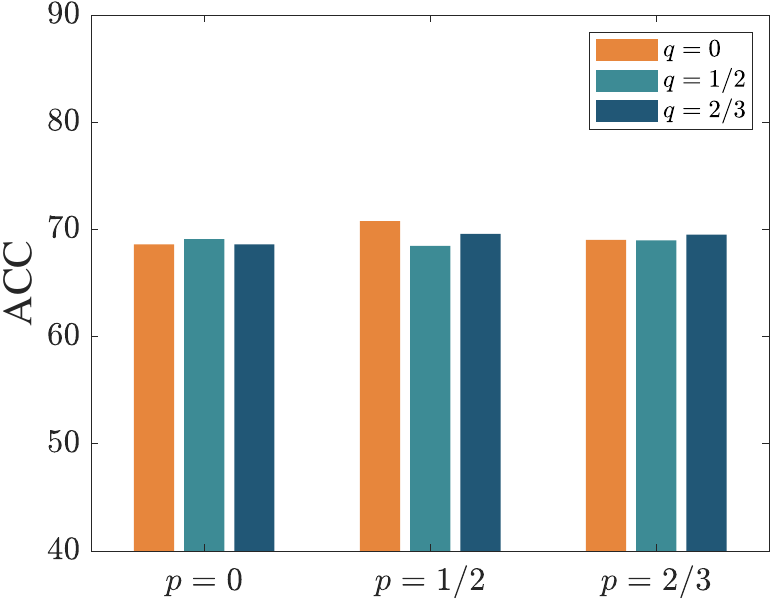}
}\hspace{-0mm}
\subfigcapskip=-1pt
\subfigure[umist]{
    \label{d}
    \centering
    \includegraphics[width=4cm]{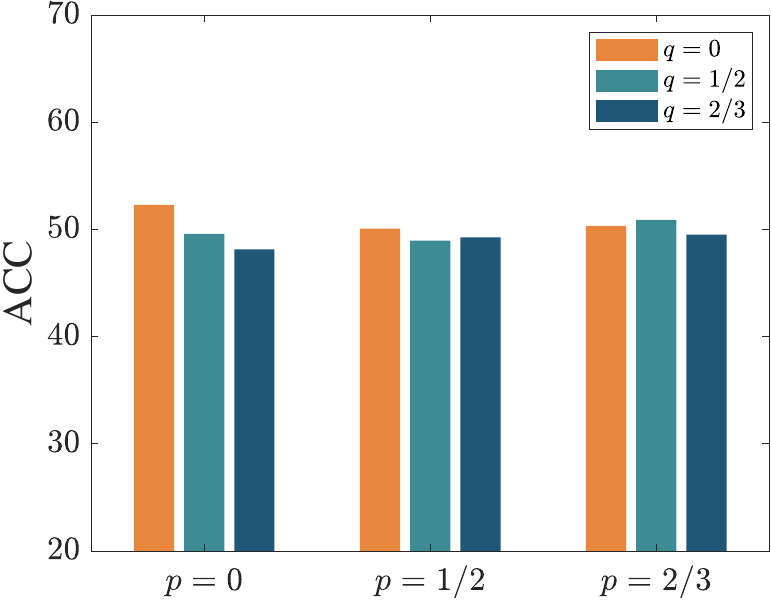}
}\hspace{-0mm}

\hspace{-0.5cm}
\subfigure[GLIOMA]{
    \label{e}
    \centering
    \includegraphics[width=4cm]{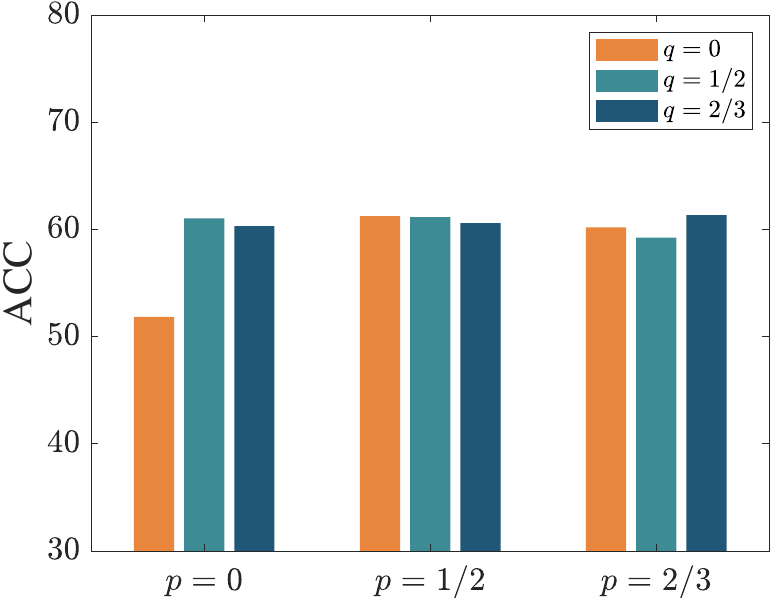}
}\hspace{-0mm}
\subfigcapskip=-1pt
\subfigure[pie]{
    \label{f}
    \centering
    \includegraphics[width=4cm]{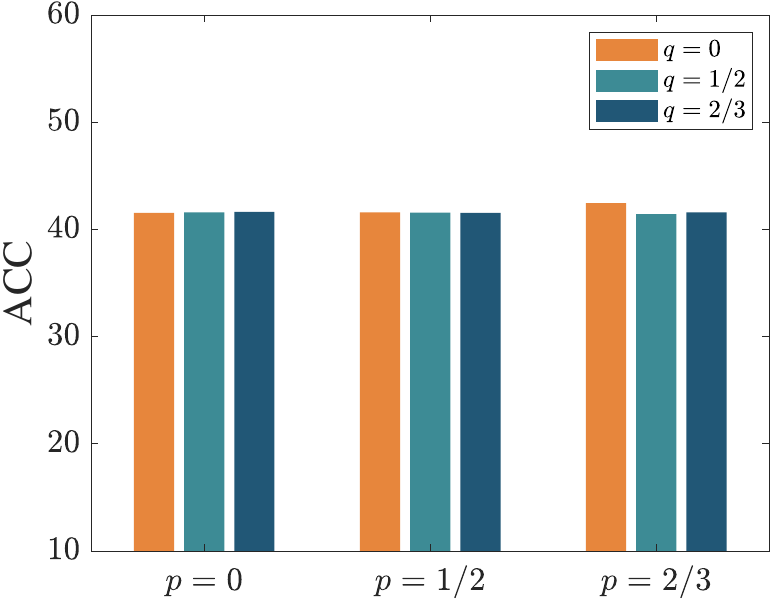}
}\hspace{-0mm}
\subfigcapskip=-1pt
\subfigure[LUNG]{
    \label{g}
    \centering
    \includegraphics[width=4cm]{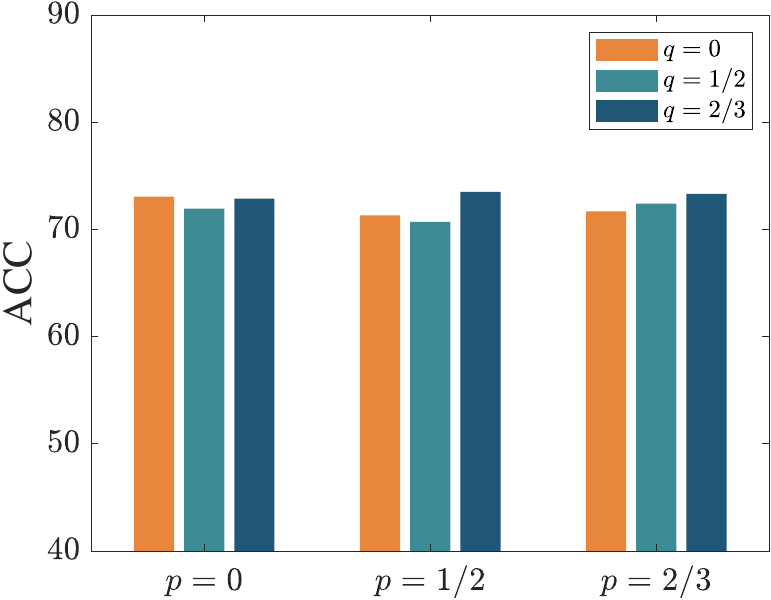}
}\hspace{-0mm}
\subfigcapskip=-1pt
\subfigure[MSTAR]{
    \label{h}
    \centering
    \includegraphics[width=4cm]{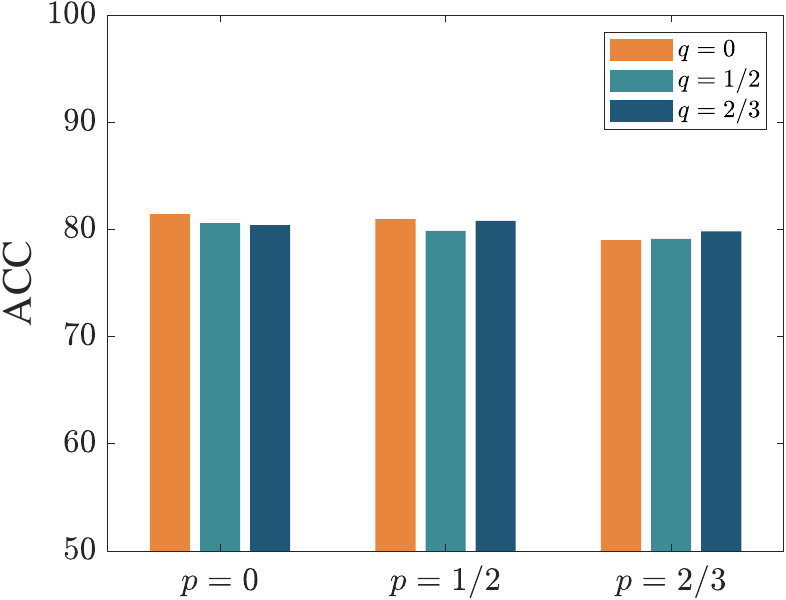}
}
\vspace{-0.1cm}
\caption{Effects of $p$ and $q$ on eight real-world datasets in terms of ACC  (\%).}
\label{pq-acc}
\end{figure*}

\begin{figure*}[t]
\centering
\hspace{-0.5cm}
\subfigcapskip=-1pt
\subfigure[COIL20]{
    \label{a}
    \centering
    \includegraphics[width=4cm]{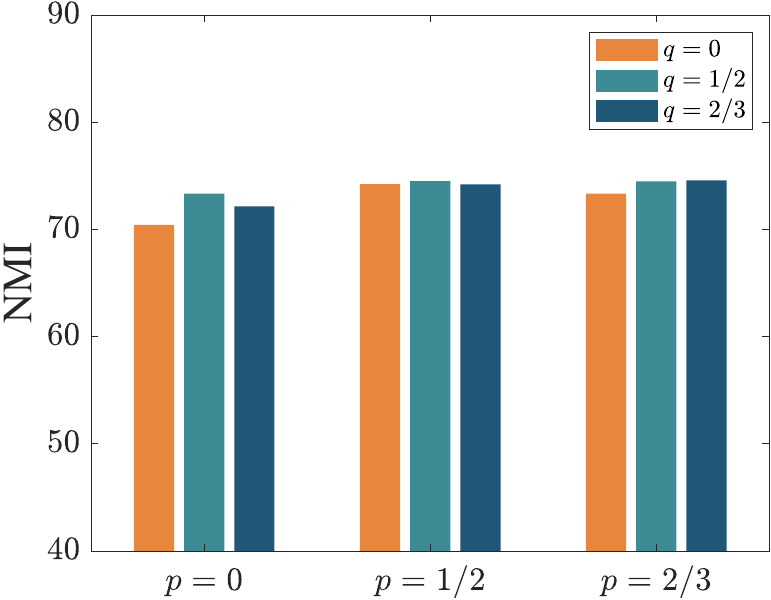}
}\hspace{-0mm}
\subfigcapskip=-1pt
\subfigure[Isolet]{
    \label{b}
    \centering
    \includegraphics[width=4cm]{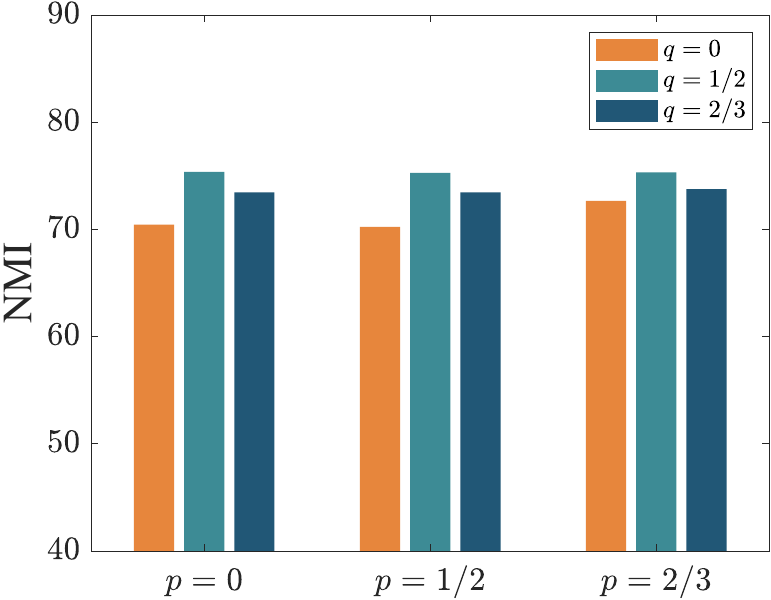}
}\hspace{-0mm}
\subfigcapskip=-1pt
\subfigure[USPS]{
    \label{c}
    \centering
    \includegraphics[width=4cm]{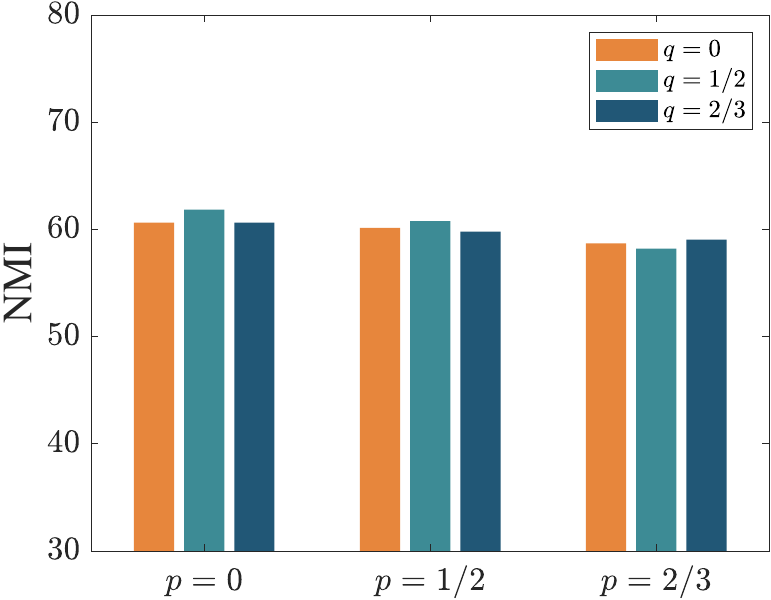}
}\hspace{-0mm}
\subfigcapskip=-1pt
\subfigure[umist]{
    \label{d}
    \centering
    \includegraphics[width=4cm]{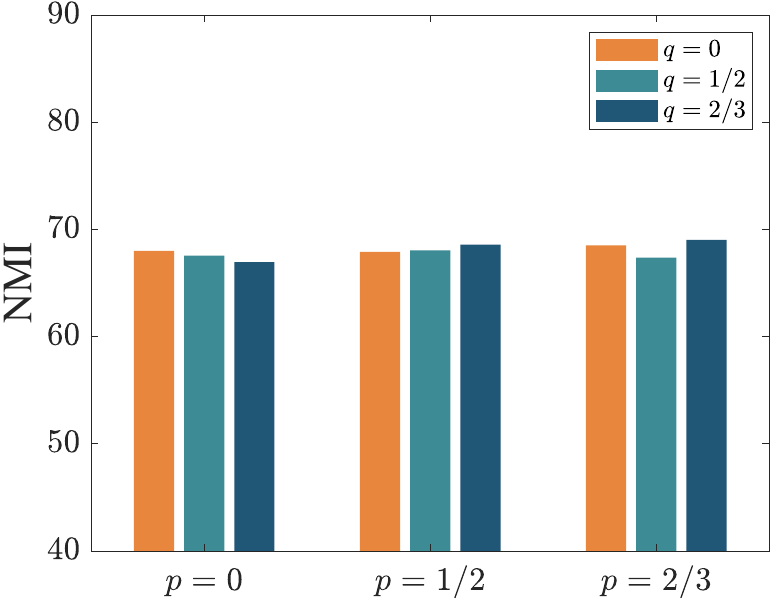}
}\hspace{-0mm}

\hspace{-0.5cm}
\subfigure[GLIOMA]{
    \label{e}
    \centering
    \includegraphics[width=4cm]{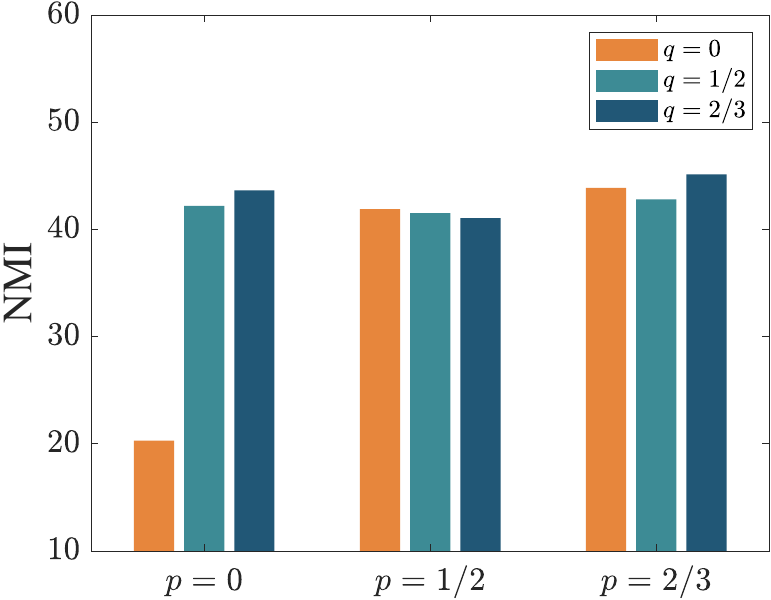}
}\hspace{-0mm}
\subfigcapskip=-1pt
\subfigure[pie]{
    \label{f}
    \centering
    \includegraphics[width=4cm]{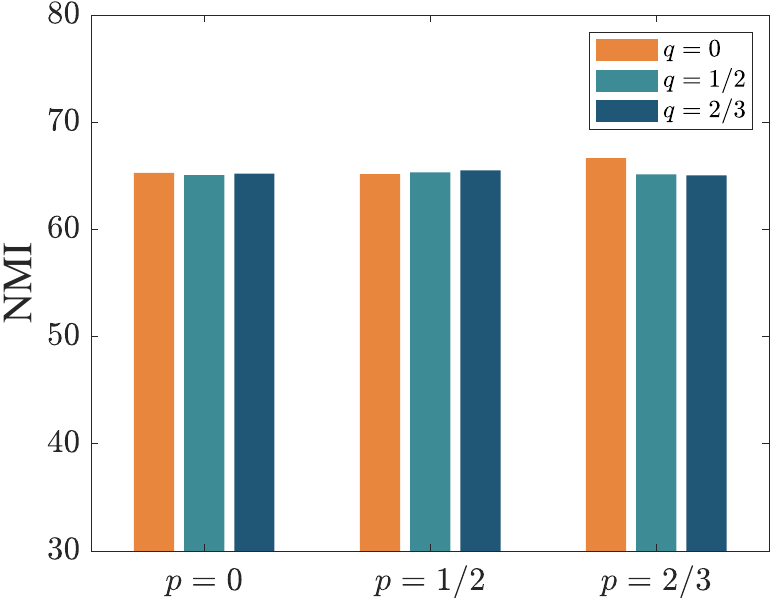}
}\hspace{-0mm}
\subfigcapskip=-1pt
\subfigure[LUNG]{
    \label{g}
    \centering
    \includegraphics[width=4cm]{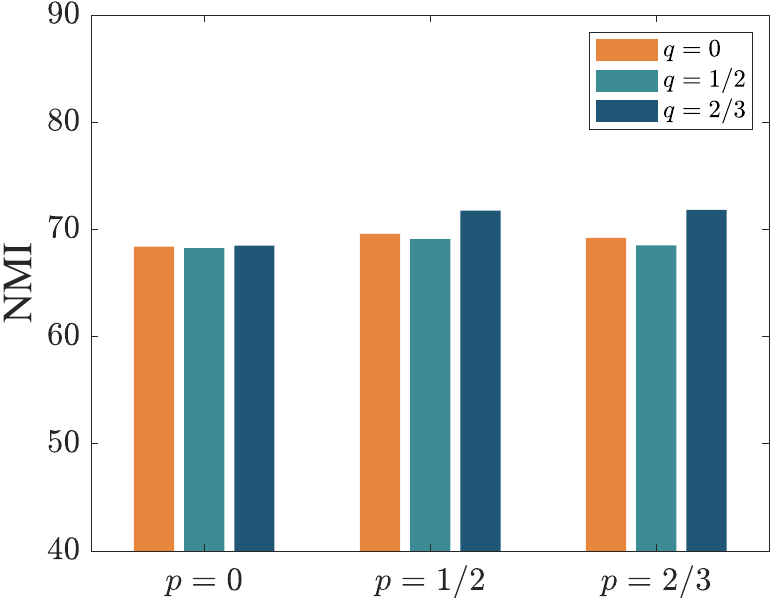}
}\hspace{-0mm}
\subfigcapskip=-1pt
\subfigure[MSTAR]{
    \label{h}
    \centering
    \includegraphics[width=4cm]{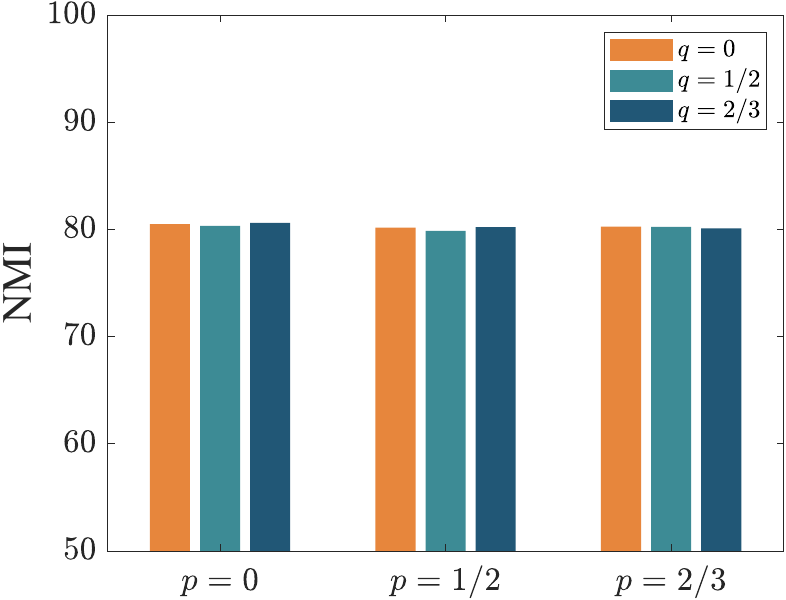}
}
\vspace{-0.1cm}
\caption{Effects of $p$ and $q$ on eight real-world datasets in terms of NMI  (\%).}
\label{pq-nmi}
\end{figure*}

\begin{figure*}[t]
\centering
\hspace{2mm}
\subfigcapskip=-1pt
\subfigure[COIL20 (SPCAFS)]{
    \label{a}
    \centering
    \includegraphics[width=4.1cm]{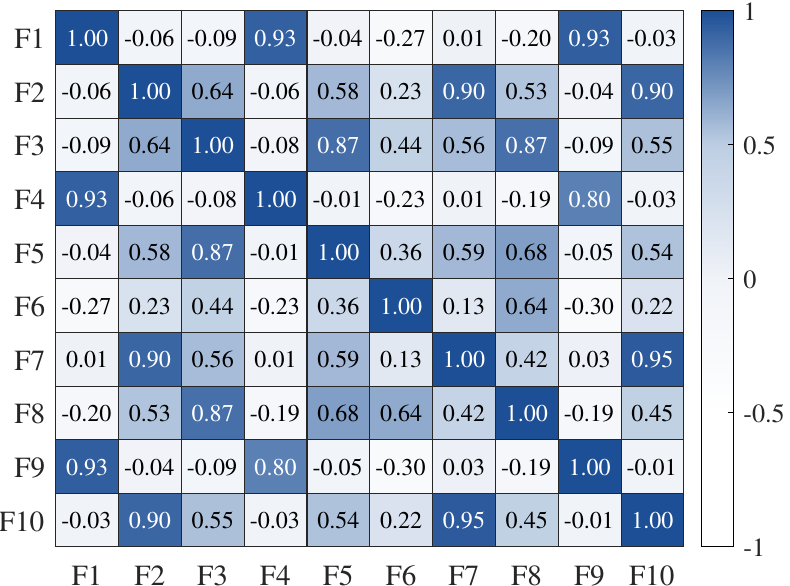}
}\hspace{-0mm}
\subfigcapskip=-1pt
\subfigure[Isolet (SPCAFS)]{
    \label{b}
    \centering
    \includegraphics[width=4.1cm]{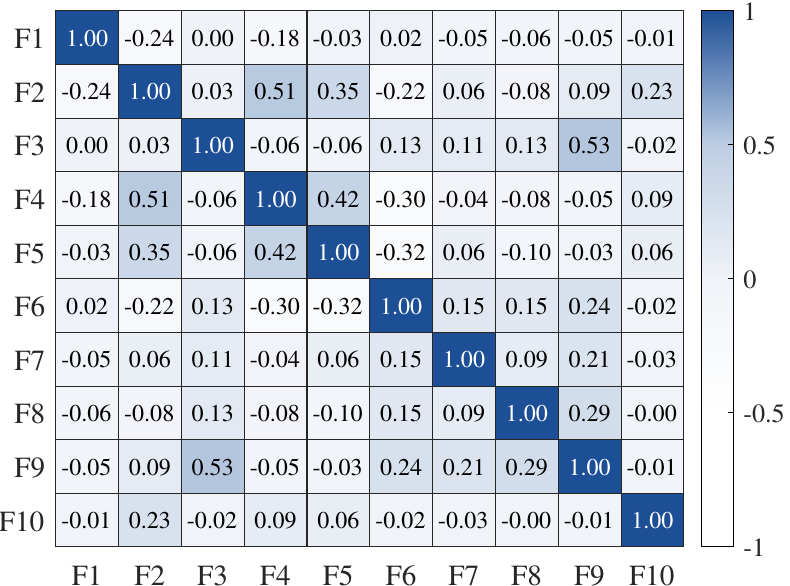}
}\hspace{-0mm}
\subfigcapskip=-1pt
\subfigure[USPS (SPCAFS)]{
    \label{c}
    \centering
    \includegraphics[width=4.1cm]{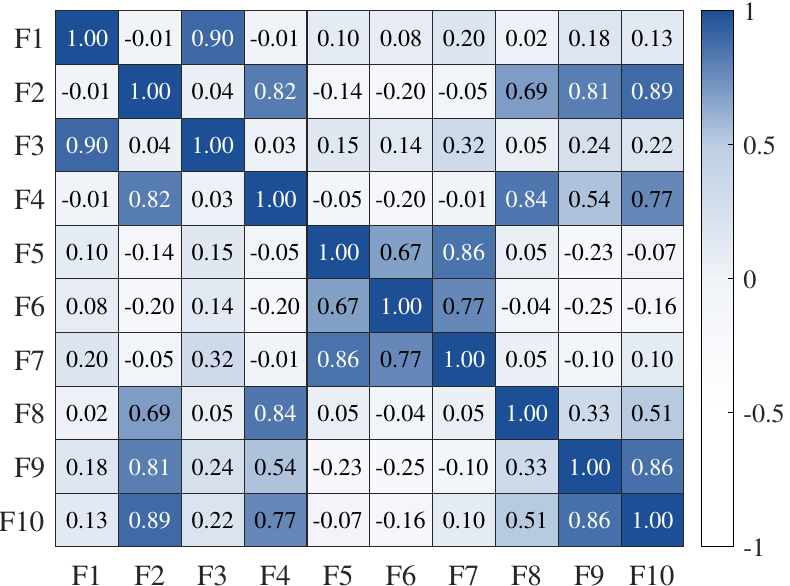}
}\hspace{-0mm}
\subfigcapskip=-1pt
\subfigure[LUNG (SPCAFS)]{
    \label{d}
    \centering
    \includegraphics[width=4.1cm]{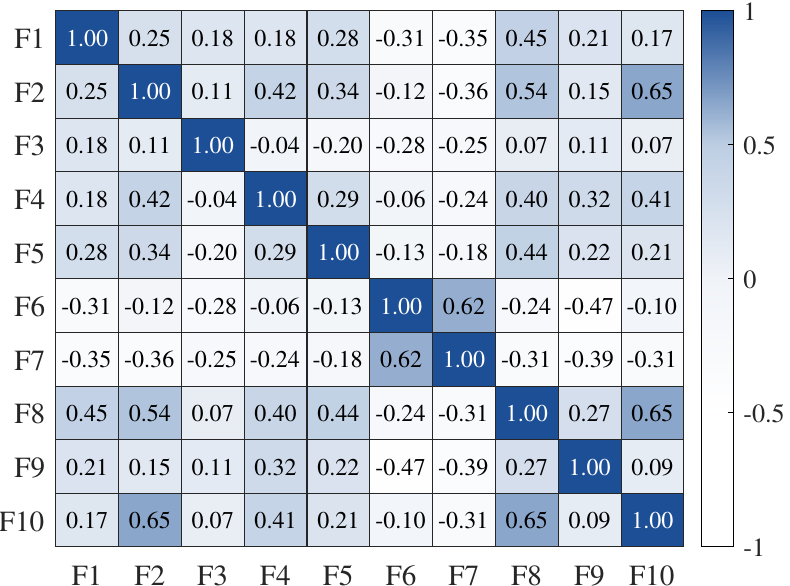}
}\hspace{-0mm}

\hspace{2mm}
\subfigure[COIL20 (BSUFS)]{
    \label{e}
    \centering
    \includegraphics[width=4.1cm]{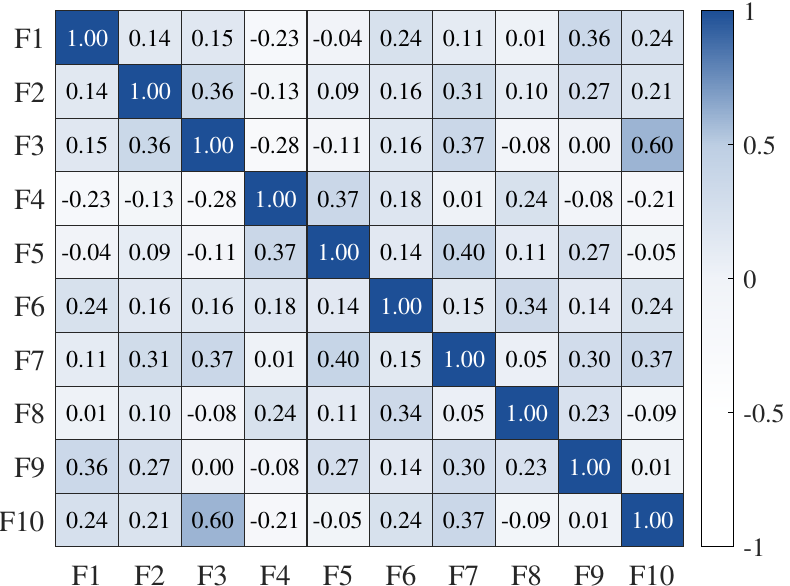}
}\hspace{-0mm}
\subfigcapskip=-1pt
\subfigure[Isolet (BSUFS)]{
    \label{f}
    \centering
    \includegraphics[width=4.1cm]{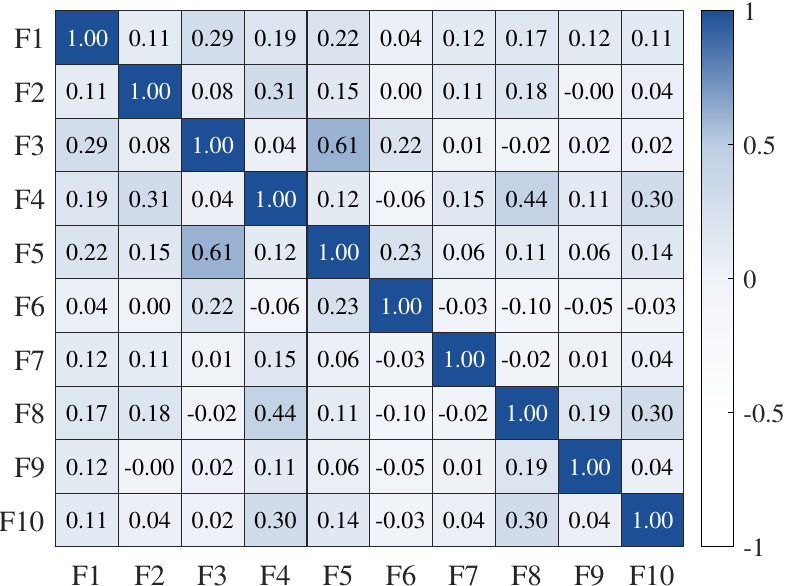}
}\hspace{-0mm}
\subfigcapskip=-1pt
\subfigure[USPS (BSUFS)]{
    \label{g}
    \centering
    \includegraphics[width=4.1cm]{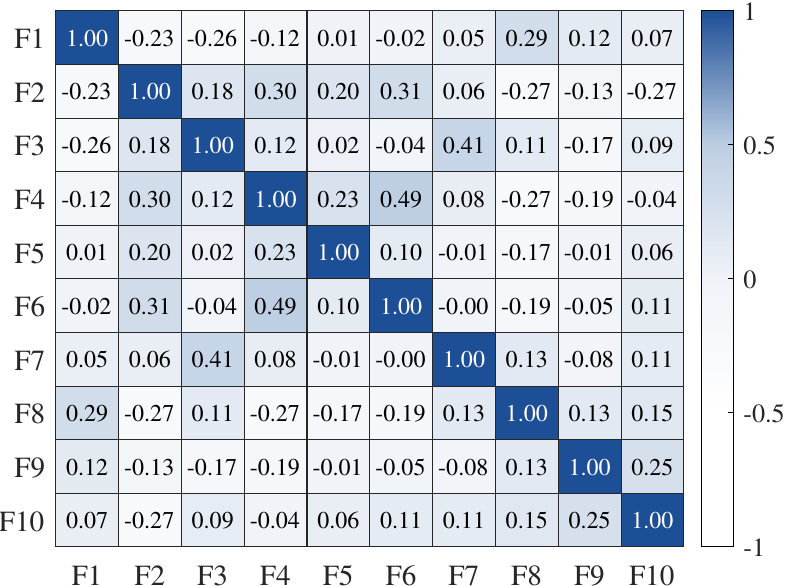}
}\hspace{-0mm}
\subfigcapskip=-1pt
\subfigure[LUNG (BSUFS)]{
    \label{h}
    \centering
    \includegraphics[width=4.1cm]{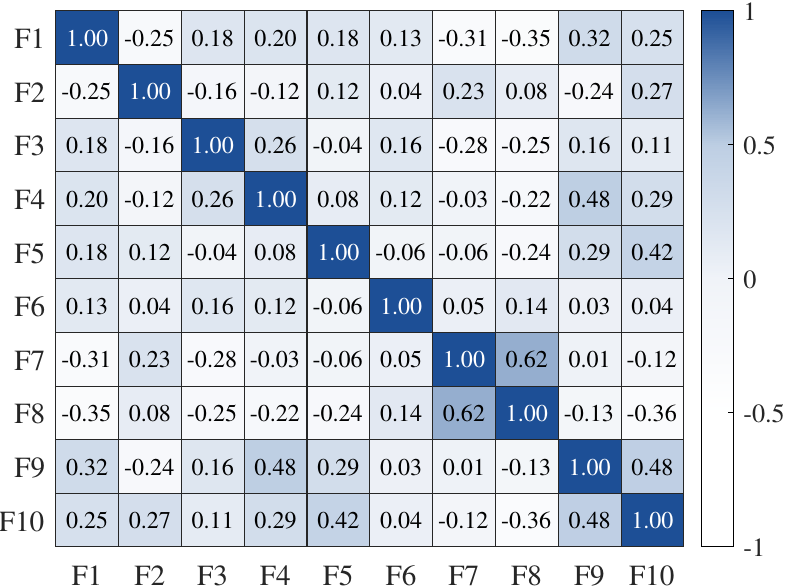}
}\hspace{-0mm}
\caption{Heatmap visualizations of correlations for 10 selected features, where (a)-(d) are the results of SPCAFS and (e)-(h) are the results of BSUFS.}
\centering
\label{correlation}
\end{figure*}

\begin{figure*}[t]
\hspace{-0.5cm}
\centering
\subfigcapskip=-2pt
\subfigure[COIL20 (ACC)]{
    \label{a}
    \centering
    \includegraphics[width=4.1cm]{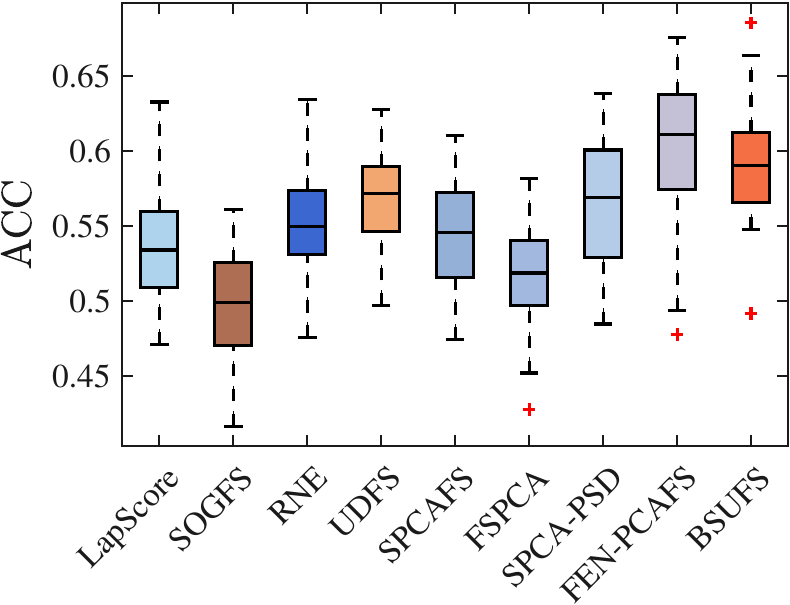}
}\hspace{0mm}
\subfigcapskip=-2pt
\subfigure[Isolet (ACC)]{
    \label{b}
    \centering
    \includegraphics[width=4.1cm]{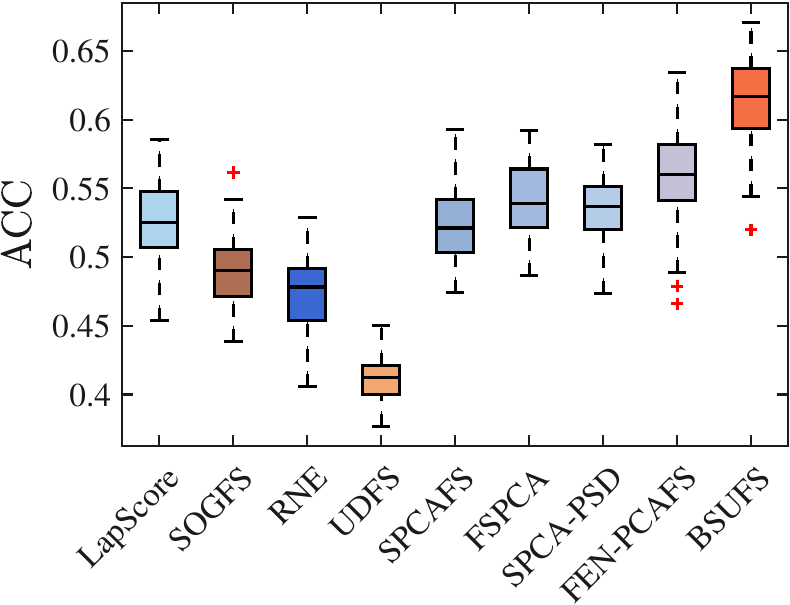}
}\hspace{0mm}
\subfigcapskip=-2pt
\subfigure[USPS (ACC)]{
    \label{c}
    \centering
    \includegraphics[width=4.1cm]{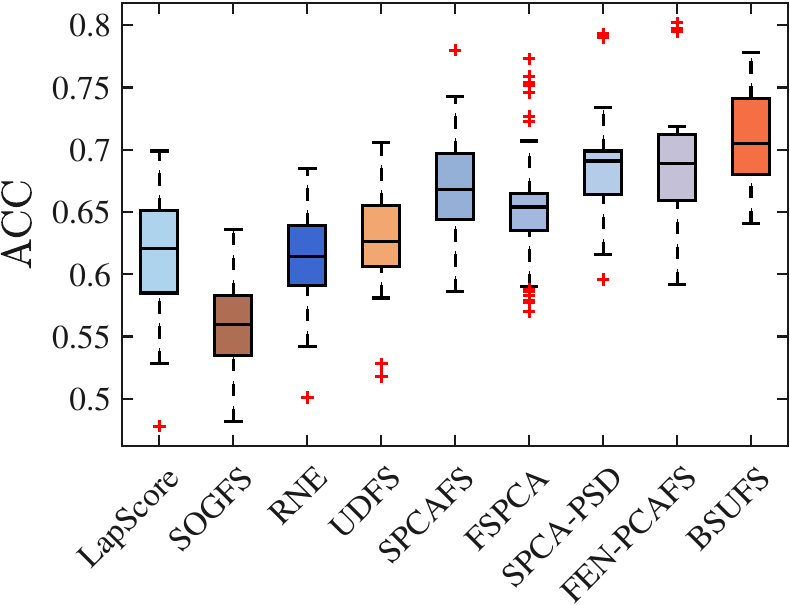}
}\hspace{0mm}
\subfigcapskip=-2pt
\subfigure[LUNG (ACC)]{
    \label{d}
    \centering
    \includegraphics[width=4.1cm]{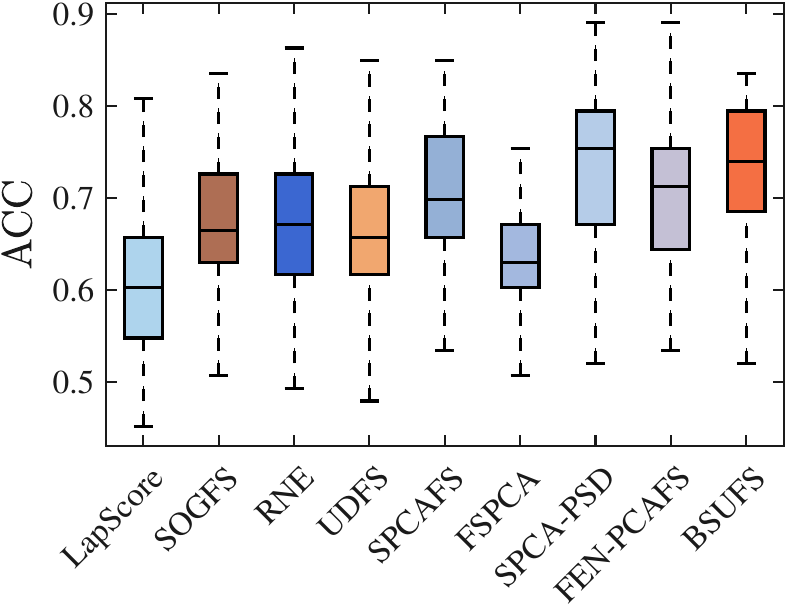}
}\hspace{0mm}

\hspace{-0.5cm}
\subfigure[COIL20 (NMI)]{
    \label{e}
    \centering
    \includegraphics[width=4.1cm]{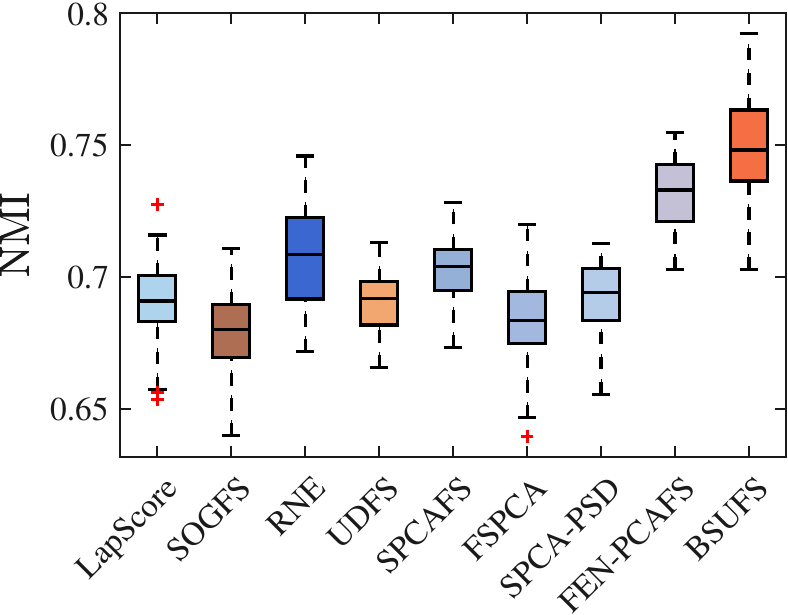}
}\hspace{0mm}
\subfigcapskip=-2pt
\subfigure[Isolet (NMI)]{
    \label{f}
    \centering
    \includegraphics[width=4.1cm]{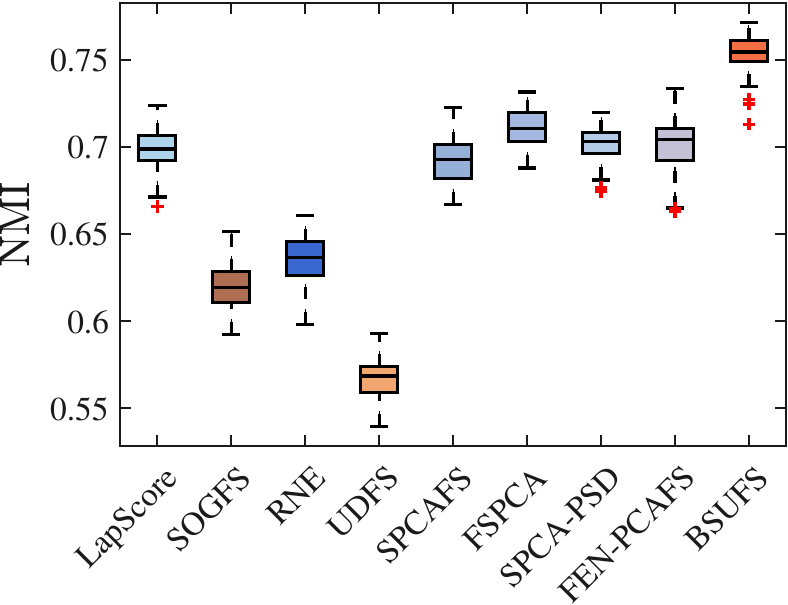}
}\hspace{0mm}
\subfigcapskip=-2pt
\subfigure[USPS (NMI)]{
    \label{g}
    \centering
    \includegraphics[width=4.1cm]{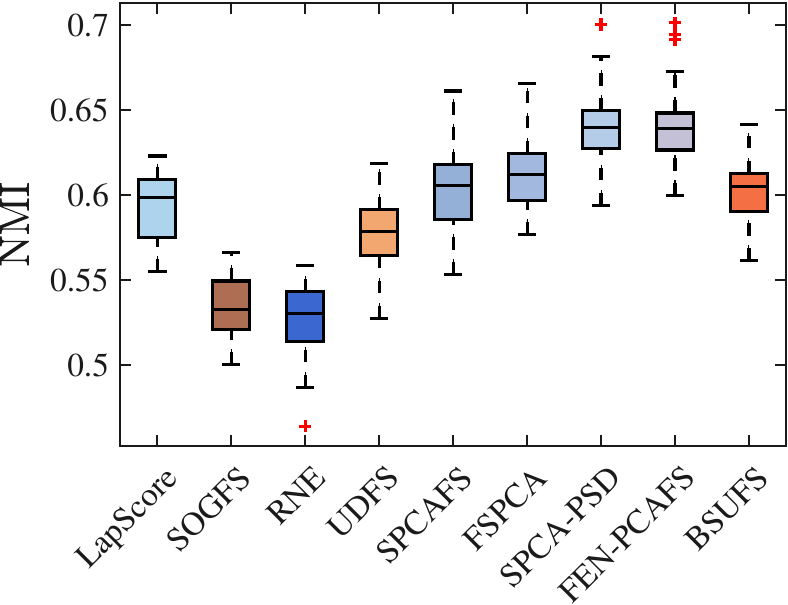}
}\hspace{0mm}
\subfigcapskip=-2pt
\subfigure[LUNG (NMI)]{
    \label{h}
    \centering
    \includegraphics[width=4.1cm]{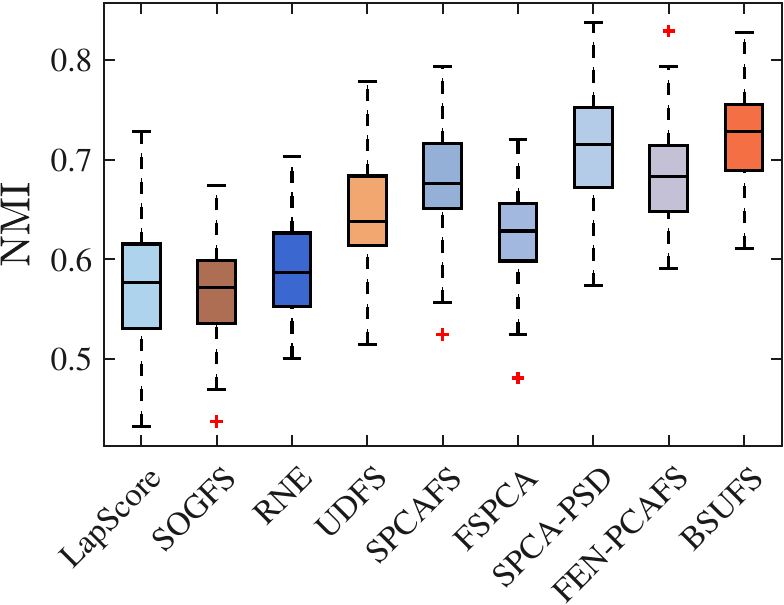}
}\hspace{0mm}
\caption{Model stability comparisons of all compared methods on four real-world datasets, where (a)-(d) are the ACC results and (e)-(h) are the NMI results.}
\label{box}
\end{figure*}

\begin{figure*}[t]
\hspace{-1cm}
\centering
\subfigcapskip=-3pt
\subfigure[COIL20 (ACC)]{
    \label{a}
    \centering
    \includegraphics[width=4.5cm]{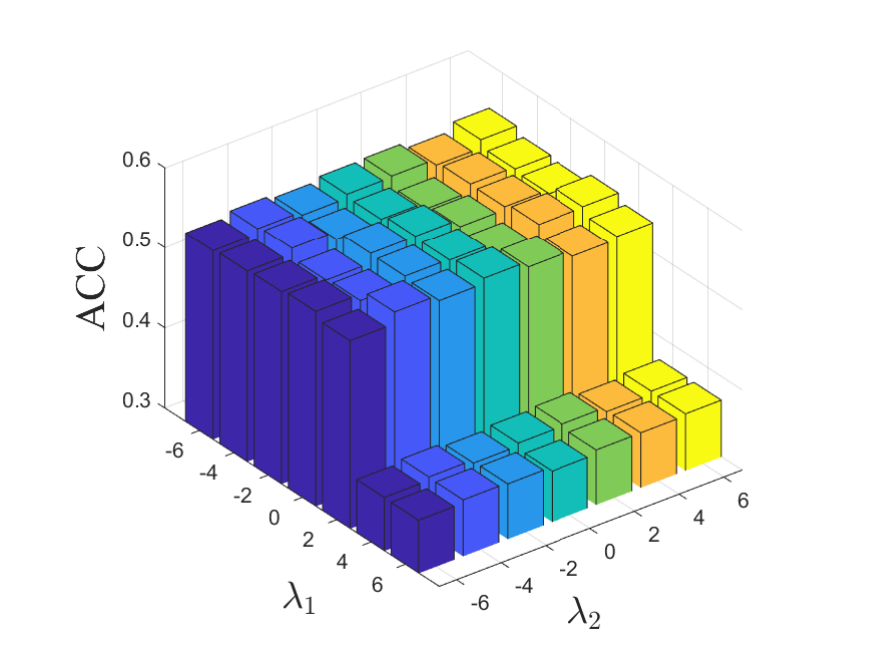}
}\hspace{-8mm}
\subfigcapskip=-3pt
\subfigure[Isolet (ACC)]{
    \label{b}
    \centering
    \includegraphics[width=4.5cm]{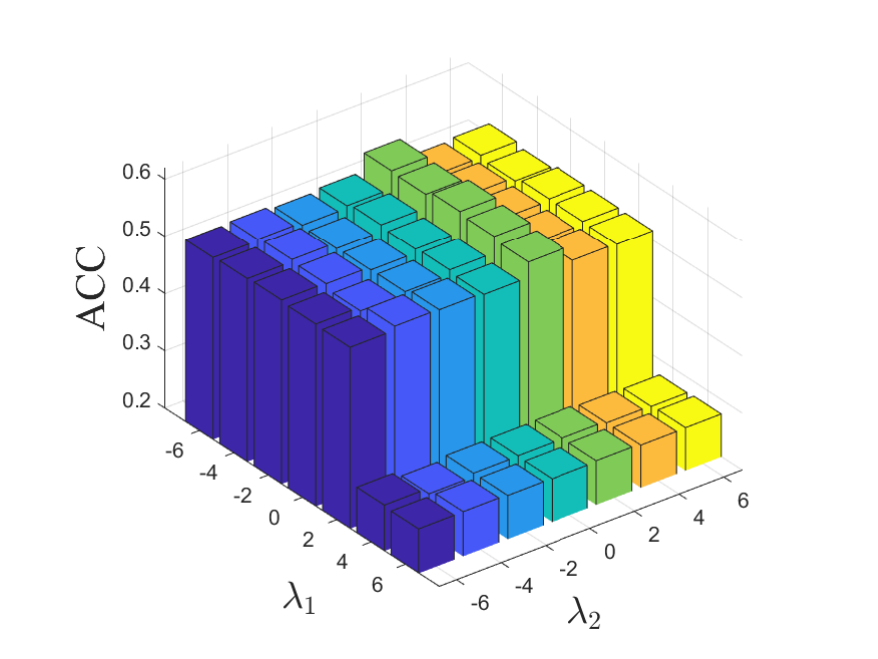}
}\hspace{-8mm}
\subfigcapskip=-3pt
\subfigure[USPS (ACC)]{
    \label{c}
    \centering
    \includegraphics[width=4.5cm]{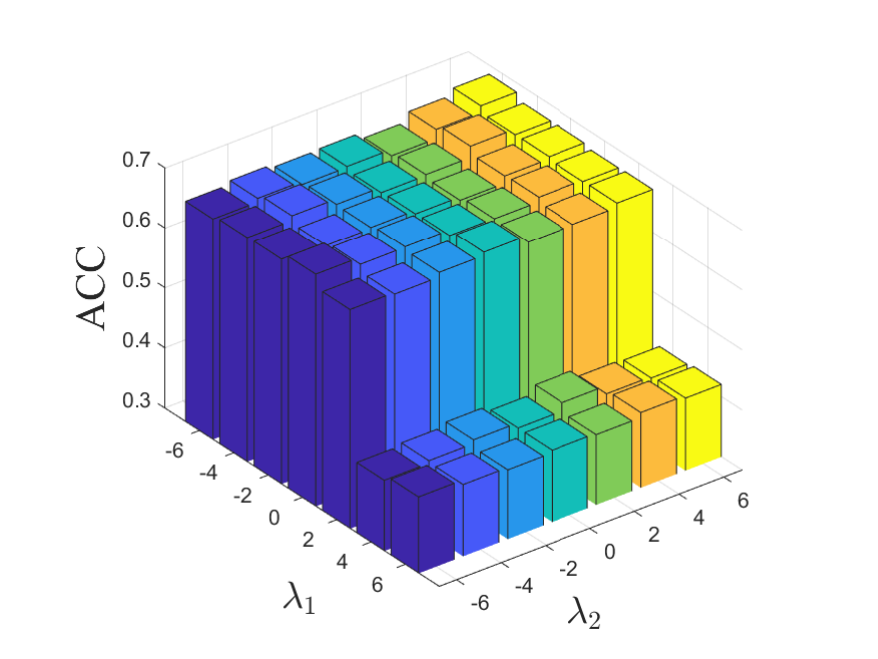}
}\hspace{-8mm}
\subfigcapskip=-3pt
\subfigure[LUNG (ACC)]{
    \label{d}
    \centering
    \includegraphics[width=4.5cm]{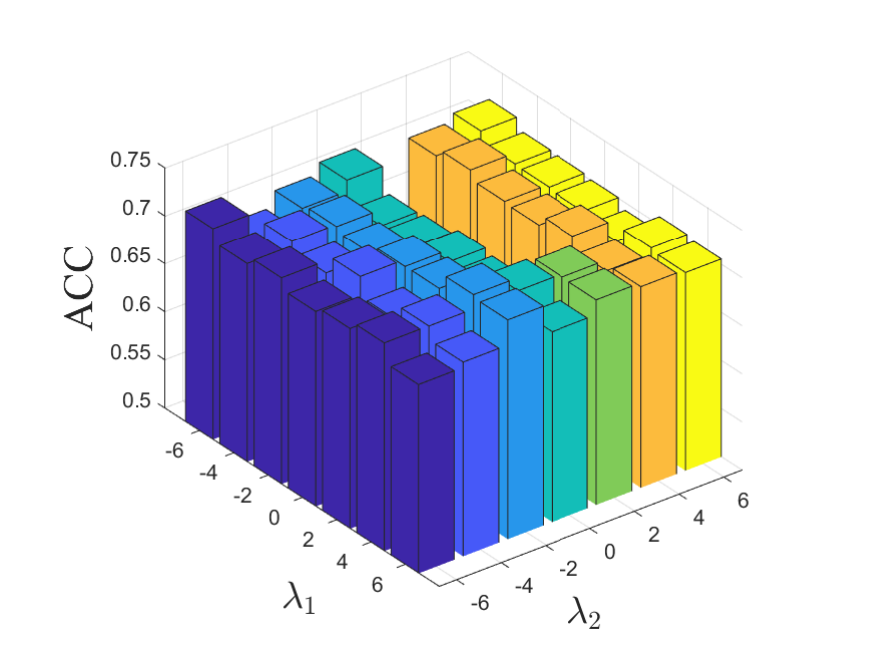}
}\hspace{-8mm}

\hspace{-1cm}
\subfigure[COIL20 (NMI)]{
    \label{e}
    \centering
    \includegraphics[width=4.5cm]{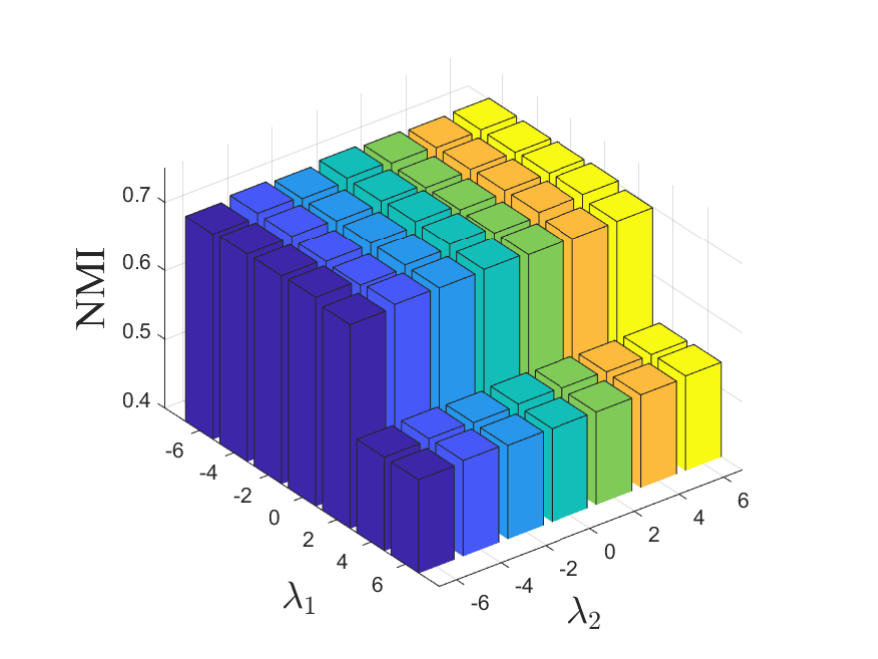}
}\hspace{-8mm}
\subfigcapskip=-3pt
\subfigure[Isolet (NMI)]{
    \label{f}
    \centering
    \includegraphics[width=4.5cm]{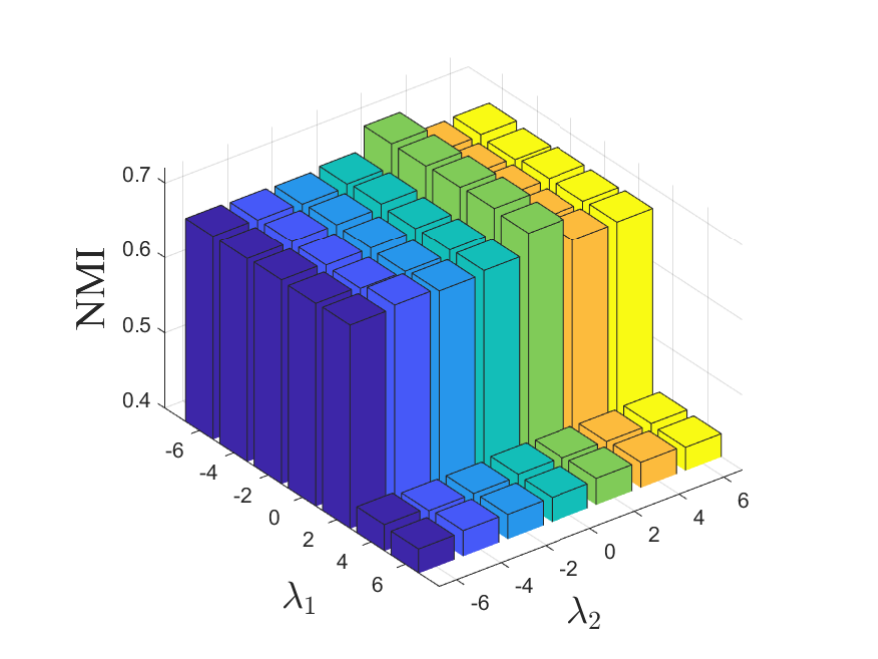}
}\hspace{-8mm}
\subfigcapskip=-3pt
\subfigure[USPS (NMI)]{
    \label{g}
    \centering
    \includegraphics[width=4.5cm]{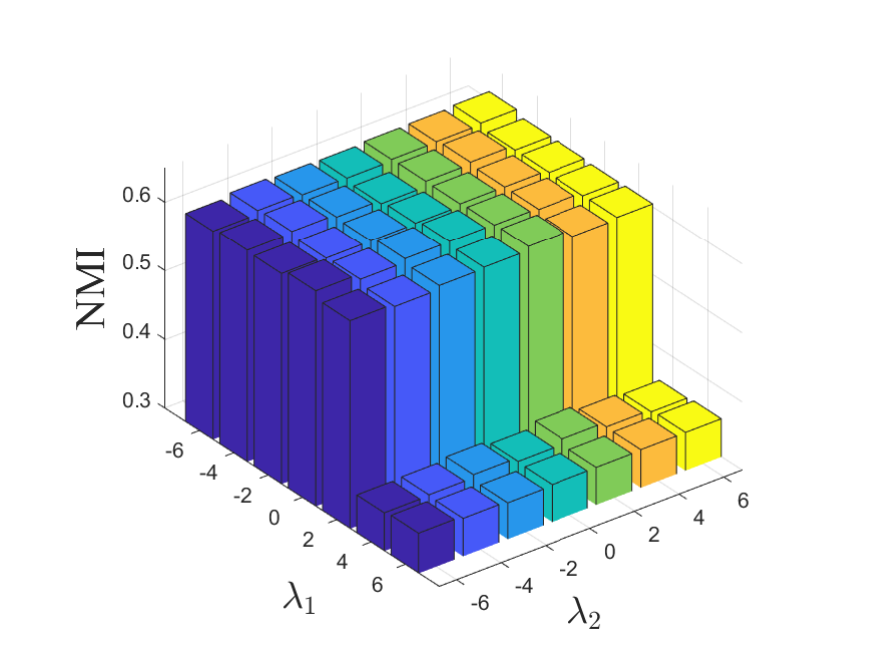}
}\hspace{-8mm}
\subfigcapskip=-3pt
\subfigure[LUNG (NMI)]{
    \label{h}
    \centering
    \includegraphics[width=4.5cm]{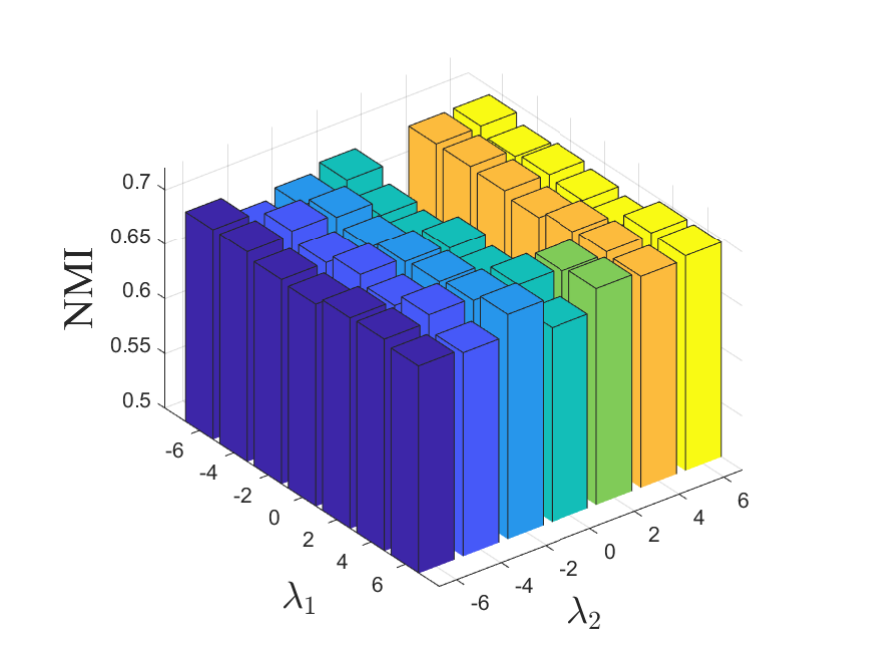}
}\hspace{-8mm}
\caption{Effects of $\lambda_1$ and $\lambda_2$ on four real-world datasets, where (a)-(d) are the ACC results and (e)-(h) are the NMI results.}
\label{plot-sen}
\end{figure*}

\begin{figure*}[t]
\centering
\hspace{-0.5cm}
\subfigcapskip=-1pt
\subfigure[COIL20]{
    \label{a}
    \centering
    \includegraphics[width=3.8cm]{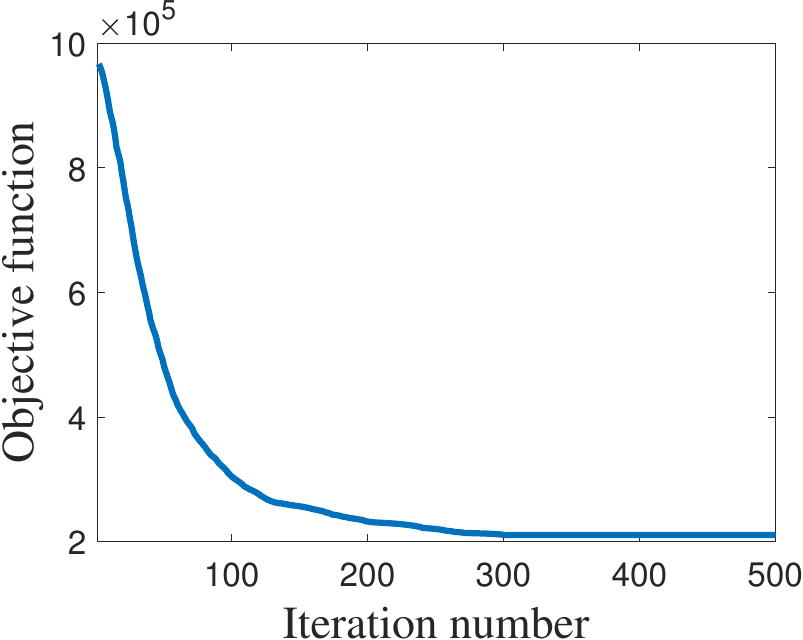}
}\hspace{-0mm}
\subfigcapskip=-1pt
\subfigure[Isolet]{
    \label{b}
    \centering
    \includegraphics[width=3.8cm]{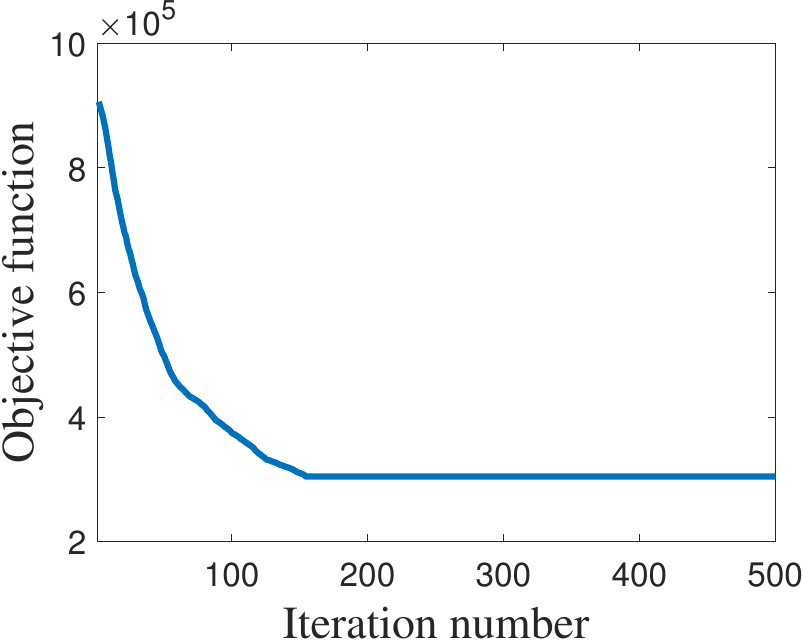}
}\hspace{-0mm}
\subfigcapskip=-1pt
\subfigure[USPS]{
    \label{c}
    \centering
    \includegraphics[width=3.8cm]{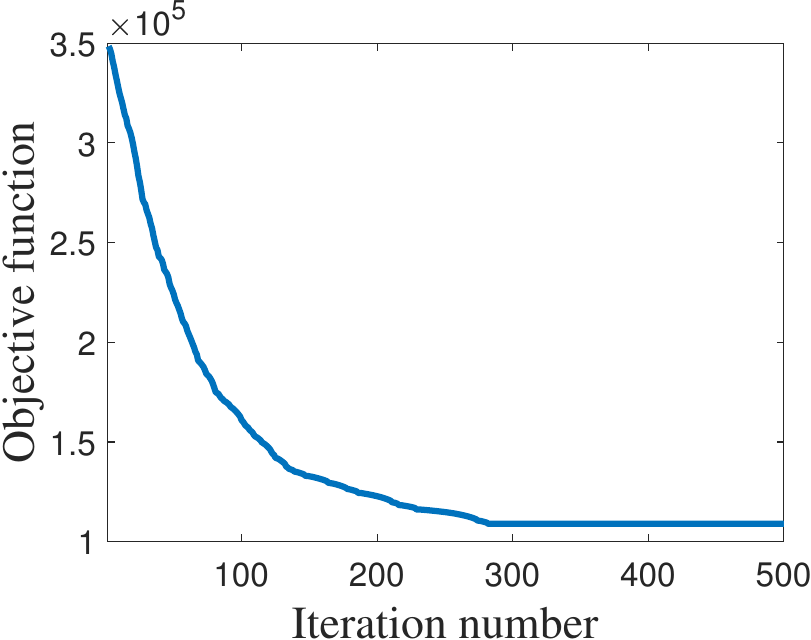}
}\hspace{-0mm}
\subfigcapskip=-1pt
\subfigure[LUNG]{
    \label{d}
    \centering
    \includegraphics[width=3.8cm]{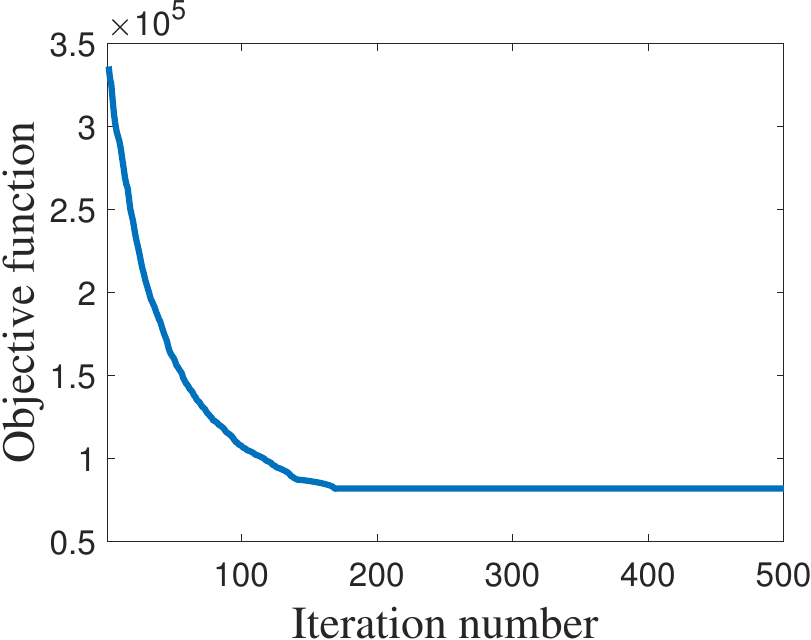}
}\hspace{-0mm}
\caption{Convergence curve of BSUFS on four real-world datasets.}
\label{con}
\end{figure*}

\subsection{Statistical Analysis}\label{stat}

To evaluate the pairwise differences between all compared methods, we employ the post-hoc Nemenyi test with the critical difference  value as a metric. The test outcomes for ACC and NMI are shown in Fig. \ref{hacc} and Fig. \ref{hnmi}, respectively.

It can be found that BSUFS demonstrates statistically significant differences when compared to SOGFS, LapScore, UDFS, and RNE, while no significant differences are observed in the performance of BSUFS relative to FSPCA, SPCAFS, SPCA-PSD, and PEN-PCAFS. This result shows that PCA-based methods are significantly different from graph-based methods, while our proposed BSUFS has no significant difference from other PCA-based methods.

\subsection{Effects of $p$ and $q$} \label{parameter}

There are three common choices of $p$ and $q$  in BSUFS, that are $0$, $1/2$, and $2/3$. In order to evaluate the importance of $p$ and $q$,  Fig. \ref{pq-acc} and Fig. \ref{pq-nmi} present the clustering performance in terms of ACC and NMI under these different values of $p$ and $q$, respectively. In these two figures, the x-axis represents  various values of $p$ and the color variations in the bars indicate different values of $q$. 

According to the clustering results, the following conclusions can be made. 
\begin{itemize}
    \item First, the optimal choices of $p$ and $q$ are different for different datasets. In specific, for the Isolet dataset, the optimal values of $p$ and $q$ are $0$ and $1/2$, respectively, while for the LUNG dataset, the optimal values are $1/2$ and $2/3$, respectively. 
    \item Second, under different values of $p$ and $q$, their ACC values and corresponding NMI values are not the same. For example, for the GLIOMA dataset, the ACC value varies greatly at $p=0$, which also shows that the selection of $q$ also affects the results.
    \item Finally, for the umist and MSTAR datasets, the best clustering performance can be observed when both $p$ and $q$ are set to 0, which illustrates that the extension from $(0,1)$ to $[0,1)$ is of importance.
\end{itemize}
Obviously, $p$ and $q$ should be tuned carefully. In practice, it is recommended to determine $p$ first and then $q$.

\subsection{Discussion} \label{discussion}

This section first visualizes the feature correlations between SPCAFS and our proposed BSUFS, then presents the model stability of all compared methods, and finally discusses the parameter sensitivity and convergence in the numerical perspective.

\subsubsection{Feature Correlation}
Fig. \ref{correlation} shows the feature correlation results using SPCAFS and BSUFS on the COIL20 and USPS datasets. Here, 10 features are selected, denoted as $F_1, F_2, \cdots,F_{10}$, and the correlation among these features is examined. It can be concluded that compared with SPCAFS, the features extracted by our proposed BSUFS are more discriminative. This means that the newly introduced $\ell_{q}$-norm can eliminate redundant features and improve feature selection results.

\subsubsection{Model Stability}

In this experiment, box plots of the 50 clustering results are shown in Fig. \ref{box}. It is obvious that in terms of ACC and NMI, the average values of BSUFS are generally larger than other methods. Especially for the Isolet data, the improvement is more obvious. This is because this dataset is highly sparse, and the advantages of double sparsity are fully demonstrated.

\subsubsection{Parameter Sensitivity}
Because there are two regularization terms, BSUFS needs to tune two regularization parameters. Fig. \ref{plot-sen} investigates the effects of two regularization parameters, i.e., $\lambda_1$ and $\lambda_2$, for ACC and NMI metrics.
Although the improvement under different $\lambda_2$ is not as significant as that under different $\lambda_1$, there are still differences, which also proves the necessity of the bi-sparse term in BSUFS.
In general, $\ell_{2,p}$-norm plays a vital role in BSUFS, while $\ell_{q}$-norm is a complementary choice in feature selection.

\subsubsection{Convergence Analysis}
Fig. \ref{con} shows the objective value curves for Algorithm \ref{algorithm} on four real-world datasets. It is seen that the objective function of our proposed BSUFS shows a consistent pattern of continuous decrease and reaches stability within finite iteration numbers. Although the convergence theorem cannot be derived as in Remark \ref{remark2}, the algorithm has good convergence in the numerical perspective.

\section{Conclusion}\label{conclusion}

In this study, a novel method called BSUFS is introduced for UFS, which integrates both $\ell_{2,p}$-norm and $\ell_{q}$-norm to PCA with $p,q\in[0,1)$. Technically, $\ell_{2,p}$-norm facilitates feature selection and  $\ell_{q}$-norm serves to remove the impact of redundant features. To our best knowledge, this is the first UFS method in a unified bi-sparse optimization framework. In algorithms, an efficient PAM optimization scheme is designed and its computational complexity is also analyzed. The feature selection results of BSUFS on synthetic and real-world datasets, respected to ACC and NMI, are more excellent than other competitors. Numerical studies also illustrate that the wide range of $p$ and $q$ is essential, and the best selection of them needs to be determined by the dataset. Furthermore, $\ell_{2,p}$-norm performs a leading role and $\ell_{q}$-norm enhances feature selection. Obviously, this bi-sparse framework can be applied to other related image processing fields.

In the future, we are interested in extending the proposed method to tensor cases \cite{zheng2024sparse} for better clustering performance. Besides, developing efficient optimization algorithms based on elegant proximal operator results \cite{liao2024subspace} and applying neural network methodologies \cite{zhang2023physics} to automatically learn parameters are worth investigating.

\bibliographystyle{IEEEtran}
\bibliography{mybib}

\end{document}